\documentclass[11pt]{amsart}
\usepackage{amssymb,amsmath,epsfig,graphics}
\theoremstyle{plain}
\newtheorem{thrm}{Theorem}[section]
\newtheorem{lemma}[thrm]{Lemma}
\newtheorem{prop}[thrm]{Proposition}
\newtheorem{cor}[thrm]{Corollary}
\newtheorem{rmrk}[thrm]{Remark}
\newtheorem{dfn}[thrm]{Definition}

\numberwithin{equation}{section}
\numberwithin{figure}{section}

\begin{document}
\newcommand{\SL}{\mathcal L^{1,p}(\Om)}
\newcommand{\Lp}{L^p(\Omega)}
\newcommand{\CO}{C^\infty_0(\Omega)}
\newcommand{\Rn}{\mathbb R^n}
\newcommand{\Rm}{\mathbb R^m}
\newcommand{\R}{\mathbb R}
\newcommand{\Om}{\Omega}
\newcommand{\Hn}{\mathbb H^n}
\newcommand{\HH}{\mathbb H^1}
\newcommand{\eps}{\epsilon}
\newcommand{\BVX}{BV_H(\Omega)}
\newcommand{\IO}{\int_\Omega}
\newcommand{\bG}{\boldsymbol{G}}
\newcommand{\bg}{\mathfrak g}
\newcommand{\p}{\partial}
\newcommand{\Xnu}{\overset{\rightarrow}{ X_\nu}}
\newcommand{\nuX}{\boldsymbol{\nu}^H}
\newcommand{\Up}{\bN^H}
\newcommand{\n}{\boldsymbol \nu}
\newcommand{\sigmau}{\boldsymbol{\sigma}^u_H}
\newcommand{\di}{\nabla_{i}^{H,\mS}}
\newcommand{\one}{\nabla_{1}^{H,\mS}}
\newcommand{\two}{\nabla_{2}^{H,\mS}}
\newcommand{\del}{\nabla^{H,\mS}}
\newcommand{\delXY}{\nabla^{H,\mS}_X Y}
\newcommand{\nui}{\nu^H_i}
\newcommand{\nuj}{\nu^H_j}
\newcommand{\dej}{\nabla_{j}^{H\mS}}
\newcommand{\cx}{\boldsymbol{c}^{H,\mathcal S}}
\newcommand{\sx}{\sigma_H}
\newcommand{\lx}{\mathcal L_H}
\newcommand{\pb}{\overline p}
\newcommand{\qb}{\overline q}
\newcommand{\ob}{\overline \omega}
\newcommand{\sbar}{\overline s}
\newcommand{\nuu}{\boldsymbol \nu_{H,u}}
\newcommand{\nuv}{\boldsymbol \nu_{H,v}}
\newcommand{\Bl}{\Bigl|_{\lambda = 0}}
\newcommand{\Bs}{\Bigl|_{s = 0}}
\newcommand{\mS}{\mathcal S}
\newcommand{\delh}{\Delta_H}
\newcommand{\delinf}{\Delta_{H,\infty}}
\newcommand{\nabh}{\nabla^H}
\newcommand{\delp}{\Delta_{H,p}}
\newcommand{\mO}{\mathcal O}
\newcommand{\delhs}{\Delta_{H,\mS}}
\newcommand{\dla}{\Delta_{H,\mS^\la}}
\newcommand{\lhs}{\hat{\Delta}_{H,\mS}}
\newcommand{\bN}{\boldsymbol{N}}
\newcommand{\bnu}{\boldsymbol \nu}
\newcommand{\la}{\lambda}
\newcommand{\nup}{(\nuX)^\perp}
\newcommand{\nuup}{\boldsymbol \nu^\perp_{H,u}}
\newcommand{\nuvp}{\boldsymbol \nu^\perp_{H,v}}
\newcommand{\ep}{\epsilon}
\newcommand{\e}{\mathfrak e}
\newcommand{\E}{\mathfrak E}
\newcommand{\BH}{\Bigl|_{H_g}}
\newcommand{\om}{\omega}
\newcommand{\se}{\sqrt{\epsilon}}
\newcommand{\oX}{\overline X}
\newcommand{\oY}{\overline Y}
\newcommand{\ou}{\overline u}
\newcommand{\fv}{\mathcal V^{H}_I(\mS;\mathcal X)}
\newcommand{\sv}{\mathcal V^{H}_{II}(\mS;\mathcal X)}


\title[Sub-Riemannian calculus on hypersurfaces in Carnot groups]
{Sub-Riemannian calculus on hypersurfaces in Carnot groups}

\author{D. Danielli}
\address{Department of Mathematics\\Purdue University \\
West Lafayette, IN 47907}
\email[Donatella Danielli]{danielli@math.purdue.edu}
\thanks{First author supported in part by NSF grants DMS-0002801 and CAREER DMS-0239771}

\author{N. Garofalo}
\address{Department of Mathematics\\Purdue University \\
West Lafayette, IN 47907}
\email[Nicola Garofalo]{garofalo@math.purdue.edu}
\thanks{Second author supported in part by NSF Grant DMS-0300477}

\author{D. M. Nhieu}
\address{Department of Mathematics\\Georgetown University \\
Washington, DC 20057-1233}
\email[Duy-Minh Nhieu]{nhieu@math.georgetown.edu}

%
%
\keywords{Horizontal Levi-Civita connection. Horizontal second
fundamental form. $H$-mean curvature. Intrinsic integration by
parts. First and second variation of the horizontal perimeter}
\subjclass{}
\date{\today}

\maketitle

\baselineskip 13.5pt

\tableofcontents


\section{\textbf{Introduction}}
\vskip 0.2in

The purpose of the present paper is to develop a sub-Riemannian
calculus on smooth hypersurfaces in a class of nilpotent Lie groups
which possess a rich geometry. Such groups arise as tangent spaces
of Gromov-Hausdorff limits of Riemannian manifolds \cite{CFG},
\cite{Be}, \cite{Mon}, and since they can be traced back to the
foundational paper of Carath\'eodory \cite{Ca} on Carnot
termodynamics, they have been christened Carnot groups by Gromov,
see \cite{Gro1}, \cite{Gro2}. Our main motivation is, in a broad
sense, the regularity theory of hypersurfaces of constant mean
curvature in such settings, as well as the applications of the
relevant calculus to the study of the Bernstein problem. These
problems have recently received increasing attention from several
groups of mathematicians and there exists nowadays a wide
literature. The following is only a partial list of references
\cite{GN1}, \cite{B}, \cite{FSS2}, \cite{FSS3}, \cite{FSS4},
\cite{GP}, \cite{CHMY}, \cite{CH}, \cite{ASV}, \cite{BC},
\cite{RR1}, \cite{HP}, \cite{DGN3}, \cite{DGNP}, \cite{CHY},
\cite{RR2}, \cite{BSV}. For an extensive bibliography we refer the
reader to the recent monographs \cite{DGNmem}, \cite{CDPT} (see also
the forthcoming book \cite{G}), and to the papers \cite{DGN3},
\cite{DGNP}. Carnot groups play a pervasive role in analysis,
geometry, and in various branches of the applied sciences, ranging
from problems in optimal control and robotics, crystallography,
mathematical finance, and
 neurophysiology of the brain.
This latter aspect, in particular, has been recently brought to
light in some very interesting works of Petitot and Tondut
\cite{Pe1}, \cite{PT}, \cite{Pe2}, and of Citti and Sarti
\cite{CS1}, \cite{CS2}, see also \cite{CMS1}, \cite{CMS2}. These
latter works have shown that there exists a close link between the
way in which the brain chooses to complete the missing visual data
in the first layer of the cerebral cortex, $V_1$, and the minimal
surfaces in a specific sub-Riemannian space, the so-called
roto-translation group, arising in the mathematical modeling of the
visual cortex $V_1$.

To describe the content of this paper we recall that during the past
century the study of minimal surfaces has been one of the main
driving forces in mathematics. Such development was prompted by the
study of the problems of Plateau and Bernstein which has led, as a
by-product, to the development of the Geometric Measure Theory, see
\cite{Fe}, \cite{Mat1}, \cite{Mat2}. Minimal surfaces also play a
central role in the positive mass theorem from relativity due to
Schoen and Yau \cite{SY}, see also the lecture notes \cite{Sch}.
Given the substantial progress which has occurred during the past
decade in the theory of subelliptic equations, and in those closely
connected aspects of geometric measure theory in sub-Riemannian
 spaces, it seems natural
at this point to direct the attention to the understanding of those
tools which are necessary for the development of a rich theory of
minimal surfaces. As we mentioned above, in this paper we solely
discuss hypersurfaces. Minimal manifolds of higher codimension are
also of interest and we hope to investigate them in future studies.

In classical geometry a central notion is that of area of a (smooth)
hypersurface. Such notion was extended by De Giorgi \cite{DG1},
\cite{DG2}, with the introduction of his variational theory of
perimeters which allowed to assign an ``area" also to sets which are
not a priori smooth. In a Carnot group $\bG$ there exists a
corresponding variational notion of perimeter adapted to the
horizontal bundle $H\bG$ (for a brief introduction to Carnot groups
we refer the reader to section \ref{S:CG}). Given a distribution of
smooth left-invariant vector fields $X = \{X_1,...,X_m\}$ which is
an orthonormal basis of the horizontal bundle (and therefore it is
bracket-generating for $T\bG$), and an open set $\Om\subset \bG$, we
let
\[
\mathcal F(\Om)\ =\ \{\zeta = \sum_{i=1}^m \zeta_i X_i  \in
C^1_0(\Om,H\bG)\ \mid\ |\zeta|_\infty\ =\ \sup_{\Om}\ |\zeta|\ \leq
1\}\ .
\]

For a  function $u\in L^1_{loc}(\Om)$, the $H$-variation of $u$ with respect to $\Om$ is defined by
\[
Var_H(u;\Om)\ =\ \underset{\zeta\in \mathcal F(\Om)}{\sup}\ \int_{\bG} u\ \sum_{i=1}^m X_i \zeta_i\ dg\ .
\]

A function $u\in L^1(\Om)$ is called of bounded $H$-variation in
$\Om$ if $Var_H(u;\Om) <\infty$. The space $BV_H(\Om)$ of
functions with bounded $H$-variation in $\Om$, endowed with the
norm
\[
||u||_{BV_H(\Om)}\ =\ ||u||_{L^1(\Om)}\ +\ Var_H(u;\Om)\ ,
\]
is a Banach space. Similarly to the classical theory (for the
latter, see for instance \cite{Gi} and \cite{Z}), such space
constitutes the appropriate replacement of the horizontal Sobolev
$W^{1,1}_H(\Om)$ space in the study of the relevant minimal
surfaces, see \cite{GN1}. Let now $E\subset \bG$ be a measurable
set, $\Om\subset \bG$ be an open set. The $H$-perimeter of $E$ with
respect to $\Om$ is defined by the equation
\begin{equation}\label{Xper}
P_H(E;\Om)\ =\ Var_H(\chi_E;\Om)\ ,
\end{equation}
where $\chi_E$ denotes the indicator function of $E$, see
\cite{CDG1}. When $E$ possesses sufficient regularity, e.g. when
$\mS = \p E$ is a hypersurface of class $C^2$, then one finds that
\begin{equation}\label{permeasure0}
P_H(E;\Om)\ =\ \int_{\Om \cap \p E} d \sigma_H\ =\ \int_{\Om \cap \p
E} \frac{|\Up|}{|\bN|}\ dH_{N-1}\ ,
\end{equation}
where we have denoted with $\bN^H$ the projection of the (non-unit)
Riemannian normal to $\p E$ onto the subbundle $H\bG$. It is
intereting to note that, in this situation, a useful alternative
understanding of the $H$-perimeter \eqref{Xper} can be obtained by
blowing-up the (suitably normalized) standard surface measure
associated with the Riemannian regularization of the sub-Riemannian
metric, see Theorem \ref{T:rr} below.

A ``minimal surface" in $\Om$ was defined in \cite{GN1} as the
boundary of a set of least $H$-perimeter, among all those with the
same boundaries outside $\Om$. The existence of such ``surfaces" (a
priori, these are just sets of locally finite $H$-perimeter), and a
measure theoretic solution of the Plateau problem, were also
established in \cite{GN1} following the classical approach of De
Giorgi \cite{DG1}, \cite{DG2}, \cite{DCP}. The natural question
arises of whether such measure theoretic minimal surfaces have, at
least when they are sufficiently smooth, vanishing ``mean
curvature". This prompts to investigate an appropriate notion of
mean curvature adapted to the horizontal bundle $H\bG$. For level
sets such a notion was proposed by one of us back in 1997, see
\cite{G1}, but its geometric content was not obvious. For the
Heisenberg group $\HH$, another notion of mean curvature was
introduced by Pauls in \cite{Pa}, who studied the solvability of the
Plateau problem by means of the Riemannian regularization of the
sub-Riemannian metric. For a surface in a three-dimensional CR
manifold, yet another notion of mean curvature has been recently
proposed in \cite{CHMY}. For instance, if the ambient manifold is
the Heisenberg group $\HH$, then the mean curvature of a surface
$\mS\subset \HH$ is defined as the standard curvature of the curve
of intersection of $\mS$ with the horizontal plane passing through
the base point. We note that, for surfaces in a Carnot group, this
same notion of curvature was also already explicitly introduced in
\cite{DGN1}. In this paper, given a $C^2$ hypersurface $\mS$ in a
Carnot group $\bG$, we introduce a second fundamental form on $\mS$
adapted to the horizontal subbundle $H\bG$, and a geometric notion
of mean curvature of $\mS$, and we show that the latter coincides
with either one of those proposed in \cite{G1}, \cite{Pa}, and
\cite{CHMY}, see Propositions \ref{P:equalMC}, \ref{P:pauls} and
\ref{P:chmy}.

In a Carnot group $\bG$, with grading of the Lie algebra $\bg = V_1
\oplus...\oplus V_r$, we define a smooth left-invariant Riemannian
metric $<\cdot,\cdot>$ by imposing that the vector fields
$X_1,...,X_m, ..., X_{r,m_r}$, defined in \eqref{genbasis}, be
orthonormal, see section \ref{S:lcc}. We can thus consider the
Riemannian connection $\nabla$ on $\bG$ induced by $<\cdot,\cdot>$.
We define the horizontal Levi-Civita connection $\nabh$ on $\bG$ by
projecting $\nabla$ onto the horizontal bundle $H\bG$, see section
\ref{S:lcc}. We note explicitly that $\nabh$ is, in essence,
Cartan's non-holonomic connection introduced in his address at the
Bologna International Congress of Mathematicians in 1928, see
\cite{C}.

In section \ref{S:XMC}, given an oriented $C^2$ hypersurface
$\mathcal S \subset \bG$, with Riemannian normal $\bN$, we define
the horizontal normal $\bN^H$ to $\mS$ as the projection of $\bN$
onto the horizontal bundle, and the horizontal Gauss map as $\nuX=
\bN^H/|\bN^H|$. Note that $|\bN^H|\not= 0$ at every point which does
not belong to the characteristic set $\Sigma_\mS$ of $\mS$. We
recall that the latter is the collection of all points $g\in \mS$ at
which $H_g\bG\subset T_g\mS$. An important notion is that of
horizontal tangent bundle $HT\mS$ to $\mS$, whose fiber $HT_g\mS$ at
each point $g\in \mS\setminus \Sigma_\mS$ is defined as the
collection of all horizontal vectors which are orthogonal to
$\bN^H$. It can be easily recognized that $HT_g\mS = T_g\mS \cap
H_g\bG$. To obtain a connection on $HT\mS$ we then project the
horizontal Levi-Civita connection $\nabh$ on the horizontal tangent
bundle $HT\mS$. More explicitly, for every $X,Y\in C^1(\mS;HT\mS)$
we define
\[
\delXY\ =\ \nabla^H_{\oX} \oY\ -\ <\nabla^H_{\oX} \oY,\nuX> \nuX\ ,
\]
where $\oX, \oY$ are any two horizontal vector fields on $\bG$ such
that $\oX = X$, $\oY = Y$ on $\mS$ (note that the above definition
does not depend on the choice of the extensions). Unlike its
Riemannian counterpart, the connection $\delXY$ is not torsion free
in general, and therefore it is not Levi-Civita in general. This is
due to the fact that, given $X,Y\in C^1(\mS;HT\mS)$, the projection
$[X,Y]^H$ onto the horizontal bundle of $[X,Y]$ does not in general
belong to the horizontal tangent space to $\mS$, $HT\mS$. We note in
passing that an interesting situation in which $\delXY$ is
Levi-Civita is that when $\bG = \HH$, the first Heisenberg group, or
when $\bG = \mathfrak E$, the four dimensional Engel group, see
section \ref{S:E}.

Inspired by the Riemannian situation we next project $\nabh$ along
the horizontal Gauss map $\nuX$. In this way we are able to
introduce the following notion of horizontal second fundamental form
on $\mS$
\begin{equation}\label{sff0} II^{H,\mS}(X,Y)\ =\ <\nabla^H_X Y,\nuX>
\nuX\ ,
\end{equation}
where $X,Y\in C^1(\mS;HT\mS)$. Since $[X,Y]^H$ is not in general in
$HT\mS$, unlike its Riemannian predecessor \eqref{sff0} is not
symmetric. One has in fact,
\[
II^{H,\mS}(X,Y)\ -\ II^{H,\mS}(Y,X)\ =\ <[X,Y]^H,\nuX> \nuX\ \not=\
 0\ .
\]

At every point $g_0\not\in \Sigma_\mS$, we define the horizontal
mean curvature $\mathcal H$ (or $H$-mean curvature) of $\mS$ as the
negative of the trace of the (symmetrized) second fundamental form.
If $\{\boldsymbol e_1,...,\boldsymbol e_{m-1}\}$ is an othonormal
basis of $HT\mS$, we thus have
\begin{equation}\label{Hmc0}
\mathcal H\ =\ -\ \sum_{i=1}^{m-1} <\nabla^H_{\boldsymbol e_i}
\boldsymbol e_i,\nuX>\ .
\end{equation}

If instead $g_0\in \Sigma_\mS$, then we define $\mathcal H(g_0)$ as
the $\underset{g\to g_0, g\not\in\Sigma_\mS}{\lim} \mathcal H(g)$,
whenever such limit exists. A $C^2$ hypersurface $\mS\subset \bG$ is
said to have constant mean-curvature $c\in \R$ if $\mathcal H \equiv
c$ as a continuous function on $\mS$. We call $\mS$ $H$-minimal if
$\mathcal H \equiv 0$ on $\mS$. We mention that recently Hladky and
Pauls \cite{HP} have introduced a notion of mean curvature for
hypersurfaces in a class of sub-Riemannian spaces which encompasses
that of Carnot groups. Their interesting approach can be seen as a
generalization of the Webster-Tanaka geometric framework for CR
manifolds, and systematically exploits the Lagrangian framework of
Bryant, Griffiths and Grossmann \cite{BGG}. Although the notion of
second fundamental form in \cite{HP} is different from the one
introduced in this paper, we notice that, specialized to Carnot
groups, their notion of mean curvature coincides with \eqref{Hmc0}.

Having introduced the notion of $H$-mean curvature, and $H$-minimal
surface, following the steps of the classical developments, it is
natural to study questions of regularity, stability, etc. It is
well-known that in the classical setting when $\mS \subset \Rn$,
with the standard surface measure $d\sigma$, an essential role in
this program is played by the following integration by parts
formula, see e.g. \cite{Gi},
\begin{equation}\label{ibpR}
 \int_\mS \nabla f\ d\sigma\ =\ (n-1)\ \int_\mS f\ H\ \n\
d\sigma\ , \end{equation} where $\nabla$ denotes the Levi-Civita
connection on $\mS$, $f\in C^2_0(\mS)$, and $H$ is the mean
curvature of $\mS$. For instance, the fundamental a priori gradient
estimates for minimal surfaces are derived from \eqref{ibpR}, see
\cite{BDM}. In section \ref{S:IBP} we establish an appropriate
generalization of \eqref{ibpR} to the case of a hypersurface in a
Carnot group. The interesting feature of such intrinsic integration
by parts formula is that the role of the surface measure is played
by the $H$-perimeter. Furthermore, it links the horizontal
connection $\nabla^{H,\mS}$ on $\mS$ to the $H$-mean curvature of
$\mS$. The relevant results states that for every $f\in
C^1_0(\mS\setminus \Sigma_\mS)$,
\begin{equation}\label{ibpSR} \int_\mS \nabla^{H,\mS} f\ d\sigma_H\
=\  \int_\mS f\ \bigg\{\mathcal H\ \nuX\ -\ \cx\bigg\}\ d\sigma_H\ ,
\end{equation}
where $\cx$ is a vector field on $\mS \setminus \Sigma_\mS$ with
values in the horizontal tangent space $HT\mS$, see Theorem
\ref{T:ibp}. This result plays a central role in the establishment
of the fundamental first and second variation formulas for the
$H$-perimeter in sections \ref{S:1&2var} and \ref{S:stab}. Although
\eqref{ibpSR} formally resembles \eqref{ibpR}, and in fact it
encompasses its Riemannian predecessor, the presence of the vector
field $\cx$ represents a new aspect which reflects the lack of
torsion freeness of the connection $\nabla^{H,\mS}$, see also
Proposition \ref{P:covder} below. In the Abelian case when $\bG
\cong \R^n$, then $\cx \equiv 0$ and we recover \eqref{ibpR}.
Another interesting situation in which $\cx\equiv 0$ is when $\mS$
is a vertical cylinder on the horizontal layer, i.e., when $\mS$ is
locally described by a defining function which depends only on the
horizontal variables.

Using the connection $\nabla^{H,\mS}$ we define two differential
operators on $\mS$, see Definition \ref{D:lb}. The former, denoted
by $\Delta_{H,\mS}$, is a sub-Riemannian version of the classical
Laplace-Beltrami operator on a manifold. The latter, indicated by
$\lhs$, contains an additional drift term, and is motivated by the
intrinsic integration
 by parts formula \eqref{ibpSR}. Its main raison
d'\^{e}tre, in fact, is that a Stokes' type theorem holds for it,
see Corollary \ref{C:lapX}. Formula \eqref{ibpSR} implies the
following identity
\[
\int_\mathcal S <\del u,\del \zeta> d\sigma_H\ =\ -\ \int_\mathcal S
u\ \lhs \zeta\ d\sigma_H\ .
\]
for every $u\in C^1(\mS)$, and every $\zeta \in C^2_0(\mS\setminus
\Sigma_\mS)$. Using this identity, we introduce a notion of
sub-harmonicity on $\mS$, see Definitions \ref{D:DI}, \ref{D:sh}. It
is an interesting open question to study the properties of
non-negative sub-harmonic functions on $\mS$. For instance, when
$\mS$ is $H$-minimal, do such functions satisfy some kind of
sub-mean value formula?

In Theorem \ref{T:mcf} we connect the operator $\Delta_{H,\mS}$ to
the flow by horizontal mean curvature recently introduced by Bonk
and Capogna \cite{BC}. We show that, similarly to its Riemannian
counterpart, such flow satisfies the following interesting partial
differential equation involving the horizontal tangential Laplacian
$\Delta_{H,\mS^t}$ on the hypersurfaces $\mS^t = F(\mS,t)$, images
of $\mS$ through the flow $F(\cdot,t)$, see Theorem \ref{T:mcf},
\[
<\frac{\p F}{\p t},\bN>\ =\ < \Delta_{H,\mS^t} F,\bN>\ .
\]

Sections \ref{S:geomid}, \ref{S:1&2var} and \ref{S:stab} are
entirely devoted to a geometric study of $C^2$ surfaces in the
Heisenberg group $\HH$. In this setting, given a $C^2$ surface $\mS$
with horizontal Gauss map $\nuX$, one easily recognizes that $HT\mS$
is spanned by the single vector field $\nup$. The triple
$\{\nup,\nuX,T\}$ forms an orthonormal moving frame on $\mS$. In
section \ref{S:geomid} we establish various geometric identities
which connect horizontal covariant differentiation along such frame
to geometric quantities such as the $H$-mean curvature and its
derivatives.

In section \ref{S:1&2var} we use such identities, in combination
with some notable integration by parts formulas which follow from
Theorem \ref{T:ibp}, see Lemma \ref{L:ibpcombined}. This lemma plays
a crucial role in establishing the first and second variation
formulas for the $H$-perimeter measure which constitute the main
results of the section, see Theorems \ref{T:variations} and
\ref{T:2varfinal}. The former allows to give a positive answer to
the question raised above: is a $C^2$ $H$-minimal surface surface a
stationary point of the $H$-perimeter? In Theorem \ref{T:variations}
we show that for $\mS\subset \HH$, the first-variation of the
$H$-perimeter for a deformation of $\mS$ along a vector field
$\mathcal X \in C^2_0(\mS\setminus \Sigma_\mS, \HH)$ is given by
\begin{equation}\label{fvSR}
\mathcal V^H_I(\mS;\mathcal X)\ =\
 \int_{\mathcal S}
\mathcal H\ \frac{<\mathcal X,\n>}{<\nuX,\n>}\ d\sigma_H\ ,
\end{equation}
where $\n=\bN/|\bN|$ represents the Riemannian Gauss map on $\mS$.
In particular, $\mathcal S$ is stationary if and only if it is
$H$-minimal (see also the less intrinsic first variation formula in
Theorem \ref{T:fv} for deformations along the normal $\bN$ and valid
for hypersurfaces in an arbitrary Carnot group).

The central result of section \ref{S:1&2var} is Theorem
\ref{T:2varfinal}, which provides a second variation formula for the
$H$-perimeter of $\mS$. The proof of such formula is considerably
more complex than that of \eqref{fvSR}, and obtaining it has
required a substantial effort. Despite such effort we notice,
however, that Theorem \ref{T:2varfinal} is in practice not as useful
as one would hope since it contains several terms whose geometric
content is not transparent, and which are very difficult to handle.
For the applications of the second variation formula to the
fundamental question of stability it is crucial to be able to
extract the geometry from Theorem \ref{T:2varfinal}. In order to do
so one needs to eliminate in the integrals involved the various
products of covariant derivatives of the projections of the testing
vector field $\mathcal X$ along the moving frame $\{\nup,\nuX,T\}$.
In this endeavor one has to choose with extreme care the terms to
play one against the other, so to be able to exploit the delicate
cancelations deriving from the various Lagrangian quantities
involved. Section \ref{S:stab} is devoted to this goal. In Theorem
\ref{T:svgeometric} we have succeeded in deriving the following
geometric second variation formula
\begin{equation}\label{svSR}
\sv\ =\ \int_\mS \bigg\{|\del F|^2\ +\ (2\mathcal A - \ob^2)
F^2\bigg\} d\sigma_H\ ,
\end{equation}
where $\mS\subset \HH$ is an $H$-minimal surface, $\mathcal X$ is as
in \eqref{fvSR}, and we have set \[ F\ =\ \frac{<\mathcal
X,\n>}{<\nuX,\n>}\ .
\]

The reader should compare \eqref{svSR} with the second variation
formula on p.153 in \cite{BGG}. The coefficient $2\mathcal A -
\ob^2$ of $F^2$ in \eqref{svSR} is a geometric quantity which
involves the projection of $\bN$ along $T$, and its horizontal
covariant derivative along the vector field $\nup$. With
\eqref{svSR} in hands, one can attack the fundamental question of
the stability. A non-characteristic $H$-minimal surface $\mS$ is
called \emph{stable} if $\sv\geq 0$ for any $\mathcal X \in
C_0^2(\mS,\HH)$. In view of \eqref{svSR} we see that a surface $\mS$
is stable if and only if the following stability inequality holds on
$\mS$
\begin{equation}\label{si}
 \int_\mS (\ob^2 - 2\mathcal A) F^2 d\sigma_H\ \leq\ \int_\mS |\del
F|^2\  d\sigma_H\ .
\end{equation}

We emphasize that one can think of \eqref{si} as a Hardy type
inequality on $\mS$. One should compare \eqref{si} with its
Riemannian counterpart, see e.g. inequality (1.105) in \cite{CM} for
normal deformations.

We conclude this introduction by emphasizing that the study of the
stability is an important new aspect in the sub-Riemannian Bernstein
problem. To clarify this point we recall the well-known fact that in
the classical Bernstein problem, stability does not apparently play
any role. This is due to the fact that the area functional for a
graph $x_{n+1} = u(x)$, $x\in \Om \subset \Rn$,
\[
A(u)\ =\ \int_\Om \sqrt{1 +|\nabla u|^2}\ dx\ ,
\]
is convex. As a consequence, a critical point of $A(u)$ is also a
local minimizer, and therefore stable. By contrast, the
sub-Riemannian area functional, the $H$-perimeter \eqref{Xper}, is
not convex, see \cite{DGNP}, and the resulting Euler-Lagrange
equation is not elliptic, but degenerate hyperbolic(-elliptic).
Using the stability inequality \eqref{si}, it has been recently
shown in \cite{DGN3} that, contrarily to what was believed by
several experts, the entire $H$-minimal graph $x = y t$ in $\HH$,
which has empty characteristic locus, is in fact unstable. This
discovery has underscored the role of the stability in the
sub-Riemannian Bernstein problem and opened the way to the solution
of the latter. Subsequently, in fact, this result has been
generalized in \cite{DGNP}, where it has been proved the instability
of every graph in $\HH$ of the type $x = y G(t)$, with $y\in \R$,
$t\in I\subset \R$, with $G\in C^2(I)$, and such that $G'>0$ on some
subinterval $J\subset I$. On the other hand, it has also been shown
in \cite{DGNP} that every entire $H$-minimal graph in $\HH$, with
empty characteristic locus, and which is not itself a vertical plane
$a x + b y =\gamma$, after possibly a left-translation and a
rotation about the $t$-axis, contains a graphical strip of the type
$x = y G(t)$, with $G'>0$ on some subinterval $J\subset \R$.
Combining these two results, the authors have obtained a solution of
the following sub-Riemannian Bernstein problem: \emph{The only
stable $H$-minimal entire graphs in $\HH$, with empty characteristic
locus, are the vertical planes}. The ideas in \cite{DGN3} have also
been used in the recent paper \cite{BSV} to prove a similar
Bernstein type theorem for the entire intrinsic graphs.

\vskip 0.6in


\section{\textbf{Carnot groups}}\label{S:CG}

\vskip 0.2in

In this section we collect some of the basic geometric facts about
Carnot groups. We particularly emphasize those properties which are
useful in this paper. For more extensive sources we refer the reader
to \cite{St1}, \cite{F2}, \cite{RS}, \cite{Str},
\cite{E1}-\cite{E3}, \cite{VSC}, \cite{Be}, \cite{Gro1},
\cite{Gro2}, \cite{Mon}, \cite{G}. A \emph{sub-Riemannian space} is
a triple $(M,HM,d)$ constituted by a connected Riemannian manifold
$M$, with Riemannian distance $d_{\mathcal R}$, a subbundle of the
tangent bundle $HM \subset TM$, and the Carnot-Carath\'eodory ($CC$)
distance $d$ generated by $HM$. Such distance is defined by
minimizing only on those absolutely continuous paths $\gamma$ whose
tangent vector $\gamma'(t)$ belongs to $H_{\gamma(t)}M$, see
\cite{NSW}, \cite{Be}. Riemannian manifolds are a special example of
sub-Riemannian spaces. They correspond to the case $HM = TM$. The
tangent space of a sub-Riemannian space is itself a sub-Riemannian
space (or a quotient of such spaces), but of a special type. It is a
graded Lie group whose Lie algebra is nilpotent. These groups, which
owe their name to the foundational paper of Charath\'eodory
\cite{Ca} on Carnot thermodynamics, occupy a central position in the
study of hypoelliptic partial differential equations, harmonic
analysis, sub-Riemannian geometry, CR geometric function theory, but
also in the applied sciences such as mathematical finance,
neurophysiology of the brain, mechanical engineering. They are
called Carnot groups.

A Carnot group of step $r$ is a connected, simply connected Lie
group $\bG$ whose Lie algebra $\bg$ admits a stratification $\bg=
V_1 \oplus \cdots \oplus V_r$ which is $r$-nilpotent, i.e.,
$[V_1,V_j] = V_{j+1},$ $j = 1,...,r-1$, $[V_j,V_r] = \{0\}$, $j =
1,..., r$. We assume henceforth that $\mathfrak g$ is endowed with a
scalar product $<\cdot,\cdot>_\bg$ with respect to which the $V_j's$
are mutually orthogonal. A trivial example of (an Abelian) Carnot
group is $\bG = \Rn$, whose Lie algebra admits the trivial
stratification $\bg = V_1 = \Rn$. The simplest non-Abelian example
of a Carnot group of step $r=2$ is the $(2n+1)$-dimensional
Heisenberg group $\Hn$, which is described in section \ref{S:E}.
Given a Carnot group $\bG$, by the above assumptions on the Lie
algebra one immediately sees that any basis of the \emph{horizontal
layer} $V_1$ generates the whole $\bg$. We will respectively denote
by
\begin{equation}\label{LT}
L_g(g')\ =\ g \ g'\ ,\quad\quad\quad\quad R_g(g')\ =\ g' \ g\ ,
\end{equation}
the operators of left- and right-translation by an element $g\in
\bG$.

The exponential mapping $\exp : \bg \to \bG$ defines an analytic
diffeomorphism onto $\bG$. We recall the important
Baker-Campbell-Hausdorff formula, see, e.g., sec.2.15 in \cite{V},
\begin{equation}\label{BCH}
exp(\xi)\  exp(\eta)\ =\ exp{\bigg(\xi + \eta + \frac{1}{2}\
[\xi,\eta] + \frac{1}{12} \big\{[\xi,[\xi,\eta]] -
[\eta,[\xi,\eta]]\big\} + ...\bigg)}\ ,
\end{equation}
where the dots indicate commutators of order four and higher. Each
element of the layer $V_j$ is assigned the formal degree $j$.
Accordingly, one defines dilations on $\bg$ by the rule
\[
\Delta_\lambda \xi\ =\ \lambda\ \xi_1\ +\ ...\ +\ \lambda^r\
\xi_r\ ,
\]
provided that $\xi = \xi_1 + ... + \xi_r \in \bg$, with $\xi_j\in
V_j$. Using the exponential mapping $\exp : \bg \to \bG$, these
anisotropic dilations are then tansferred to the group $\bG$ as
follows
\[
\delta_\lambda(g)\ =\ \exp\ \circ\ \Delta_\lambda\ \circ
\exp^{-1}\ g\ .
\]

Throughout the paper we will indicate by $dg$ the bi-invariant
Haar measure on $\bG$ obtained by lifting via the exponential map
$exp$ the Lebesgue measure on $\bg$. We let $m_j = dim\ V_j$, $j=
1,...,r$, and denote by $N = m_1 + ... + m_r$ the topological
dimension of $\bG$. One easily checks that
\[
(d\circ\delta_\lambda)(g)\ =\ \lambda^Q\ dg ,
\quad\quad\text{where}\quad Q\ =\ \sum_{j=1}^r j\, m_j\ .
\]

The number $Q$, called the \emph{homogeneous dimension} of $\bG$,
plays an important role in the analysis of Carnot groups. In the
non-Abelian case $r>1$, one clearly has $Q>N$.

We denote by $d(g,g')$ the \emph{CC distance} on $\bG$ associated
with the system $X$. It is well-known that $d(g,g')$ is equivalent
to the \emph{gauge pseudo-metric} $\rho(g,g')$ on $\bG$, i.e.,
there exists a constant $C=C(\bG)>0$ such that
\begin{equation}\label{equivalence}
C\ \rho(g,g')\ \leq\ d(g,g')\ \leq\ C^{-1}\ \rho(g,g'),
\quad\quad\quad\quad g,g'\in \bG,
\end{equation}
see \cite{Ca}, \cite{Chow}, \cite{Ra}, \cite{NSW}, \cite{VSC}. The
pseudo-distance $\rho(g,g')$ is defined as follows, see \cite{F2}.
Let $|\cdot |$ denote the Euclidean distance to the origin  on
$\mathfrak g$. For $\xi = \xi_1 + \cdots + \xi_r \in \mathfrak g$,
$\xi_j\in V_j$, one lets
\begin{equation}\label{gauge}
|\xi|_{\mathfrak g}\ =\ \left(\sum_{j=1}^r
|\xi_j|^{2r!/j}\right)^{1/2r!}\ ,\quad\quad\quad |g|_{\bG}\ =\
|\exp^{-1}\ g|_{\mathfrak g}, \quad\quad g\in \bG ,
\end{equation}
and defines
\begin{equation}\label{pseudod}
\rho(g, g')\ =\ |g^{-1}\ g'|_{\bG}.
\end{equation}

Both $d$ and $\rho$ are invariant under left-translations
\begin{equation}\label{isometry}
d(L_g(g'),L_g(g''))\ =\ d(g',g'')\ ,\quad\quad\quad\quad
\rho(L_g(g'),L_g(g''))\ =\ \rho(g',g'')\ .
\end{equation}
and homogeneous of degree one
\begin{equation}\label{scaleinv}
d(\delta_\lambda(g'),\delta_\lambda(g''))\ =\ \lambda\ d(g',g'')\
,\quad\quad\quad\quad
\rho(\delta_\lambda(g'),\delta_\lambda(g''))\ =\ \lambda\
\rho(g',g'')\ .
\end{equation}

Denoting respectively with
\begin{equation}\label{balls}
B(g,R)\ =\ \{g'\in \bG\mid d(g',g)<R\}, \quad\quad\quad
B_\rho(g,R) = \{g'\in \bG\mid \rho(g',g)<R\} ,
\end{equation}
the $CC$ ball and the gauge pseudo-ball centered at $g$ with radius
$R$, one easily recognizes that there exist $\omega =
\omega(\bG)>0$, and $\alpha = \alpha(\bG)>0$ such that
\begin{equation}\label{growth}
|B(g,R)|\ =\ \omega\ R^Q, \quad\quad\quad |B_\rho(g,R)|\ =\
\alpha\ R^Q ,\quad\quad\quad g\in \bG, R > 0 .
\end{equation}

Let $\pi_j : \bg \to V_j$ denote the projection onto the $j$-th
layer of $\bg$. Since the exponential map $exp : \bg \to \bG$ is a
global analytic diffeomorphism, we can define analytic maps
$\xi_j:\bG\to  V_j$, $j=1,...,r$, by letting $\xi_j = \pi_j \circ
\exp^{-1}$. As a rule, we will use letters $g, g', g'', g_0$ for
points in $\bG$, whereas we will reserve the letters $\xi, \xi',
\xi'', \xi_0, \eta$, for elements of the Lie algebra $\bg$.  The
notation $\{e_{j,1},...,e_{j,m_j}\}$, $j = 1,...,r,$ will indicate
a fixed orthonormal basis of the $j-th$ layer $V_j$. For $g\in
\bG$, the projection of the \emph{exponential coordinates} of $g$
onto the layer $V_j$, $j = 1,..., r$, are defined as follows
\begin{equation}\label{coord1}
x_{j,s}(g)\ =\ <\xi_j(g), e_{j,s}>_\bg , \quad\quad\quad s =
1,...,m_j .
\end{equation}

The vector $\xi_j(g) \in V_j$, $j = 1,...,r$, will be routinely
identified with the point
\[
(x_{j,1}(g),...,x_{j,m_j}(g))\ \in\ \mathbb R^{m_j}\ .
\]

Since Carnot groups of step $r=2$ often play a special role in
analysis and geometry, it will be convenient to have a simplified
notation for objects in the horizontal layer $V_1$, and in the first
vertical layer $V_2$. For simplicity, we set $m = m_1$, $k = m_2$,
and let
\begin{equation}\label{xandy}
\{e_1,\ ...\ ,\ e_m\}\ =\ \{e_{1,1},\ ...\ ,\ e_{1,m_1}\}\ ,\quad
\quad\quad\{\epsilon_1,...,\ep_k\}\ =\ \{e_{2,1},\ ...\ ,\
e_{2,m_1}\}\ .
\end{equation}

We indicate with
\begin{equation}\label{coord2}
x_i(g)\ =\ <\xi_1(g), e_i>_\bg ,\quad\quad i = 1,...,m\
,\quad\quad\quad t_s(g)\ =\ <\xi_2(g), \ep_s>_\bg ,\quad\quad s =
1,...,k\ ,
\end{equation}
the projections of the exponential coordinates of $g$ onto $V_1$
and $V_2$ respectively. Whenever convenient, we will identify
$g\in \bG$ with its exponential coordinates
\begin{equation}\label{ec}
 x(g)\ \overset{def}{=}\ (x_1(g),...,x_m(g), t_1(g),...,t_k(g),...,x_{r,1}(g),...,x_{r,m_r}(g))\ \in\ \mathbb R^N\ ,
\end{equation}
and we will ordinarily drop in the latter the dependence on $g$,
i.e., we will write $g = (x_1,...,x_{r,m_r})$.

For later purposes it will be useful to introduce the horizontal
group constants of $\bG$. By the grading assumption on the Lie
algebra, we have $[V_1,V_1] = V_2$. Therefore, if $e_i, e_j \in
\{e_1,...,e_m\}$, we let
\begin{equation}\label{gc}
b^s_{ij}\ \overset{def}{=}\ <[e_i,e_j],\ep_s>_\bg\ ,
\quad\quad\text{so that}\quad\quad [e_i,e_j]\ =\ \sum_{s=1}^k
b^s_{ij}\ \ep_s\ , \quad\quad\quad i , j = 1,...,m\ .
\end{equation}

Consider the orthonormal basis
$\{e_1,...,e_m,\ep_1,....,\ep_k,...,e_{r,1},...,e_{r,m_r}\}$ of
$\bg$. Using \eqref{LT} we define left-invariant vector fields on
$\bG$ by letting
\begin{equation}\label{genbasis}
X_{j,s}(g)\ =\ (L_g)_*(e_{j,s})\ ,\quad\quad\quad j=1,...,r, \quad
s = 1,...,m_j\ ,
\end{equation}
where $(L_g)_*$ indicates the differential of $L_g$. As in
\eqref{xandy} we use a special notation for the first two layers,
and let
\begin{equation}\label{X}
X_i(g)\ =\ (L_g)_*(e_i)\ ,\quad  i = 1,...,m, \quad\quad T_s(g)\ =\
(L_g)_*(\ep_s)\quad s = 1,...,k,\quad\ g\in \bG\ .
\end{equation}

Using the Baker-Campbell-Hausdorff formula \eqref{BCH} we can
express \eqref{X} using the exponential coordinates \eqref{ec},
obtaining the following lemma.

\medskip

\begin{lemma}\label{L:zero}
For each $i = 1,...,m$, and $g=(x_1,...,x_{r,m_r})$, we have
\begin{align}\label{fdG}
X_i\ & =\ \frac{\partial }{\partial{x_i}}\ +\
\sum_{j=2}^{r}\sum_{s=1}^{m_j}\
b^s_{j,i}(x_1,...,x_{{j-1},m_{(j-1)}})\ \frac{\partial
}{\partial{x_{j,s}}}
\\
& =\ \frac{\partial }{\partial{x_i}}\ +\
\sum_{j=2}^{r}\sum_{s=1}^{m_j}\ b^s_{j,i}(\xi_1,...,\xi_{j-1})\
\frac{\partial }{\partial{x_{j,s}}}\ , \notag
\end{align}
where each $b^s_{j,i}$ is a homogeneous polynomial of
\emph{weighted degree} $j-1$. In particular, if $\bG$ has step
$r=2$, then for every $i = 1,...,m$, one has
\begin{align}\label{firstd}
X_i\ &=\ \frac{\partial }{\partial x_i}\ +\ \frac{1}{2}\ \sum_{s =
1}^k <[\xi_1,e_i], \epsilon_s>_\bg \frac{\partial}{\partial t_s}
\\
\ &=\ \frac{\partial }{\partial x_i}\ +\ \frac{1}{2}\ \sum_{s = 1}^k
\sum_{\ell = 1}^m b^s_{\ell i}\ x_{\ell}\ \frac{\partial}{\partial
t_s}\ , \notag
\end{align}
where $b^s_{\ell i}$ are the group constants defined by \eqref{gc}.
We notice that an immediate consequence of \eqref{fdG} is that
\begin{equation}\label{divfreeE}
div_E\ X_i\ =\ 0\ ,\quad\quad\quad\quad i = 1,...,m\ ,
\end{equation}
where $div_E\ X_i$ indicates the Euclidean divergence of $X_i$
with respect to the exponential coordinates.
\end{lemma}

\medskip

By weighted degree in the statement of Lemma \ref{L:zero} we mean
that, as previously mentioned, the layer $V_j$, $j=1,...,r,$ in the
stratification of $\bg$ is assigned the formal degree $j$.
Correspondingly, each homogeneous monomial $\xi_1^{\alpha_1}
\xi_2^{\alpha_2}...\xi_r^{\alpha_r}$, with multi-indices $\alpha_j =
(\alpha_{j,1},...,\alpha_{j,m_j}), j=1,...,r,$ is said to have
weighted degree $k$ if
\[
\sum_{j=1}^r j\ (\sum_{s=1}^{m_j} \alpha_{j,s})\ =\ k\ .
\]

\vskip 0.6in

\section{\textbf{Two basic models}}\label{S:E}

\vskip 0.2in

In this section we describe two basic models of Carnot groups. The
first example is the Heisenberg group $\Hn$ with step $r=2$. Such
group plays an ubiquitous role in analysis and geometry, see e.g.
\cite{St1}, \cite{FS}, \cite{Gav}, \cite{Ko1}-\cite{Ko3},
\cite{K1}-\cite{K3}, \cite{KaR}, \cite{KoR1}-\cite{KoR2}, \cite{Be},
\cite{Mon}, \cite{CDPT}. From the standpoint of geometry $\Hn$
constitutes the central prototype of a pseudoconvex CR  manifold,
with vanishing Webster-Tanaka curvature. In fact, via the Caley
transform it can be identified with the boundary of the Siegel upper
half-space \[ \mathcal D^+\ =\ \{(z,z_{n+1})\in \mathbb C^{n+1}\mid
Im\ z_{n+1}\
>\ 2 \sum_{j=1}^n |z_j|^2\}\ , \]
see Ch.12 in \cite{St2}.  The second example is the cyclic, or Engel
group $\mathfrak E$, of step $r=3$, see \cite{CGr}, \cite{Mon}. This
is an interesting example to keep in mind since it represents the
basic prototype of a group of step $r=3$, and thereby constitutes
the next level of difficulty with respect to the Heisenberg group.
Some fundamental analytical and geometric properties are true for
Carnot groups of step $r=2$, but fail for groups of step $r\geq 3$.
In this respect, $\mathfrak E$ is the simplest sub-Riemannian model
in which to test whether conjectures which are true in step two
continue to be valid in step three or higher.

\medskip

\noindent \textbf{The Heisenberg group $\Hn$.}  The underlying
manifold of this Lie group is simply $\mathbb R^{2n +1}$, with the
non-commutative group law
\begin{equation}\label{Hn}
g\ g'\ =\ (x,y,t)\ (x',y',t')\ =\ (x + x', y + y', t + t' +
\frac{1}{2} ( <x,y'> - <x',y> ))\ ,
\end{equation}
where we have let $x,x', y, y' \in \Rn$, $t,t'\in \mathbb R$. Let
$(L_g)_*$ be the differential of the left-translation \eqref{Hn}.
A simple computation shows that
\begin{align}\label{Hn2}
& (L_g)_*\left(\frac{\partial}{\partial x_i}\right)\
\overset{def}{=}\ X_i\ =\ \frac{\partial}{\partial x_i}\ -\
\frac{y_i}{2}\ \frac{\partial}{\partial t}\ , \quad\quad i =
1,...,n\ ,
\\
& (L_g)_*\left(\frac{\partial}{\partial y_i}\right)\
\overset{def}{=}\ X_{n+i}\ =\ \frac{\partial}{\partial y_i}\ +\
\frac{x_i}{2}\ \frac{\partial}{\partial t}\ , \quad\quad i =
1,...,n\ ,
\notag\\
& (L_g)_*\left(\frac{\partial}{\partial t}\right)\
\overset{def}{=}\ T\ =\ \frac{\partial}{\partial t} \notag
\end{align}

We note that the only non-trivial commutator is
\[
[X_i,X_{n+j}]\ =\ \delta_{ij}\ T\ , \quad\quad\quad i, j =
1,...,n\ ,
\]
therefore the vector fields $\{X_1,...,X_{2n}\}$ generate the Lie
algebra $\mathfrak h_n = \mathbb R^{2n+1} = V_1 \oplus V_2$, where
$V_1 = \mathbb R^{2n} \times \{0\}_t$, $V_2 = \{0\}_{(x,y)} \times
\mathbb R$. We notice that the sub-Laplacian (see \eqref{sl})
associated with the orthonormal basis $\{\frac{\partial}{\partial
x_1}, ... ,\frac{\partial}{\partial x_1}, \frac{\partial}{\partial
y_1}, ... , \frac{\partial}{\partial y_n} \}$ of $V_1$ is
\begin{equation}\label{SLHn}
\delh\ =\ \sum_{j=1}^{2n} X_j^2\ =\ \Delta_{x,y}\ +\
\frac{1}{4}(|x|^2 + |y|^2)\ \frac{\partial^2}{\partial t^2}\ -\
\frac{\partial}{\partial t}\ \sum_{j=1}^n\ \bigg\{y_j
\frac{\partial}{\partial x_j}\ -\ x_j \frac{\partial}{\partial
y_j}\bigg\}\ ,
\end{equation}
which coincides with the real part of the complex Kohn-Spencer
Laplacian, see \cite{St2}. The non-isotropic group dilations are
\begin{equation}\label{Hn3}
\delta_\lambda(g)\ =\ (\lambda x, \lambda y, \lambda^2 t)\ ,
\end{equation}
with homogeneous dimension $Q = 2n + 2$. A convenient
renormalization of the gauge \eqref{gauge} is given by
\begin{equation}\label{gaugeHn}
N(g)\ =\ \left((|x|^2 + |y|^2)^2\ +\ 16\ t^2\right)^{1/4}\ .
\end{equation}

The importance of such function is connected with the discovery
due to Folland \cite{F1} that the fundamental solution of
\eqref{SLHn} is given by
\begin{equation}\label{FSHn}
\Gamma(g)\ =\ \Gamma(g,e)\ =\ \frac{C_Q}{N(g)^{Q-2}}\ ,
\end{equation}
where $C_Q<0$ is an explicit constant.

As a useful illustration, we compute the metric tensor $g_{ij}d\xi_i
\otimes d\xi_j$ associated with the smooth Riemannian product on
$\HH$ with respect to which $\{X_1,X_2,T\}$ is an orthonormal basis.
From \eqref{Hn2} we obtain
\begin{equation}\label{sb}
\frac{\p}{\p x}\ =\ X_1\ +\ \frac{y}{2}\ T\ ,\ \frac{\p}{\p y}\ =\
X_2\ -\ \frac{x}{2}\ T\ ,\ \frac{\p}{\p t}\ =\ T\ ,
\end{equation}
and therefore the metric coefficients are given by $g_{11} =
<\frac{\p}{\p x},\frac{\p}{\p x}> = 1 + \frac{y^2}{4}$, $g_{12} =
<\frac{\p}{\p x},\frac{\p}{\p y}> = - \frac{xy}{4}$, etc. One
easily finds
\begin{equation}\label{gij}
(g_{ij})\ =\ \begin{pmatrix} 1 + \frac{y^2}{4} & - \frac{xy}{4} &
\frac{y}{2}
\\
- \frac{xy}{4} & 1 + \frac{x^2}{4} & - \frac{x}{2}
\\
\frac{y}{2} & - \frac{x}{2} & 1
\end{pmatrix}
\end{equation}

Notice that, since $det(g_{ij}) = 1$, the volume form is given by
the standard (Lebesgue) volume form $dx\wedge dy \wedge dt$ in
$\R^3$. The inverse $(g^{ij})$ of the matrix \eqref{gij} is given by
\begin{equation}\label{gijinverse}
(g^{ij})\ =\ \begin{pmatrix} 1  & 0 & - \frac{y}{2}
\\
0 & 1  &  \frac{x}{2}
\\
- \frac{y}{2} &  \frac{x}{2} & 1 + \frac{x^2 + y^2}{4}
\end{pmatrix}
\end{equation}

Recall now the expression of the Riemannian gradient in local
coordinates, see for instance \cite{He} p.387,
\begin{equation}\label{Rgrad}
\nabla u\ =\ \sum_{i,j=1}^N g^{ij}\ \frac{\p u}{\p \xi_i}\
\frac{\p }{\p \xi_j}\ ,
\end{equation}
where we have denoted by $N=dim(\bG)$. Keeping in mind
\eqref{Hn2}, a simple calculation gives
\[ (g^{ij})\ \begin{pmatrix} u_x \\
u_y \\
u_t \end{pmatrix}\ =\ \begin{pmatrix} X_1 u \\
X_2 u \\
\frac{x u_y - y u_x}{2} + \left(1 + \frac{x^2 + y^2}{4}\right) u_t
\end{pmatrix}\ =\ \begin{pmatrix} X_1 u \\
X_2 u \\
\frac{x}{2} X_2u - \frac{y}{2} X_1 u + Tu
\end{pmatrix}\
\]

From this formula, and from \eqref{sb}, \eqref{gijinverse}, we
finally obtain
\begin{equation}\label{RgradH}
\nabla u\ =\ <(g^{ij})\ \begin{pmatrix} u_x \\
u_y \\
u_t \end{pmatrix}, \begin{pmatrix} \frac{\p}{\p x} \\
\frac{\p}{\p y} \\
\frac{\p}{\p t} \end{pmatrix}>_{\R^3} \ =\  X_1u\ X_1\ +\ X_2 u\
X_2\ +\ Tu\ T\ ,
\end{equation}
which verifies \eqref{grad}. It is worth observing that the
Laplace-Beltrami operator is given by
\[
\Delta u\ =\ X_1X_1u + X_2X_2u + T Tu\ .
\]

\vskip 0.2in

\noindent \textbf{The four-dimensional Engel group.}\ We next
describe the four-dimensional cyclic or Engel group. This group is
important in many respects since it represents the next level of
difficulty with respect to the Heisenberg group and provides an
ideal framework for testing whether results which are true in step
$2$ generalize to step $3$ or higher. The reader unfamiliar with
the cyclic group can consult \cite{CGr}, or also \cite{Mon}. The
Engel group $\mathfrak E=K_3$, see ex. 1.1.3 in \cite{CGr}, is the
Lie group whose underlying manifold can be identified with
$\mathbb R^4$, and whose Lie algebra is given by the grading,
\[
\mathfrak e\ =\ V_1 \oplus V_2 \oplus V_3\ ,
\]
where $V_1 = span\{e_1,e_2\}$, $V_2 = span\{e_3\}$, and $V_3 =
span\{e_4\}$, so that $m_1=2$ and $m_2=m_3=1$. We assign the
bracket relations
\begin{equation}\label{comm}
[e_1,e_2]\ =\ e_3\ \quad\quad \quad [e_1,e_3]\ =\ e_4\ ,
\end{equation}
all other brackets being assumed trivial. For the corresponding
left-invariant vector fields on $\mathfrak E$ given by $X_i(g) =
(L_g)_*(e_i)$, $i=1,2$, $T(g) = (L_g)_*(e_3)$, $S(g) =
(L_g)_*(e_4)$, we obtain the corresponding commutator relations
\begin{equation}\label{commengel}
[X_1,X_2]\ =\ T\ \quad\quad \quad [X_1,T]\ =\ [X_1,[X_1,X_2]]\ =\ S\
,
\end{equation}
all other commutators being trivial. We observe that the
homogeneous dimension of $\mathfrak E$ is
\[
Q\ =\ m_1\ +\ 2\ m_2\ +\ 3\ m_3\ =\ 7\ .
\]

We will denote with $(x,y)$, $t$ and $s$ respectively the
variables in $V_1$, $V_2$ and $V_3$, so that any $\xi \in
\mathfrak e$ can be written as $\xi = x e_1 + y e_2 + t e_3 + s
e_4$. If $g = \exp(\xi)$, we will identify $g = (x,y,t,s)$. The
group law in $\mathfrak E$ is given by the
Baker-Campbell-Hausdorff formula \eqref{BCH}. In exponential
coordinates, if $g = \exp(\xi)$, $g' = \exp(\xi')$, we have
\[
g\ \circ\ g'\ =\ \xi\ +\ \xi'\ + \frac{1}{2}\ [\xi,\xi']\ +\
\frac{1}{12}\ \big\{[\xi,[\xi,\xi']]\ -\ [\xi',[\xi,\xi']]\big\} .
\]

A computation based on \eqref{comm} gives (see also ex. 1.2.5 in
\cite{CGr})
\[
g\ \circ\ g'\ =\ \bigg(x + x', y + y', t + t' + P_3, s + s' +
P_4\bigg)\ ,
\]
where
\[
P_3\ =\ \frac{1}{2}\ (x y' - y x')\ ,
\]
\[
P_4\ =\ \frac{1}{2}(x t' - t x') + \frac{1}{12} \bigg(x^2 y' - x
x' (y + y') + y x^{'2}\bigg)\ .
\]

Using the Baker-Campbell-Hausdorff formula we find the following
expressions for the vector fields $X_1,...,X_4$
\begin{equation}\label{vfE}
\begin{cases}
X_1\ =\ \frac{\p}{\p x}\ -\ \frac{y}{2}\ \frac{\p}{\p t}\ -\ \left(\frac{t}{2} + \frac{x y}{12} \right)\ \frac{\p}{\p s}\ ,  \\
X_2\ =\ \frac{\p}{\p y}\ +\ \frac{x}{2}\ \frac{\p}{\p t}\ +\ \frac{x^2}{12}\ \frac{\p}{\p s}\ ,  \\
T\ =\ \frac{\p}{\p t}\ +\ \frac{x}{2}\ \frac{\p}{\p s}\ , \\
S\ =\ \frac{\p}{\p s}\ .
\end{cases}
\end{equation}

We note that the action of $X_1,X_2,T$ on a function on $\mathfrak
E$ which is independent of the variable $s$ reduces to the action of
the corresponding vector fields in $\HH$.

\vskip 0.6in


\section{\textbf{The subbundle of horizontal planes}}\label{S:HP}

\vskip 0.2in

Consider a Carnot group $\bG$, with Lie algebra $\bg = V_1 \oplus
\cdots\oplus V_r$, with an orthonormal basis $\{e_1,...,e_m\}$ of
the horizontal layer $V_1$, and corresponding system $X =
\{X_1,...,X_m\}$ of generators, where $X_i(g) = (L_g)_*(e_i)$, $g\in
\bG$. Henceforth, the fiber $H_g\bG$ of the horizontal bundle at a
point $g\in \bG$ will be denoted by $H_g$, so that $H\bG =
\cup_{g\in \bG} H_g$. We note explicitly that $H_g = g \exp(V_1)$,
where $\exp :\bg \to \bG$ denotes the exponential mapping. We will
call $H_g$ the \emph{horizontal plane} through $g$. For example,
when $\bG$ is the Heisenberg group $\Hn$, then a simple computation
shows that the horizontal plane through a point $g_0 =
(x_0,y_0,t_0)$ is given by the hyperplane
\begin{equation}\label{hp}
H_g\ =\ \left\{(x,y,t)\in \Hn \mid t = t_0 +
\frac{1}{2}\left(<x_0,y> - <y_0,x>\right)\right\}\ .
\end{equation}

More in general, we have the following result which is Proposition
4.3 in \cite{DGN1}.

\medskip

\begin{prop}\label{P:Hplane2}
Let $\bG$ be a Carnot group of step $2$, then for any given
$g_0\in \bG$ the horizontal plane passing through $g_0$ is the
collection of all points $g\in \bG$ whose exponential coordinates
$(x,t) = (x(g),t(g))$ verify the $k$ linear equations
\[
\Psi_s(g)\ =\ t_s(g)\ -\ t_s(g_0)\ -\ \frac{1}{2}\ \sum_{i, j =
1}^m b^s_{i j}\ x_i(g_0)\ x_{j}(g)\ =\ 0\ ,\quad\quad\quad s =
1,...,k\ ,
\]
where $b^s_{ij}$ represent the horizontal group constants defined
by \eqref{gc}.
\end{prop}

\medskip

Another interesting example is provided by the four-dimensional
Engel group $\mathfrak E$ described in the previous section.
Identifying $\E$ with $\R^4$, with coordinates $g=(x,y,t,s)$, given
a point $g_0 = (x_0,y_0,t_0,s_0)$ we have that $H_{g_0} =
span\{X_1(g_0),X_2(g_0)\}$. A simple computation based on
\eqref{vfE} shows that $H_{g_0}$ is described by the two equations
\begin{equation}\label{hpe}
\begin{cases}
\Psi_1(x,y,t,s)\ =\ t - t_0 + \frac{x y_0 - x_0 y}{2}\ =\ 0\ ,
\\
\Psi_2(x,y,t,s)\ =\ s - s_0 + \frac{x (6 t_0 + x_0 y_0) - x_0^2 y
- 6 x_0 t_0}{12}\ =\ 0\ .
\end{cases}
\end{equation}

From \eqref{fdG} in Lemma \ref{L:zero} we see that for a Carnot
group $\bG$ of step $r$, with $N = dim(\bG)$, the horizontal plane
$H_{g_0}$ is described by a system of $N-m$ linear equations for
the exponential variables, see \eqref{ec},
\begin{equation}\label{defeq}
\begin{cases}
\Psi_1(g)\ =\ t_1(g) - t_1(g_0) - B_1(g)\ =\ 0\ ,
\\
.........
\\
 \Psi_k(g)\ =\ t_k(g) - t_k(g_0) - B_k(g)\ =\ 0\ ,
\\
.........
\\
 \Psi_{r,1}(g)\ =\ x_{r,1}(g) - x_{r,1}(g_0) - B_{r,1}(g)\ =\ 0\ ,
\\
.........
\\
\Psi_{r,m_r}(g)\ =\ x_{r,m_r}(g) - x_{r,m_r}(g_0) - B_{r,m_r}(g)\
=\ 0\ ,
\end{cases}
\end{equation}
with $B_j(g_0) = 0$ for $j=1,...,k$, ... ,$B_{r,j}(g_0) = 0$, $j =
1,...,m_r$.

\medskip

\begin{dfn}\label{D:surface}
We say that $\mS \subset \bG$ is a $C^k$ hypersurface if $\mS$ is a
co-dimension one immersed manifold of class $C^k$. If, in addition,
$\mS$ is embedded, then we say say that it is an embedded
hypersurface.
\end{dfn}

\medskip

We note explicitly that, by the implicit function theorem, for every
$g_0\in \mS$ there exists an open set $\mathcal O\subset \bG$ and a
function $\phi\in C^k(\mathcal O)$ such that: (i) $|\nabla \phi(g)|
\not= 0$ for every $g\in \mathcal O$; (ii) $\mS\cap \mathcal O =
\{g\in \mathcal O\mid \phi(g) = 0\}$. When we will need to use this
local representation, we will always assume that $\mS$ is oriented
in such a way that for every $g_0\in \mS$ and $\phi$ as in (ii), one
has $\bN(g_0) = \nabla \phi(g_0)$, where $\boldsymbol N$ denotes the
non-unit Riemannian normal to $\mathcal S$. The following notion
plays a pervasive role in sub-Riemannian geometry, as well as in the
study of subelliptic equations.

\medskip

\begin{dfn}\label{D:char}
Given a $C^1$ hypersurface $\mS\subset \bG$, a point $g_0\in
\mathcal S$ is called \emph{characteristic} if one has $H_{g_0}
\subset T_{g_0} \mathcal S$. Notice that this is equivalent to
saying that
\begin{equation}\label{charperp}
X_j(g_0)\ \in\ T_{g_0} \mathcal S\ ,\quad\quad\quad\quad j =
1,...,m\ .
\end{equation}
The \emph{characteristic locus of $\mathcal S$}, $\Sigma_{\mathcal
S}$, is the collection of all characteristic points of $\mathcal S$.
\end{dfn}

\medskip

Although we will not use in this paper the following two results, we
recall them because of their interest. The first theorem is a
special case of a result due to Derridj \cite{De1}, \cite{De2}.

\medskip

\begin{thrm}\label{T:derridj}
Let $\mS$ be a $C^\infty$ hypersurface in a sub-Riemannian space
$M$ of dimension $N$, then denoting with $H^s$ the $s$-dimensional
Hausdorff measure
 constructed with the Riemannian distance one has
\[
H^{N-1}(\Sigma_\mS)\ =\ 0\ .
\]
\end{thrm}

\medskip

For Carnot groups one has the following sharper result first
proved in codimension one by Balogh for the Heisenberg group
\cite{B}, and subsequently extended to arbitrary Carnot groups and
codimension by Magnani \cite{Ma1}, \cite{Ma2}.

\medskip

\begin{thrm}\label{T:magnani}
Let $\bG$ be a Carnot group and denote by $\mathcal H^s$ the
$s$-dimensional Hausdorff measure constructed with the
Carnot-Carath\'eodory distance. For any $C^1$ manifold of
codimension $k$ one has
\[
\mathcal H^{Q-k}(\Sigma_\mS)\ =\ 0\ .
\]
In particular, the characteristic set of a $C^1$ hypersurface has
zero $\mathcal H^{Q-1}$-measure.
\end{thrm}

\medskip

Since for a $C^2$ hypersurface in a Carnot group $\bG$ it was
proved in \cite{DGN1}, see also \cite{DGN2}, that   the
$H$-perimeter
 measure $P_H(\Om;\cdot)$, introduced in section \ref{S:pm} below, is mutually absolutely continuous with
respect to the Hausdorff measure $\mathcal H^{Q-1}$, we conclude
from Theorem \ref{T:magnani} that for such domains the $H$-perimeter
measure of the characteristic set is zero, i.e.,
\begin{equation}\label{magnani}
\sigma_H(\Sigma_\mS)\ =\ 0\ .
\end{equation}

\medskip

\begin{prop}\label{P:submanifold}
Let $\bG$ be a Carnot group, and $\mS\subset \bG$ be a $C^k$
hypersurface. If $g_0\in \mS \setminus \Sigma_\mS$, denote by $\mS_0
= \mS \cap H_{g_0}$. There exists a sufficiently small open
neighborhood $\mathcal O$ of $g_0$ such that $\mS_0 \cap \mathcal O$
is a $C^k$ submersed manifold of $\bG$ of dimension $m-1$.
\end{prop}

\begin{proof}[\textbf{Proof}]
According to Definition \ref{D:surface}, there exists a
neighborhood $\mathcal O$ of $g_0$ such that $\mS\cap \mathcal O =
\{g\in \mathcal O\mid \phi(g) = 0\}$. Using the exponential
coordinates \eqref{ec}, we now introduce the $C^k$ function $F :
exp^{-1}(\mathcal O) \subset \R^N \to \R^{N-m +1}$ defined by
\[
F\ \overset{def}{=}\
(\phi,\Psi_1,...,\Psi_k,...,\Psi_{r,1},...,\Psi_{r,m_r})\ ,
\]
where the $N-m$ functions $\Psi_j$ are as in \eqref{defeq}.
Clearly, we have $F(g_0) = 0 \in \R^{N-m+1}$. Denoting by $J_F$
the Jacobian matrix of $F$, we now claim that the hypothesis
\[
g_0\ \in\ \mS\setminus \Sigma_\mS\ \quad\quad
\Longrightarrow\quad\quad rank\ J_F(g_0) = N - m + 1\ . \]

Taking the claim for granted, we see that the conclusion of
Proposition \ref{P:submanifold} immediately follows from the
implicit function theorem (of course, by possibly restricting the
neighborhood $\mathcal O$), since the latter guarantees that,
locally around $g_0$, the set $\mS_0$ is a submersed manifold of
class $C^k$ of dimension $N - (N-m+1) = m-1$.

We now prove the claim in two special situations, namely that of a
Carnot group of step $r=2$, and that of the Engel group $\E$,
leaving it to the interested reader to provide the (lengthy) details
for a general Carnot group. Suppose then that $\bG$ has step $r=2$.
Since $g_0\not\in \Sigma$, we know that $\nabla_H \phi(g_0) \not=
0$. Therefore, there exists $i\in \{1,...,m\}$ such that
$X_i\phi(g_0)\not= 0$. Without loss of generality, let us assume
that $X_m\phi(g_0) \not=0$. According to \eqref{firstd} we thus have
\begin{equation}\label{nca} \phi_{x_m}(g_0)\ +\ \frac{1}{2}\
\sum_{s = 1}^k \sum_{j = 1}^m b^s_{jm}\ x_{j}(g_0)\
\phi_{t_s}(g_0)\ \not=\ 0\ .
\end{equation}

The Jacobian matrix of $F = (\phi,\Psi_1,...,\Psi_k)$ at $g_0$ is
now given by
\[
J_F(g_0)\ =\ \begin{pmatrix} \phi_{x_1} & ... & \phi_{x_m} &
\phi_{t_1} & \phi_{t_2} & ... & \phi_{t_k}
\\
- \frac{1}{2} \sum_{i=1}^m b^1_{i1} x_i(g_0) & ... & - \frac{1}{2}
\sum_{i=1}^m b^1_{im} x_i(g_0) & 1 & 0 & ... & 0
\\
- \frac{1}{2} \sum_{i=1}^m b^2_{i1} x_i(g_0) & ... & - \frac{1}{2}
\sum_{i=1}^m b^2_{im} x_i(g_0) & 0 & 1 & ... & 0
\\
. & ... & . & . & . & ... & .
\\
. & ... & . & . & . & ... & .
\\
- \frac{1}{2} \sum_{i=1}^m b^k_{i1} x_i(g_0) & ... & - \frac{1}{2}
\sum_{i=1}^m b^k_{im} x_i(g_0) & 0 & 0 & ... & 1
\end{pmatrix}\ ,
\]
where all derivatives of $\phi$ are evaluated at $g_0$. We
consider the $(k+1)\times (k+1)$ minor
\[
\tilde{J}_F(g_0)\ =\ \begin{pmatrix} \phi_{x_m} & \phi_{t_1} &
\phi_{t_2} & ... & \phi_{t_k}
\\
- \frac{1}{2} \sum_{i=1}^m b^1_{im} x_i(g_0) & 1 & 0 & ... & 0
\\
- \frac{1}{2} \sum_{i=1}^m b^2_{im} x_i(g_0) & 0 & 1 & ... & 0
\\
. & . & . & ... & .
\\
 . & . & . & ... & .
\\
- \frac{1}{2} \sum_{i=1}^m b^k_{im} x_i(g_0) & 0 & 0 & ... & 1
\end{pmatrix}\ ,
\]
of the matrix $J_F(g_0)$. A careful examination of the special
structure of the matrix $\tilde{J}_F(g_0)$, and the cofactor
expansion of its determinant, allow to conclude that
\[
det\ \tilde{J}_F(g_0)\ =\ \phi_{x_m} + \frac{1}{2} \sum_{i=1}^m
b^1_{im} x_i(g_0) \phi_{t_1} + ... + \frac{1}{2} \sum_{i=1}^m
b^k_{im} x_i(g_0) \phi_{t_k}\ \not=\ 0\ ,
\]
where in the last equation we have used \eqref{nca}. This proves
that $rank\ J_F(g_0) = k + 1 = N - m + 1$, and therefore the claim
follows for groups of step $r=2$.

If instead $\mS \subset \bG = \E$ is a hypersurface in the Engel
group, with $g_0 = (x_0,y_0,t_0,s_0) \in \mS\setminus \Sigma$,
then we can assume for instance that we have at $g_0$
\begin{equation}\label{ncE}
X_2\phi(g_0)\ =\ \phi_y(g_0) + \frac{x_0}{2} \phi_t(g_0) +
\frac{x_0^2}{12} \phi_s(g_0)\ \not=\ 0\ . \end{equation}

We consider the function $F = (\phi,\Psi_1,\Psi_2) : \R^4 \to
\R^3$, where $\Psi_i$, $i=1,2$ are as in \eqref{hpe}. Its Jacobian
matrix is given by
\[
J_F(g_0)\ =\ \begin{pmatrix} \phi_x & \phi_y & \phi_t & \phi_s
\\
\frac{y_0}{2} & - \frac{x_0}{2} & 1 & 0
\\
\frac{6 t_0 + x_0 y_0}{12} & - \frac{x_0^2}{12} & 0 & 1
\end{pmatrix}
\]

One readily sees that the $3\times 3$ minor
\[
\tilde{J}_F(g_0)\ =\ \begin{pmatrix} \phi_y & \phi_t & \phi_s
\\
 - \frac{x_0}{2} & 1 & 0
\\
- \frac{x_0^2}{12} & 0 & 1
\end{pmatrix}
\]
has determinant given by $X_2\phi(g_0)$. From \eqref{ncE} we
conclude that $rank\ \tilde{J}_F(g_0) = 3 = N - m + 1$, and again
the claim follows.

\end{proof}

 \vskip 0.6in



\section{\textbf{Horizontal Levi-Civita connection}}\label{S:lcc}

\vskip 0.2in

Let $\bG$ be a Carnot group of step $r$. Henceforth in this paper we
will assume that $\bG$ is endowed with a left-invariant Riemannian
metric $<\boldsymbol u,\boldsymbol v> = g_{ij}u^i v^j$, where
$\boldsymbol u , \boldsymbol v\in T\bG$, with respect to which the
left-invariant vector fields defined in \eqref{genbasis}
\[
\{X_1,...,X_m,T_1,...,T_k,...,X_{r,1},...,X_{r,m_r}\} \] constitute
an orthonormal frame for $T\bG$. No other inner product will be used
on $T\bG$, thereby when we write $<\cdot,\cdot>$ there will be no
risk of confusion. We denote with by $\nabla$ the corresponding
Levi-Civita connection on $\bG$. Recall that $\nabla$ is torsion
free, \begin{equation}\label{torsionfree} \nabla_X Y - \nabla_Y X\
=\ [X,Y]\ , \end{equation}
 and that
it is metric preserving, i.e., $\nabla g = 0$ or, equivalently,
\begin{equation}\label{compcond}
X<Y,Z>\ =\ <\nabla_X Y,Z>\ +\ <Y,\nabla_X Z>\ .
\end{equation}

Permuting cyclically the roles of $X,Y,Z$ in \eqref{compcond}, one
obtains the basic Koszul identity, see e.g. (1.13) on p.28 of
\cite{Sa},
\begin{align}\label{Koszul}
2 <\nabla_X Y,Z>\ & =\ X<Y,Z>  + Y<X,Z> - Z<X,Y>
\\
& -\ <Y,[X,Z]> - <X,[Y,Z]> + <Z,[X,Y]>\ . \notag
\end{align}

Using \eqref{Koszul} it is easy to check that
\begin{equation}\label{delii}
\nabla_{X_i}X_i\ =\ 0\ , i=1,...,m, \ ... , \
\nabla_{X_{j,m_j}}X_{j,m_j}\ =\ 0\ , j = 1,...,r\ .
\end{equation}

In addition to \eqref{delii}, we can easily verify from
\eqref{Koszul} and the grading of the Lie algebra, that
\begin{equation}\label{deij}
<\nabla_{X_i} X_j, X_\ell> \ =\ 0\ ,\quad\quad\quad i , j , l =
1,...,m\ .
\end{equation}

The remaining covariant derivatives and the Christoffel symbols can
be determined from the group constants. For instance,  we have the
following proposition.

\medskip

\begin{prop}\label{P:covder}
Let $\bG$ be a Carnot group of step $r$, then
\begin{equation}\label{delijZ}
\nabla_{X_i} X_j\ =\ \frac{1}{2} \sum_{s=1}^k b^s_{ij} T_s\
,\quad\quad i,j=1,...,m\ . \end{equation}
\begin{equation}\label{delTsp}
 \nabla_{T_p} T_s\ =\ 0\ ,\quad\quad p,s
= 1,...,k\ .
\end{equation}
\begin{equation}\label{delXT}
\nabla_{X_i} T_s  = - \frac{1}{2} \sum_{j=1}^m b^s_{ij} X_j +
\frac{1}{2} \sum_{p=1}^{m_3} <[X_i,T_s],X_{3,p}> X_{3,p}\ ,\
i=1,...,m,\ s=1,...,k\ .
\end{equation}
In particular, when $\bG = \Hn$ one has for $i,j=1,...,n,$
\[
\nabla_{X_i} X_{n+j}\ =\ \frac{\delta_{ij}}{2} \ T\ ,\quad\quad
\nabla_{X_i} T\ =\ \nabla_{T} X_i\ =\ -\ \frac{1}{2}\ X_{n+i}\ ,\
\quad \nabla_{X_{n+i}} T\ =\ \nabla_{T} X_{n+i}\ =\  \frac{1}{2}\
X_{i}\ .
\]
\end{prop}

\begin{proof}[\textbf{Proof}]
Using \eqref{Koszul}, for any vector field \[ Z = \sum_{\ell=1}^m
a_\ell X_\ell + \sum_{s=1}^k b_s T_s + \sum_{h=3}^r \sum_{p=1}^{m_h}
c_p X_{h,p}\ , \] we obtain
\begin{align*}
<\nabla_{X_i} X_{j},Z>\ & =\ \sum_{\ell =1}^m a_\ell <\nabla_{X_i}
X_{j},X_\ell> + \sum_{s=1}^k b_s <\nabla_{X_i} X_{j},T_s>
\\
& +\ \sum_{h=3}^r \sum_{p=1}^{m_h} c_p <\nabla_{X_i} X_{j},X_{h,p}>
\end{align*}

Now \eqref{deij} gives $<\nabla_{X_i} X_{j},X_\ell> = 0$, whereas
using \eqref{Koszul} again, we find
\[
2 <\nabla_{X_i} X_{j},T_s>\ =\ -\ <T_s,[X_j,X_i]>\ =\ \sum_{p=1}^k
b^p_{ij} \delta_{sp}\ =\ b^s_{ij}\ .
\]

Similarly, for $h\in \{3,...,r\}$ we have
\[
2 <\nabla_{X_i} X_{j},X_{h,p}>\ =\ \sum_{s=1}^k b^s_{ij}
<X_{h,p},T_s>\ =\ 0\ .
\]

From these equations we obtain
\[
<\nabla_{X_i} X_{j},Z>\ =\ <\frac{1}{2} \sum_{s=1}^k b^s_{ij}
T_s,Z>\ .
\]

From the arbitrariness of $Z$ we conclude that \eqref{delijZ} holds.
In a similar way, one obtains \eqref{delTsp}, and \eqref{delXT}. We
leave the details to the reader.

\end{proof}

\medskip

Next, we want to introduce a connection on the horizontal bundle. We
do this by projecting onto $H\bG$ the Levi-Civita connection
$\nabla$.

\medskip

\begin{dfn}\label{D:horconn}
If $X$ is a vector field on $\bG$, and $Y$ is a horizontal vector
field on $\bG$, then we define the (Levi-Civita) \emph{horizontal
connection} on $H\bG$ as follows
\begin{equation}\label{horconn}
\nabla^H_X Y\ \overset{def}{=}\ \sum_{i=1}^m <\nabla_X Y,X_i> X_i\ .
\end{equation}
\end{dfn}

\medskip

Let us notice that $\nabla^H$ satisfies the metric compatibility
condition
\begin{equation}\label{compa}
X<Y,Z>\ =\ <\nabla^H_X Y,Z>\ +\ <Y,\nabla^H_X Z>\ ,
\end{equation}
for every triple of vector fields $X, Y, Z$ on $\bG$, such that $Y$
and $Z$ are horizontal. This follows from the corresponding
compatibility condition \eqref{compcond} satisfied by the
Levi-Civita connection $\nabla$, and from the definition of
$\nabla^H$. From Proposition \ref{P:covder} and Definition
\ref{D:horconn} we obtain.

\medskip

\begin{prop}\label{P:hcovder}
Let $\bG$ be a Carnot group of step $r$, then
\begin{equation}\label{delijZ}
\nabla^H_{X_i} X_j\ =\ 0\ ,\quad\quad i,j=1,...,m\ .
\end{equation}
\begin{equation}\label{delTsp}
 \nabla^H_{T_p} T_s\ =\ 0\ ,\quad\quad p,s
= 1,...,k\ .
\end{equation}
\begin{equation}\label{delXT}
\nabla^H_{X_i} T_s\ =\  -\ \frac{1}{2} \sum_{j=1}^m b^s_{ij} X_j\ ,\
i=1,...,m,\ s=1,...,k\ .
\end{equation}
\end{prop}

\medskip

\begin{rmrk}\label{R:horcon}
We mention that the horizontal Levi-Civita connection $\nabla^H_X Y$
is intimately connected with the notion of non-holonomic connection
introduced by Cartan in his address at the 1928 International
Congress of Mathematicians in Bologna \cite{C}. In this respect we
refer the reader to the interesting re-visitation of Cartan's
address by Koiller, Rodrigues and Pitanga, see \cite{KRP1},
\cite{KRP2}, where the authors generalize some of the ideas in
\cite{C} and also introduce a non-holonomic connection (see their
Definition 1.1 in \cite{KRP1}) which, for a Carnot group, gives
precisely our Definition \ref{D:horconn}. A general framework has
been recently set forth by Hladky and Pauls in \cite{HP} for what
they call \emph{vertically rigid spaces}. These are sub-Riemannian
manifolds which include, in particular, Carnot groups. When
specialized to Carnot groups, the adapted connection in \cite{HP}
coincides with the horizontal connection in Definition
\ref{D:horconn}.
\end{rmrk}

\medskip

Hereafter, for a given vector field $X$ we indicate with $X^H =
\sum_{i=1}^m <X,X_i> X_i$ the projection of $X$ on the horizontal
bundle $H\bG$.

\medskip

\begin{prop}\label{P:symmetry}
Given horizontal vector fields $X$ and $Y$, one has
\[
\nabla^H_X Y\ -\ \nabla^H_Y X\ =\ [X,Y]^H\ \overset{def}{=}\
\sum_{i=1}^m <[X,Y],X_i> X_i\ .
\]
\end{prop}

\begin{proof}[\textbf{Proof}]
From Definition \ref{D:horconn} and the torsion freeness
\eqref{torsionfree} of the Levi-Civita connection we obtain
\begin{align*}
\nabla^H_X Y - \nabla^H_Y X\ & =\ \sum_{i=1}^m <\nabla_X Y -
\nabla_Y X,X_i> X_i \\
& =\ \sum_{i=1}^m <[X,Y],X_i> X_i\ =\ [X,Y]^H\ .
\end{align*}

\end{proof}

\medskip

If we define the \emph{horizontal torsion} as follows
 \[
 T^H(X,Y)\ =\ \nabla^H_X Y\ -\ \nabla^H_Y X\ -\ [X,Y]^H\ ,
 \]
 then Proposition \ref{P:symmetry} asserts that the horizontal
 connection is torsion free, and this is why we call it the
 horizontal Levi-Civita connection. Permuting cyclically the roles of $X,Y,Z$ in \eqref{compa}, and
using Proposition \ref{P:symmetry}, we obtain the following
\emph{horizontal Koszul identity} for $\nabla^H$.

\medskip

\begin{prop}\label{P:KoszulH}
Let $X,Y,Z$ be horizontal vector fields on $\bG$, then
\begin{align}\label{KoszulH}
2 <\nabla^H_X Y,Z>\ & =\ X<Y,Z>  + Y<X,Z> - Z<X,Y>
\\
& -\ <Y,[X,Z]^H> - <X,[Y,Z]^H> + <Z,[X,Y]^H>\ \ . \notag
\end{align}
\end{prop}

\medskip

Proposition \ref{P:KoszulH} shows in particular that $\nabla^H$ is
completely determined by the Riemannian inner product in $\bG$ and
by the horizontal bundle $H\bG$. Given a function $u\in C^1(\bG)$,
its Riemannian gradient with respect to the inner product
$<\cdot,\cdot>$ is given by
\begin{equation}\label{grad}
\nabla u\ =\ X_1u\ X_1 + ... + X_mu\ X_m + T_1u\ T_1 + ... + T_ku\
T_k + ... + X_{r,1}u\ X_{r,1} + ... + X_{r,m_r}u\ X_{r,m_r}\ ,
\end{equation}

If we let $G = det(g_{ij})$, then as a consequence of
\eqref{divfreeE}, and of the fact that $G \equiv 1$ (see \cite{F2},
or \cite{CGr}), we obtain for the divergence of $X_i$ (see \cite{He}
p.387)
\begin{equation}\label{divfree}
div\ X_i = \frac{1}{\sqrt{G}} \sum_{k=1}^N \frac{\p}{\p \xi_k}
\left(\sqrt{G} (X_i)_k\right) = div_{E} X_i + \sum_{k=1}^N (X_i)_k
\frac{\p}{\p \xi_k} (\log \sqrt{G}) = 0\ ,
\end{equation}
for every $i = 1,...,m$. The horizontal gradient of $u$ is obtained
by projecting $\nabla u$ on the subbundle $H\bG$ (see Definition
\ref{D:horconn}). The resulting horizontal vector field on $\bG$ is
nothing but the horizontal connection acting on $u$
\begin{equation}\label{hg}
\nabh u\ =\ <\nabla u,X_1> X_1 + ... + <\nabla u, X_m> X_m\ =\ X_1u\
X_1 + ... +  X_mu\ X_m\ .
\end{equation}

If $\zeta = \zeta_1 X_1 + ... + \zeta_m X_m \in C^1(\bG,H\bG)$, then
the horizontal divergence of $\zeta$ is given by
\begin{equation}\label{hordiv}
div_H \zeta\ =\ X_1 \zeta_1 + ... + X_m\zeta_m\ .
\end{equation}

The horizontal Laplacian (also known as sub-Laplacian) of a function
$u\in C^2(\bG)$ is given by
\begin{equation}\label{sl}
\Delta_H u\ =\ div_H \nabh u\ =\  \sum_{i=1}^m X_i^2\ .
\end{equation}

Except for the Abelian case when the step $r=1$ and $\delh$ is just
the standard Laplacian $\Delta = \sum_{i=1}^m \p^2/\p x_i^2$, such
operator fails to be elliptic at every point of $\bG$. We notice
that $\Delta_H u = trace (\nabla_H^2 u)$, where we have denoted by
$\nabla_H^2 u$ the $m\times m$ matrix-valued function on $\bG$
defined by \begin{equation}\label{horhess} \nabla_H^2 u\ =\ u_{,ij}\
=\ \frac{X_iX_ju + X_jX_iu}{2}\ ,\quad\quad\quad i,j=1,...,m\ .
\end{equation}

The following proposition contains a useful property of Carnot
groups.

\medskip

\begin{prop}\label{P:coord}
Let $\bG$ be a Carnot group, then
\begin{equation}\label{coord1-1}
X_i x_j\ =\ \delta_{ij}\ ,\quad\quad\quad\quad\quad \delh x_j\ =\ 0\
,\  i , j = 1, ... , m\ .
\end{equation}
As a consequence, we find
\begin{equation}\label{coord2-1}
|\nabh(|x|^2)|^2\ =\ 4\ |x|^2\ .
\end{equation}
One also has
\begin{equation}\label{coord3}
X_i t_s\ =\ \frac{1}{2}\ <[\xi_1, e_i],\ep_s>\ =\ \frac{1}{2}\
\sum_{j=1}^m x_j b^{s}_{ji},\ \quad\quad\quad X_j X_i t_s\ =\
\frac{1}{2}\ b^s_{ji}\ .
\end{equation}
In particular, we obtain $\nabla^2_H(t_s) = 0$, and therefore $\delh
t_s = 0$, $s = 1, ... , k$.
\end{prop}

\vskip 0.6in


\section{\textbf{Horizontal Gauss map and tangent space to a hypersurface}}\label{S:XMC}

\vskip 0.2in

In this section we introduce two basic geometric concepts for an
hypersurface in a Carnot group $\bG$ which are adapted to the
horizontal subbundle of $\bG$. We consider the Riemannian manifold
$M = \bG$ with the metric tensor with respect to which
$X_1,...,X_m, ... , X_{r,m_r}$ is an orthonormal basis, the
corresponding Levi-Civita connection $\nabla$ on $\bG$, and the
horizontal Levi-Civita connection $\nabla^H$ introduced in
Definition \ref{D:horconn}. Let $\mathcal S \subset \bG$ be a
$C^k$ oriented hypersurface, with $k\geq 2$. We will denote by
$\boldsymbol N$ the non-unit Riemannian normal to $\mathcal S$,
and will indicate with $\boldsymbol \nu = \boldsymbol
N/|\boldsymbol N|$ the Riemannian Gauss map of $\mathcal S$. It
will be convenient to introduce the following notation
\begin{equation}\label{ps}
p_j\ =\ <\boldsymbol N, X_j>\ ,\quad\quad\quad i = j,...,m\ ,\quad\quad\quad W\ =\ \sqrt{p_1^2\ +\ ...\ +\ p_m^2}\ .
\end{equation}

We now set
\begin{equation}\label{pbargen}
\pb_j\ =\ \frac{p_j}{W}\ ,\quad\quad\quad \text{so that}\quad\quad
\pb_1^2\ +\ ...\ +\ \pb_m^2\ \equiv\ 1\quad\quad \text{on}\quad\quad
\mathcal S\setminus \Sigma_\mS\ .
\end{equation}

We also define
\begin{equation}\label{omegas}
\om_s\ =\ <\bN,T_s>\ ,\ \quad\quad\ob_s\ =\ \frac{\om_s}{W}\ ,\ s
= 1,...,k\ ,
\end{equation}
\[
\om_{j,s}\ =\ <\bN,X_{j,s}>\ ,\quad\quad \ob_{j,s}\ =\
\frac{\om_{j,s}}{W}\ ,\quad\quad\quad j = 1,...,r,\quad s =
1,...,m_j\ .
\]

If $g_0\in \mathcal S$ is characteristic, then we have $p_j(g_0) =
0$, $j=1,...,m$, and therefore we have the alternative
characterization of $\Sigma_\mS$ as the zero set of the continuous
function $W$
\begin{equation}\label{csaf}
\Sigma_\mS\ =\ \{g\in \mathcal S\mid W(g) = 0\}\ ,
\end{equation}
which shows that $\Sigma_\mS$ is a closed subset of $\mS$. The next
definition plays a basic role in the sequel.

\medskip

\begin{dfn}\label{D:HGauss}
We define the \emph{horizontal normal} $\bN^H : \mathcal S
\rightarrow H\bG$ by the formula
\begin{equation}\label{up1}
\bN^H\ =\ \sum_{j=1}^m <\bN,X_j> X_j\ =\ \sum_{j=1}^m p_j\ X_j\ .
\end{equation}
The \emph{horizontal Gauss map} $\nuX$ is defined by
\begin{equation}\label{up2}
\nuX\ =\ \frac{\bN^H}{|\bN^H|}\ =\ \sum_{j=1}^m \pb_j\ X_j\
,\quad\quad\quad\quad \text{on}\quad \mathcal S \setminus
\Sigma_\mS\ .
\end{equation}

\end{dfn}

\medskip

We note that $\Up$ is the projection of the Riemannian normal $\bN$
on the horizontal subbundle $H\bG \subset T \bG$. Such projection
vanishes only at characteristic points, and this is why the
horizontal Gauss map is not defined on $\Sigma_\mS$. A trivial
consequence of the definition which, however, will be important in
the sequel is
\begin{equation}\label{normone}
|\nuX|^2\  =\ \pb_1^2\ +\ ...\ +\ \pb_m^2\ \equiv \ 1\
,\quad\quad\quad\quad \text{in}\quad \mathcal S \setminus
\Sigma_\mS\ ,
\end{equation}
which is of course a re-formulation of the second equation in \eqref{pbargen}.
One also has
\begin{equation}\label{Ynu}
<\nuX , \Up>\ =\ |\Up|\ ,\quad\quad\quad \Up\ -\ <\Up,\nuX>\ \nuX\
=\ 0\ .
\end{equation}

We note explicitly that, with these quantities in place, the
Riemannian (non-unit) normal to $\mS$ is given at every $g\in \mS
\setminus \Sigma_\mS$ by
\begin{align}\label{rn}
\bN\ & =\ \bN^H\ +\ \om_1 T_1 + ... + \om_k T_k + ... + \om_{r,m_r}
X_{r,m_r}
\\
& =\ W \bigg\{\pb_1 X_1 + ... + \pb_m X_m + \ob_1 T_1 + ... +
\ob_k T_k + ... + \ob_{r,m_r} X_{r,m_r}\bigg\} \notag\\
& =\ W \bigg\{\nuX + \ob_1 T_1 + ... + \ob_k T_k + ... +
\ob_{r,m_r} X_{r,m_r}\bigg\}\  . \notag
\end{align}

Since $<\nuX,T_s> = <\nuX,X_{j,m_j}> = 0$ for $s=1,...,k$, and $j
= 3,...,r$, it is obvious from \eqref{rn} that
\begin{equation}\label{ip}
<\bN,\nuX>\  =\ <\bN^H,\nuX>\ =\ W\ ,\ \quad\quad
\quad\text{hence}\quad\quad \cos(\nuX \angle \bN)\ =\
\frac{W}{|\bN|}\ .
\end{equation}

Because of \eqref{ip}, the function $W$ is also called the
\emph{angle function}.

\medskip

\begin{rmrk}\label{R:pq}
To help the reader's comprehension, we sometimes give proofs or
examples in the special case when $\bG = \HH$, the first Heisenberg
group. Furthermore, Sections \ref{S:geomid}, \ref{S:1&2var} and
\ref{S:stab} are devoted to this special setting. It will thus be
convenient to simplify the notation introduced above as follows. For
surfaces $\mS\subset \HH$ we will let
\begin{equation}\label{pqH} p = p_1\ ,\  q = p_2\ ,\ \om = \om_1\
,\ W = \sqrt{p^2 + q^2}\ ,\ \pb = \pb_1\ ,\ \qb = \pb_2\ ,\ \ob =
\ob_1\ . \end{equation}

Consequently, in this setting the normal $\bN$ and the horizontal
Gauss map $\nuX$ will always be respectively written as
\begin{equation}\label{n&hgm}
 \bN\ =\ p\ X_1\ +\ q\ X_2\ +\ \om\ T\  =\ \bN^H\ +\ \om\ T\  ,\quad\quad\quad \nuX\ =\ \pb\ X_1\ +\ \qb\ X_2\
 ,
\end{equation}
so that \eqref{rn} becomes \[ \bN = W \big\{\nuX\ +\ \ob\ T\big\}\
.
\]

The horizontal vector field defined on $\mS\setminus \Sigma_\mS$ by
\begin{equation}\label{nupH}
\nup\ =\ \qb\ X_1\ -\ \pb\ X_2\ ,
\end{equation}
is perpendicular to $\nuX$, but it is also orthogonal to the
Riemannian normal $\bN$ to $\mS$.
\end{rmrk}

\medskip

\begin{dfn}\label{D:hts}
At a point $g \in \mathcal S \setminus \Sigma_\mS$ the
\emph{horizontal tangent space} is defined as follows
\[
HT_{g}\mS\ \overset{def}{=}\ \left\{\boldsymbol v \in H_g \mid
<\boldsymbol v,\nuX>_g\ =\ 0\right\}\ .
\]
The \emph{horizontal tangent bundle} of $\mathcal S$ is defined by
\[
HT \mS\ =\ \underset{g\in \mathcal S\setminus \Sigma_\mS}{\bigcup}\
HT_{g} \mathcal S\ .
\]
\end{dfn}

\medskip

One can check that $HT \mS$ has the structure of a vector bundle.
It is clear that, since $dim\ H_g = m$, then $dim\ HT_{g} \mathcal
S = m -1$, and one has in fact
\begin{equation}\label{split}
H_g\ =\ HT_{g}\mathcal S \ \oplus\ span\ \{\nuX(g)\}\ .
\end{equation}

For instance, when $\bG = \HH$, then if for a $C^2$ surface $\mS
\subset \HH$ we consider the unit vector field on $\mS$ given by
\eqref{nupH},
 then it is clear that at every
point $g\in \mS\setminus \Sigma$, one has
\begin{equation}\label{spick&span}
HT_{g} \mathcal S\ =\ span\{(\nuX)^\perp(g)\}\ .
\end{equation}

If we consider $\Gamma = \mS \cap H_g$, then from Proposition
\ref{P:submanifold} we know that $\Gamma$ is submersed manifold of
dimension one (a curve). Its Riemannian tangent space in $g$ can be
identified in a canonical way with $HT_{g}$. We also observe that an
orthonormal basis for the Riemannian tangent space $T_g\mS$ of $\mS$
at $g$ is given by
\begin{equation}\label{basists}
T_g \mS\ =\ span \left\{\nup\ ,\ \frac{\omega}{|\bN|} \nuX -
\frac{W}{|\bN|} T\right\}\ .
\end{equation}

\medskip

\begin{prop}\label{P:proj}
Let $g\in \mathcal S \setminus \Sigma_\mS$, then one has \[
HT_{g}\mathcal S\ =\ T_g \mS\ \cap\ H_g\ .
\]
\end{prop}

\begin{proof}[\textbf{Proof}]
We begin by observing that, since by hypothesis $g\not\in
\Sigma_\mS$, then $H_g \not\subset T_g \mathcal S$, and therefore
$\bN^H \not= 0$. Now, from \eqref{rn} and the fact that
$X_1,...,X_m, T_1,...,X_{r,m_r}$ constitute an orthonormal basis of
$T_g\mS$ at every $g\in \bG$, one sees from \eqref{rn} that $\bN -
\bN^H \perp H_g$. Therefore, if $\boldsymbol v\in HT_{g}\mathcal S$,
then we have
\begin{equation}\label{ortho}
<\boldsymbol v,\bN>\ =\ <\boldsymbol v,\bN - \bN^H>\ +\ <\boldsymbol
v, \bN^H>\ =\ 0\ ,
\end{equation}
which shows $\boldsymbol v \in T_g\mS$. We thus have the inclusion
$HT_{g} \mS \subset T_g\mathcal S \cap H_g$. To establish the
opposite inclusion, let $\boldsymbol v\in T_g\mathcal S \cap H_g$.
We thus have that the left-hand side of \eqref{ortho} is zero, and
since $<\boldsymbol v,\bN - \bN^H> = 0$ because of the fact that
$\boldsymbol v \in H_g$, we conclude that it must be $<\boldsymbol
v, \bN^H> = 0$, hence $\boldsymbol v \in HT_{g} \mS$.

\end{proof}

\vskip 0.6in


\section{\textbf{Horizontal connection on a
hypersurface}}\label{S:hc}

\vskip 0.2in

We recall the classical definition of the Levi-Civita connection
of a $n$-dimensional immersed submanifold $N=N^n$ of an
$m$-dimensional Riemannian manifold $M=M^m$. Denoting with $i : N
\hookrightarrow M$ the immersion, and having endowed $N$ with the
induced Riemannian metric $i^* g$, let $i_*:TN \to TM$ be the
differential of $i$. We identify $T_p N$ with the subspace
$(i_*)_p(T_p N)$ of $T_pM$, and denote by $T_pN^\perp$ its
orthogonal complement. $TN^\perp = \bigcup_{p\in N} T_p N^\perp$
has the structure of a $(m-n)$-dimensional vector bundle,
traditionally referred to as the normal bundle of $N$. We can thus
write $TM\mid N \cong TN \oplus TN^\perp$, and for every
$\boldsymbol u\in T_p M$, we indicate with $\boldsymbol u^\top$
its $T_pN$ component, and with $\boldsymbol u^\perp$ its
$T_pN^\perp$ component. Since $p \to (\nabla^M_X Y)^\top(p)$
satisfies all the assumptions of a Levi-Civita connection on $N$,
by the uniqueness of the latter we obtain
\begin{equation}\label{lc}
 \nabla^N_X Y\ =\ (\nabla^M_X Y)^\top\ .
\end{equation}

Before proceeding we need to say a few words concerning
\eqref{lc}. First of all, since $X, Y$ are only initially defined
on the submanifold $N$, we need to give a meaning to right-hand
side. Using a partition of unity argument, we can extend $X, Y$ to
smooth vector fields $\oX, \oY$ on $M$, and therefore interpret
the right-hand side as follows
\begin{equation}\label{lc2}
(\nabla^M_X Y)^\top\ =\ (\nabla^M_{\oX} \oY)^\top\ .
\end{equation}

This immediately raises the question of whether \eqref{lc2} is a
good definition, in other words, whether it is independent of the
particular extensions of $X, Y$ that we have picked. Since the
value of $\nabla^M_X Y$ at $p\in N$ depends only on $X_p$, it is
clear that \eqref{lc2} is independent of the extension of $X$. On
the other hand, $(\nabla^M_X Y)_p$ depends only on the values of
$Y$ along any curve on $M$ whose initial tangent vector is $X_p$.
By picking a curve which lies entirely on $N$, we see that
\eqref{lc2} is also independent of the extension of $Y$. This
fact, can be also recognized by the following observations, which
also establish the torsion freeness of the connection $(\nabla^M_X
Y)^\top$. Denoting with $i_*$ the differential of the immersion,
we have $i_*(X) = \oX$, $i_*(Y) = \oY$, and therefore, see Theorem
7.9 in \cite{Bo}, $i_*[X,Y] = [\oX,\oY]$. This implies, in
particular, that $[X,Y] = [\oX,\oY]^\top$. From the torsion
freeness of $\nabla^M$, we thus conclude that
\begin{equation}\label{lc3}
(\nabla^M_{\oX} \oY)^\top\ -\ (\nabla^M_{\oY} \oX)^\top\ =\
[\oX,\oY]^\top\ =\ [X,Y]\ .
\end{equation}

We notice that $[X,Y]_p$ only depends on the values of $X,Y$ in a
neighborhood of $p$ in $N$, and therefore \eqref{lc3} shows at
once that \eqref{lc2} is a good definition, and that $(\nabla^M_X
Y)^\top$ is torsion free. The remaining properties of a
Levi-Civita connection are checked easily.

Inspired by the Riemannian situation we now introduce a notion of
horizontal connection on a hypersurface $\mS\subset \bG$ by
projecting the horizontal Levi-Civita connection $\nabla^H$ in the
ambient Lie group $\bG$ onto the horizontal tangent space $HT\mS$.

\medskip

\begin{dfn}\label{D:horconS}
Let Let $\mS\subset \bG$ be a non-characteristic, $C^k$
hypersurface, $k\geq 2$, then we define the \emph{horizontal
connection} on $\mS$ as follows. Let $\nabla_H$ denote the
horizontal Levi-Civita connection introduced in Definition
\ref{D:horconn}. For every $X,Y\in C^1(\mS;HT\mS)$ we define
\[
\delXY\ =\ \nabla^H_{\oX} \oY\ -\ <\nabla^H_{\oX} \oY,\nuX> \nuX\
,
\]
where $\oX, \oY$ are any two horizontal vector fields on $\bG$
such that $\oX = X$, $\oY = Y$ on $\mS$.
\end{dfn}

\medskip

Arguing as above one can check that Definition \ref{D:horconS} is
well-posed, i.e., it is independent of the extensions $\oX, \oY$
of the vector fields $X, Y$. From Proposition \ref{P:symmetry} we
immediately obtain.

\medskip

\begin{prop}\label{P:nonsymmetry}
For every $X,Y\in C^1(\mS;HT\mS)$ one has
\[
\delXY\ -\ \nabla^{H,\mS}_Y X\ =\ [X,Y]^H\ -\ <[X,Y]^H,\nuX>\nuX\
.
\]
\end{prop}

\medskip

It is clear from this proposition that the horizontal connection
$\nabla^{H,\mS}$ on $\mS$ is not necessarily torsion free. This
depends on the fact that it is not true in general that, if $X,Y\in
C^1(S;HT\mS)$, then $[X,Y]^H \in C^1(\mS;HT\mS)$. In the special
case of the first Heisenberg group this fact is true, and we have
the following result.

\medskip

\begin{prop}\label{P:symmH1}
Given a $C^k$ non-characteristic surface $\mS \subset \HH$, $k\geq
2$, one has $[X,Y]^H\in HT\mS$ for every $X,Y\in C^1(\mS;HT\mS)$,
and therefore the horizontal connection on $\mS$ is torsion free.
\end{prop}

\begin{proof}[\textbf{Proof}]
According to \eqref{spick&span}, for every $g\in \mS$ we have
$HT_{g} \mathcal S = span\{\boldsymbol e_1(g)\}$, where $\boldsymbol
e_1 = (\nuX)^\perp$. Therefore, if we take two vector fields $X,Y
\in C^1(\mS;HT\mS)$, then we can write $X= a \boldsymbol e_1$, $Y= b
\boldsymbol e_1$, for appropriate $C^{k-1}$ functions $a$ and $b$.
We thus have
\[
[X,Y]\ =\ [a \boldsymbol e_1,b \boldsymbol e_1]\ =\ \left\{a
\boldsymbol e_1(b) - b \boldsymbol e_1(a)\right\} \boldsymbol e_1
\]

This shows that $[X,Y]\in C(\mS;HT\mS)$, and therefore Proposition
\ref{P:nonsymmetry} gives
\[
\delXY\ -\ \nabla^{H,\mS}_Y X\ =\ [X,Y]\ .
\]

This gives the desired conclusion.

\end{proof}

\medskip

\begin{dfn}\label{D:delta}
Let $\mS$ be as in Definition \ref{D:horconS}. Consider a function
$u \in C^1(\mS)$. We define the \emph{tangential horizontal
gradient} of $u$ as follows
\[
\del u\ \overset{def}{=}\ \nabh \ou\ -\ <\nabh \ou,\nuX>\ \nuX\ ,
\]
where $\ou\in C^1(\bG)$ is such that $\ou = u$ on $\mS$.
\end{dfn}

\medskip

We note that $\del u = \sum_{i=1}^m \di u\ X_i$, where
\[
\di u\ =\ X_i \ou\ -\  <\nabh \ou,\nuX>\ \nui\ =\  X_i \ou\ -\
\pb_i\ \pb_j\ X_j\ou\  . \]

Observe also that $\del u \in HT\mS$. One has in fact from
\eqref{normone} and Definition \ref{D:delta}
\begin{equation}\label{perp}
<\del u, \nuX>\ \equiv\ 0\ \quad\quad\quad\text{in}\quad \mathcal
S\setminus \Sigma\ ,
\end{equation}
and therefore
\begin{equation}\label{squares}
|\del u|^2\ =\ |\nabh u|^2\ -\ <\nabh u,\nuX>^2\ .
\end{equation}

\vskip 0.6in


\section{\textbf{Perimeter measure and horizontal first fundamental form}}\label{S:pm}

\vskip 0.2in

In a Carnot group $\bG$, given an open set $\Om\subset \bG$, we
let
\[
\mathcal F(\Om)\ =\ \{\zeta = \sum_{i=1}^m \zeta_i X_i \in
C^1_0(\Om,H\bG)\ \mid\ |\zeta|_\infty\ =\ \sup_{\Om}\
\bigg(\sum_{i=1}^m \zeta_i^2\bigg)^{1/2}\ \leq 1\}\ .
\]

For a  function $u\in L^1_{loc}(\Om)$, the $H$-variation of $u$
with respect to $\Om$ is defined by
\[
Var_H(u;\Om)\ =\ \underset{\zeta\in \mathcal F(\Om)}{\sup}\
\int_{\bG} u\ \text{div}_H \zeta\ dg\ .
\]

We say that $u\in L^1(\Om)$ has bounded $H$-variation in $\Om$ if
$Var_H(u;\Om) <\infty$. The space $BV_H(\Om)$ of functions with
bounded $H$-variation in $\Om$, endowed with the norm
\[
||u||_{BV_H(\Om)}\ =\ ||u||_{L^1(\Om)}\ +\ Var_H(u;\Om)\ ,
\]
is a Banach space.

\medskip

\begin{dfn}\label{D:perimeter}
Let $E\subset \bG$ be a measurable set, $\Om$ be an open set. The
$H$-perimeter of $E$ with respect to $\Om$ is defined by
\[
P_H(E;\Om)\ =\ Var_H(\chi_E;\Om)\ ,
\]
where $\chi_E$ denotes the indicator function of $E$. We say that
$E$ is a $H$-Caccioppoli set if $\chi_E \in BV_H(\Om)$ for every
$\Omega \subset \subset \bG$.
\end{dfn}

\medskip

The above definitions are taken from \cite{CDG1}, see also
\cite{GN1}. Following classical arguments \cite{Z}, \cite{EG}, one
obtains from the Riesz representation theorem.

\medskip

\begin{thrm}\label{T:structure}
Given an open set $\Om\subset \bG$, let $E\subset \bG$ be a
$H$-Caccioppoli set in $\Om$. There exist a Radon measure
$||\partial^H E||$ in $\Om$, and a $||\partial^H E||$-measurable
function $\nuX_E : \Om \to H\bG$, such that
\[
|\nuX_E(g)|\ =\ 1 \quad\quad\quad\text{for}\quad ||\partial^H
E||-a.e.\quad g\in \Om\ ,
\]
and for which one has for every $\zeta\in C^1_0(\Om;H\bG)$
\[
\int_E \text{div}_H \zeta\ dg\ =\ \int_\Om\ < \zeta, \nuX_E>\
d||\partial^H E||\ =\ \int_\Om\ < \zeta,d [\partial^H E]>\  .
\]
\end{thrm}

\medskip

Let $E\subset \bG$ be a $C^1$ domain, with Riemannian outer unit
normal $\n$. If $\zeta\in C^1_0(\Om;H\bG)$, we have
\begin{equation*}
\int_E\ \text{div}_H \zeta\ dg\ =\ \int_{\partial E \cap \Om}\
\sum_{i=1}^m\ \zeta_i\ <X_i,\n>\  d H_{N-1}\ .
\end{equation*}

From this observation, and from Theorem \ref{T:structure}, we
conclude the following result.

\medskip

\begin{prop}\label{P:smoothXvar}
Let $E\subset \bG$ be a $C^1$ domain. For every open set $\Om
\subset \bG$, and any $\zeta\in C^1_0(\Om;H\bG)$, one has
\[
\int_\Om\ < \zeta, \nuX_E>\ d||\partial^H E||\ =\ \int_{\partial E
\cap \Om}\ <\zeta, \frac{\Up}{|\bN|}>\  d H_{N-1}\ ,
\]
where $\Up$ is defined in \eqref{up1}. Moreover,
\begin{equation}\label{peri1}
d||\partial^H E||\ =\ |\Up|\ d \left(H_{N-1}\lfloor \partial
E\right)\ ,
\end{equation}
and one has
\begin{equation}\label{peri2}
||\partial^H E||(\Om)\ =\ P_H(E;\Om)\ =\ \int_{\partial E \cap \Om}\
\frac{|\Up|}{|\bN|}\ d H_{N-1}\ =\ \int_{\partial E \cap \Om}\
\frac{W}{|\bN|}\ d H_{N-1}\ ,
\end{equation}
where $W$ is the angle function defined in \eqref{ps}.
\end{prop}

\medskip

\begin{dfn}\label{D:permeas}
Given an oriented $C^2$ hypersurface $\mathcal S\subset \bG$, we
will denote by
\begin{equation}\label{permeasure}
d \sigma_H\ =\ \frac{|\Up|}{|\bN|}\ dH_{N-1}\lfloor \mathcal S\ =\
\frac{W}{|\bN|}\ dH_{N-1}\lfloor \mathcal S\ ,
\end{equation}
the $H$-perimeter measure supported on $\mathcal S$ (see
\eqref{peri2} and \eqref{ip}).
\end{dfn}

\medskip

For a detailed local study of such measure the reader should see
\cite{DGN1}, \cite{DGN2}, \cite{Ma1}, \cite{Ma2}. An interesting
interpretation of the $H$-perimeter measure is that the latter is
obtained by blowing-up the Riemannian regularizations of the
sub-Riemannian metric of the group $\bG$. In a different context,
this idea was first exploited systematically by Kor\'anyi \cite{Ko1}
in his computations of the sub-Riemannian geodesics in $\Hn$. For
simplicity, and to illustrate the main idea, we will state the
relevant result in the case when $\bG = \Hn$.

\medskip

\begin{thrm}\label{T:rr}
Consider in the Heisenberg group $\Hn$ the left-invariant Riemannian
metric tensor $\{g^\ep_{ij}\}_{i,j=1,...,2n +1}$ with respect to
which $\{X_1,...,X_{2n}, \sqrt{\ep} T\}$ constitutes an orthonormal
frame of $T\Hn$. Let $\mS\subset \Hn$ be a $C^2$ hypersurface, with
$\Sigma_\mS = \varnothing$, and denote by $I^\mS_\ep(\cdot,\cdot)$
the first fundamental form in the Riemannian metric on $\mS$ induced
by $\{g^\ep_{ij}\}_{i,j=1,...,2n +1}$. Denote by $\sigma^\ep$ the
corresponding surface area on $\mS$, then for any bounded open chart
$U\subset \mS$ one has
\[
\sigma_H(U)\ =\ \underset{\ep \to 0}{\lim}
\frac{\sigma^\ep(U)}{\sqrt{det(g^\ep_{ij})}}\ . \]
\end{thrm}

\begin{proof}[\textbf{Proof}]
For simplicity, we present the proof in the case $n=1$. Let $T_\ep
= \sqrt \epsilon T$, and consider in $\HH$ the one-parameter
family of left-invariant Riemannian metrics
$\{(g^\ep_{ij})_{i,j=1,2,3}\}_{\ep>0}$ with respect to which
$\{X_1,X_2,T_\ep\}$ constitute an orthonormal basis of $T\HH$.
Similarly to \eqref{gij}, we find
\begin{equation}\label{gijep}
(g^\ep_{ij})\ =\ \begin{pmatrix} 1 + \frac{y^2}{4\ep} & -
\frac{xy}{4\ep} & \frac{y}{2\ep}
\\
- \frac{xy}{4\ep} & 1 + \frac{x^2}{4\ep} & - \frac{x}{2\ep}
\\
\frac{y}{2\ep} & - \frac{x}{2\ep} & \frac{1}{\ep}
\end{pmatrix}\ .
\end{equation}

One easily verifies that \[ G^\ep\ =\ det(g^\ep_{ij})\ =\
\ep^{-1}\ , \] and that letting $((g^\ep)^{ij}) =
(g^\ep_{ij})^{-1}$, then
\begin{equation}\label{gijinverseep}
((g^\ep)^{ij})\ =\ \begin{pmatrix} 1  & 0 & - \frac{y}{2}
\\
0 & 1  &  \frac{x}{2}
\\
- \frac{y}{2} &  \frac{x}{2} & \ep + \frac{x^2 + y^2}{4}
\end{pmatrix}
\end{equation}

Let $\Om \subset \R^2$ be a bounded open set such that $U$ is
represented by $\theta : \Om \to U$, with $\theta\in C^2(\Om)$. We
have $\theta(u,v) = x(u,v) X_1 + y(u,v) X_2 + t(u,v) T$, see
\eqref{thetaH}. We now use some of the computations from Section
\ref{S:1&2var}. We rewrite \eqref{dertheta} as follows
\begin{equation}\label{derthetaep}
\begin{cases}
\theta_u\ =\  x_u X_1 + y_u X_2 + \frac{1}{\sqrt{\ep}}\left(t_u +
\frac{y x_u - x
y_u}{2}\right) T_\ep\ , \\
\theta_v\ =\  x_v X_1 + y_v X_2 + \frac{1}{\sqrt{\ep}}\left(t_v +
\frac{y x_v - x y_v}{2}\right)  T_\ep\ . \end{cases}
\end{equation}

Denoting by $\wedge_\ep$ the wedge product with respect to the
orthonormal frame $\{X_1,X_2,T_\ep\}$, similarly to \eqref{nun} we
obtain for the non-unit Riemannian normal to $\mathcal S$ with
respect to $I_\ep(\cdot,\cdot)$
\begin{align}\label{nunep}
\bN^\ep\ =\  \theta_u\ \wedge_\ep\ \theta_{v}\ & =\
\frac{1}{\se}\left(y_u t_v - y_v t_u\ -\ \frac{y}{2} (x_u y_v -
x_v
y_u)\right)\ X_1 \\
& +\ \frac{1}{\se}\left(x_v t_u - x_u t_v\ +\ \frac{x}{2} (x_u y_v
- x_v
y_u)\right)\ X_2 \notag\\
& +\ (x_u y_v - x_v y_u)\ T_\ep \notag\\
& =\  \frac{1}{\se} p X_1 + \frac{1}{\se} q X_2 + \om  T_\ep\ ,
\notag
\end{align}
where in the last equality we have used \eqref{ppar}. From
\eqref{nunep} we conclude that
\begin{equation}\label{pparep}
\frac{\sigma_\ep(U)}{\sqrt{g^\ep}}\ =\ \se \int_U d\sigma_\ep\ =\
=\ \se \int_\Om \sqrt{I_\ep(\bN^\ep,\bN^\ep)} du\wedge dv\ = \
\int_\Om \sqrt{p^2 + q^2 + \ep \om^2}\ du\wedge dv\ .
\end{equation}

Letting $\ep \to 0$ in \eqref{pparep}, we conclude
\[
\underset{\ep \to 0}{\lim} \frac{\sigma^\ep(U)}{\sqrt{g^\ep}}\ =\
\int_\Om W\ du\wedge dv\ =\ \int_U \frac{W}{|\bN|} d\sigma\ =\
\sigma_H(U)\ ,
\]
where in the last equality we have used \eqref{permeasure}. This
completes the proof.

\end{proof}

\medskip

We close this section by collecting two basic properties of the
$H$-perimeter. The former is a trivial consequence of the
left-invariance on the vector fields $X_1,...,X_m$, and of the
definition of $H$-perimeter.

\medskip

\begin{prop}\label{P:invariance2}
For any $H$-Caccioppoli set $E$ in a Carnot group $\bG$, and any
open set $\Om\subset \bG$, one has
\begin{equation}\label{t1}
P_H(L_{g_0}(E) ;L_{g_0}(\Om))\ =\ P_H(E;\Om)\ ,\quad\quad\quad g_0
\in \bG\ ,
\end{equation}
where $L_{g_0}g = g_0 g$ is the left-translation on the group. In
particular,
\begin{equation}\label{l2}
P_H(L_{g_0}(E) ;\bG)\ =\ P_H(E;\bG)\ ,\quad\quad\quad g_0 \in \bG\
.
\end{equation}
\end{prop}

\medskip

\medskip

\begin{prop}\label{P:blowupP}
In a Carnot group $\bG$ one has for every $H$-Caccioppoli set
$E\subset \bG$, any open set $\Om\subset \bG$, and every $\la>0$
\begin{equation}\label{Pscaling}
P_H(\delta_{\la} E ; \delta_{\la} \Om)\ =\ \la^{Q-1}\ P_H(E;\Om)\
.
\end{equation}
In particular,
\begin{equation}\label{dilations2}
P_H(\delta_{\la} E ; \bG)\ =\ \la^{Q-1}\ P_H(E;\Om)\ .
\end{equation}
\end{prop}

\begin{proof}[\textbf{Proof}]
We observe that if $\zeta \in C^1_0(\bG,H\bG)$, then $\zeta \circ
\delta_{1/\la} \in C^1_0(\delta_{\la} \Om;H\bG)$. Furthermore,
$\zeta\in \mathcal F(\Om)$ if and only if $\zeta\circ
\delta_{1/\la} \in \mathcal F(\delta \Om)$. The divergence
theorem, and a rescaling now give
\begin{equation}\label{ibyparts1}
\int_E\ div_H \zeta\ dg\ =\ \int_E\ \sum_{j=1}^m X_j \zeta_j\ dg\
=\ \la^{-Q}\ \int_{\delta_\la E}\ \sum_{j=1}^m X_j
\zeta_j(\delta_{1/\la} g)\ dg\ .
\end{equation}

Since
\[
X_j(\zeta_j \circ \delta_{1/\la})\ =\ \la^{-1}\ (X_j \zeta_j)
\circ \delta_{1/\la}\ ,
\]
we conclude from \eqref{ibyparts1}
\[
\int_E\ div_H \zeta\ dg\ =\ \la^{-(Q-1)}\ \int_{\delta_\la E}\
\sum_{j=1}^m X_j(\zeta_j\circ \delta_{1/\la})\ dg\ .
\]

Taking the supremum on all $\zeta \in \mathcal F(\Om)$ in the
latter equation, we reach the desired conclusion.

\end{proof}

\medskip

Combining Propositions \ref{P:invariance2} and \ref{P:blowupP} we
obtain the following result.

\medskip

\begin{cor}\label{blowup}
Let $\mS \subset \bG$ be a $C^2$ hypersurface with finite
$H$-perimeter, then for every $g_0\in \bG$, and every $\lambda >0$,
one has
\[
\sigma_H(L_{g_0}(\mS))\ =\ \sigma_H(\mS)\ ,
\]
\[
\sigma_H(\delta_\lambda \mS)\ =\ \lambda^{Q-1}\ \sigma_H(\mS)\ .
\]
\end{cor}

\vskip 0.6in


\section{\textbf{Horizontal second fundamental form and mean curvature}}\label{S:hmc}

\vskip 0.2in

We open this section by computing the first variation of the
$H$-perimeter for deformations of a hypersurface $\mathcal S$
along the Riemannian normal $\bN$ to $\mS$. This will provide a
first motivation for the introduction of the notion of $H$-mean
curvature.

\medskip

\begin{thrm}\label{T:fv}
Let $\mathcal U \subset \bG$ be a bounded open set and consider
$\phi \in C^2(\mathcal U)$ with $|\nabla \phi|\geq \alpha
> 0$ in $\mathcal U$, and for small $\la \in [-\la_0,\la_0]$ consider
the one-parameter family of $\mS^\la = \p \mathcal U_\la$, where
we have let $\mathcal U_\la = \{g\in \mathcal U \mid \phi(g) <
\la\}$. Assume that each of the $\mS^\la$ be a $C^2$
non-characteristic hypersurface. Let $\mS = \mS^0$ and define a
function $\mathcal H : \mS \to \R$ by letting
\begin{equation}\label{mcmot}
\mathcal H\ \overset{def}{=}\ \sum_{i=1}^m X_i \pb_i\ ,
\end{equation}
where the $\pb_i$ are the components of the horizontal Gauss map
introduced in \eqref{pbargen}. We then have
\[
\frac{d}{d\la}\ P_H(\mS^\la)\Bl\ \overset{def}{=}\ \frac{d}{d\la}\
P_H(\mathcal U_\la;\bG)\Bl\ =\ \int_{\mS} \frac{\mathcal
H}{|\bN|}\ dH_{N-1}\ ,
\]
In particular, $\mS$ is a critical point of the $H$-perimeter with
respect to the deformations $\mS\to \mS^\la$ if and only of
$\mathcal H \equiv 0$.
\end{thrm}

\begin{proof}[\textbf{Proof}]
Using Federer's coarea formula \cite{Fe} we can write
\begin{equation}\label{per1gen}
\int_{\mathcal U_\la} |\nabh\phi|\ dg\ =\ \int_{- \infty}^\la
\int_{\p \mathcal U_\tau} \frac{W}{|\bN|}\ dH_{N-1}\ d\tau\ =\
\int_{- \infty}^\la P_H(\mathcal U_\tau;\bG)\ d\tau\ ,
\end{equation}
where the second equality is a consequence of \eqref{peri2}. The
identity \eqref{per1gen} gives
\begin{equation}\label{per2gen}
P_H(\mathcal U_\la;\bG)\ =\ \frac{d}{d\la}\ \int_{\mathcal U_\la}
|\nabh\phi|\ dg\ .
\end{equation}

Using the summation convention over repeated indices, and
integration by parts, we now compute
\begin{align}\label{per3}
\int_{\mathcal U_\la} |\nabh\phi|\ dg\ & =\ \int_{\mathcal U_\la}
X_i \phi\ \nui\ dg
\\
& =\ \int_{\p \mathcal U_\la} \phi <X_i,\n> \nui\ dH_{N-1}\ -\
\int_{\mathcal U_\la} \phi\ X_i \nui\ dg
\notag\\
& =\ \int_{\p \mathcal U_\la} \phi\ |\Up| dH_{N-1}\ -\
\int_{\mathcal U_\la} \phi\ X_i \nui\ dg
\notag\\
& =\ \la\ P_H(\mathcal U_\la;\bG)\ -\ \int_{\mathcal U_\la} \phi\
X_i \nui\ dg\ , \notag
\end{align}
where we have used \eqref{peri2}. From \eqref{per2gen},
\eqref{per3}, and the coarea formula again, we find
\begin{equation}\label{per4}
P_H(\mathcal U_\la;\bG)\ =\ \frac{d}{d\la}\left\{\la\ P_H(\mathcal
U_\la;\bG)\right\}\ -\ \la\ \int_{\p \mathcal U_\la}  \frac{X_i
\pb_i}{|\bN|} dH_{N-1}\ .
\end{equation}

Equation \eqref{per4} easily implies the desired conclusion.

\end{proof}

\medskip

\begin{rmrk}\label{R:fv}
As we will see in this section, the critical points of the
$H$-perimeter are precisely the so-called $H$-minimal hypersurfaces.
We will return to this question in Section \ref{S:1&2var}, where we
will develop more precise intrinsic first and second variation
formulas of the $H$-perimeter in the setting of the Heisenberg
group.
\end{rmrk}

\medskip

We are now ready to introduce the central notions of sub-Riemannian,
or horizontal second fundamental form, and of $H$-mean curvature.
According to Theorem \ref{T:fv}, hypersurfaces for which the
function $\mathcal H$ in \eqref{mcmot} vanishes identically on $\mS$
are critical points of the $H$-perimeter functional with respect to
deformations of the surface in the direction of the Riemannian
normal $\bN$ to $\mS$. This suggests a notion of horizontal mean
curvature of $\mS$ based on the equation \eqref{mcmot}. Such notion
was proposed by one of us back in 1997, see \cite{G1}, and it
produces precisely the function in \eqref{mcmot}. We will in fact
introduce a more intrinsic notion which is based on that of
horizontal second fundamental form, and then recognize that such
definition coincides with \eqref{mcmot}. This closely parallels the
classical development of the subject. We recall the classical
definition of the mean curvature $H$ of a $n$-dimensional immersed
submanifold $N=N^n$ of an $m$-dimensional Riemannian manifold
$M=M^m$. Denoting with $i : N \hookrightarrow M$ the immersion, we
recall that the Levi-Civita connection of $N$ is given by the
equation \eqref{lc}. The second fundamental form of $N$ is defined
by
\[
II_N(X,Y)\ =\ (\nabla^{M}_X Y)^\perp\ ,
\]
where $\nabla^{M}$ is the Levi-Civita connection of $M$. Since for
vector fields on $N$ one has
\[
II_N(X,Y) - II_N(Y,X)\ =\  (\nabla^{M}_X Y - \nabla^{M}_Y X)^\perp\
=\ [X,Y]^\perp\ =\ 0\ , \] $II_N$ defines a symmetric tensor field
on $N$ of type $(0,2)$, which takes values in $TN^\perp$. The mean
curvature of $N$ at a point $p\in N$ is defined by \[ H\ =\ -\
\frac{1}{n}\ trace(II_N)\ =\ -\ \frac{1}{n} \sum_{i=1}^n
II_N(\boldsymbol e_i,\boldsymbol e_i)\ , \] where $\{\boldsymbol
e_1,...,\boldsymbol e_n\}$ is an orthonormal basis of $T_pN$.

We now consider the Riemannian manifold $M = \bG$ with the metric
tensor with respect to which $X_1,...,X_m, ... , X_{r,m_r}$ is an
orthonormal basis, and the corresponding Levi-Civita connection
$\nabla$ on $\bG$. Let $\nabla^H$ denote the horizontal
Levi-Civita connection introduced in Definition \ref{D:horconn}.
Let $\mS\subset \bG$ be a $C^2$ hypersurface. Inspired by the
Riemannian situation we introduce a notion of horizontal second
fundamental on $\mS$ as follows.

\medskip

\begin{dfn}\label{D:sff}
Let $\mS\subset \bG$ be a $C^2$ hypersurface, with $\Sigma_\mS =
\varnothing$, then we define a tensor field of type $(0,2)$ on
$HT\mS$, as follows: for every $X,Y\in C^1(\mS;HT\mS)$
\begin{equation}\label{sff} II^{H,\mS}(X,Y)\ =\ <\nabla^H_X Y,\nuX>
\nuX\ .
\end{equation}
We call $II^{H,\mS}(\cdot,\cdot)$ the \emph{horizontal second
fundamental form} of $\mS$. We also define $\mathcal A^{H,\mS} :
HT \mS \to HT \mS$ by letting for every $g\in \mS$ and
$\boldsymbol u, \boldsymbol v \in HT_{g}$
\begin{equation}\label{shape}
<\mathcal A^{H,\mS} \boldsymbol u,\boldsymbol v>\ =\ -\
<II^{H,\mS}(\boldsymbol u,\boldsymbol v),\nuX>\ =\ -\ <\nabla_X^H
Y,\nuX>\ ,
\end{equation}
where $X, Y \in C^1(\mS,HT\mS)$ are such that $X_g = \boldsymbol u$,
$Y_g = \boldsymbol v$. We call the endomorphism $\mathcal A^{H,\mS}
: HT_{g}\mS \to HT_{g}\mS$ the \emph{horizontal shape operator}. If
$\boldsymbol e_1,...,\boldsymbol e_{m-1}$ denotes a local
orthonormal frame for $HT\mS$, then the matrix of the horizontal
shape operator with respect to the basis $\boldsymbol
e_1,...,\boldsymbol e_{m-1}$ is given by the $(m-1)\times(m-1)$
matrix $ \big[- <\nabla_{\boldsymbol e_i}^H \boldsymbol
e_j,\nuX>\big]_{i,j=1,...,m-1}$.

\end{dfn}

\medskip

Using the horizontal Koszul identity \eqref{KoszulH}, one easily
verifies that
\[
<\nabla_{\boldsymbol e_i}^H \boldsymbol e_j,\nuX>\ =\ -\
<\nabla_{\boldsymbol e_i}^H \nuX, \boldsymbol e_j>\ .
\]

Using Proposition \ref{P:symmetry} in Definition \ref{D:sff} we
immediately recognize that
\begin{equation}\label{nonsymm}
II^{H,\mS}(X,Y)\ -\ II^{H,\mS}(Y,X)\ =\ <[X,Y]^H,\nuX> \nuX\ ,
\end{equation}
and therefore, unlike its Riemannian counterpart, the horizontal
second fundamental form of $\mS$ is not necessarily symmetric. This
depends on the fact, already observed, that if $X,Y\in
C^1(\mS;HT\mS)$, then it is not necessarily true that $[X,Y]^H\in
C(\mS;HT\mS)$. The next proposition gives a necessary and sufficient
condition for the symmetry of $II^{H,\mS}$.

\medskip

\begin{prop}\label{P:symm}
The horizontal second fundamental form $II^{H,\mS}(\cdot,\cdot)$
is a $(0,2)$ symmetric tensor field on $HT\mS$ if and only if for
any local orthonormal basis $\boldsymbol e_1,...,\boldsymbol
e_{m-1}$ of $HT\mS$, one has
\[
[\boldsymbol e_i,\boldsymbol e_j]^H\ \in\ HT\mS\ ,\quad\quad i, j
=1,...,m-1\ .
\]
\end{prop}

\begin{proof}[\textbf{Proof}]
Let $X= \sum_{i=1}^{m-1} a_i \boldsymbol e_i$, $Y=
\sum_{j=1}^{m-1} b_j \boldsymbol e_j$, then
\[
[X,Y]\ =\ \sum_{j=1}^{m-1} \left\{\sum_{i=1}^{m-1} \big(a_i
\boldsymbol e_i(b_j) - b_i \boldsymbol e_i(a_j)\big)\right\}
\boldsymbol e_j\ +\ \sum_{i,j=1}^{m-1} a_i b_j [\boldsymbol
e_i,\boldsymbol e_j]\ .
\]

This identity gives
\begin{align*}
<[X,Y]^H,\nuX>\  & =\ \sum_{i,j=1}^{m-1} \big[a_i \boldsymbol
e_i(b_j) - b_i \boldsymbol e_i(a_j)\big]<\boldsymbol e_j,\nuX>\ +\
\sum_{i,j=1}^{m-1} a_i b_j <[\boldsymbol e_i,\boldsymbol
e_j]^H,\nuX>
\\
& =\ \sum_{i,j=1}^{m-1} a_i b_j <[\boldsymbol e_i,\boldsymbol
e_j]^H,\nuX>\ =\ 0\ ,
\end{align*}
provided that $[\boldsymbol e_i,\boldsymbol e_j]^H\in HT\mS$.
Therefore, under the assumption $[\boldsymbol e_i,\boldsymbol
e_j]^H \in HT\mS$, we finally obtain from Definition \ref{D:sff}
\[
II^{H,\mS}(X,Y)\ -\ II^{H,\mS}(Y,X)\ =\ <\nabla^H_X Y - \nabla^H_Y
X,\nuX> \nuX\ =\ <[X,Y]^H,\nuX> \nuX\ =\ 0\ ,
\]
which proves the symmetry of $II^{H,\mS}$. Vice-versa, suppose
that $II^{H,\mS}$ be symmetric, then applying the latter identity
with $X=\boldsymbol e_i$, $Y = \boldsymbol e_j$, we reach the
conclusion that $<[\boldsymbol e_i,\boldsymbol e_j]^H,\nuX> = 0$.

\end{proof}

\medskip

\begin{cor}\label{C:sffH}
The horizontal second fundamental form of a $C^2$
non-characteristic surface $\mS \subset \HH$ is symmetric.
\end{cor}

\begin{proof}[\textbf{Proof}]
In this situation the assumption of Proposition \ref{P:symm} is
trivially satisfied since a basis of $HT\mS$ is given by the
single vector field $\boldsymbol e_1 = \nup$, and therefore
$[\boldsymbol e_1,\boldsymbol e_1]^H = 0 \in HT\mS$, see also
Proposition \ref{P:symmH1}.

\end{proof}

\medskip

Another situation in which the assumption of Proposition
\ref{P:symm} is fulfilled is that when $\mathcal S$ is a
cylindrical hypersurface over the first layer of the Lie algebra.

\medskip

\begin{prop}\label{P:rigate}
Suppose that the hypersurface $\mathcal S$ is a \emph{vertical
cylinder}, i.e., it can be represented in the form
\begin{equation}\label{vertical}
\mathcal S\ =\ \{g\in \bG \mid \mathfrak h(x_1(g),...,x_m(g)) =
0\}\ ,
\end{equation}
where $\mathfrak h \in C^2(\R^m)$, and there exist an open set
$\omega\subset \R^m$ and $\alpha>0$ such that $|\nabla \mathfrak
h|\geq \alpha$ in $\omega$. Under these assumptions, we have
$\Sigma_\mS = \varnothing$, and the horizontal second fundamental
form is symmetric.
\end{prop}

\begin{proof}[\textbf{Proof}]
The function $\phi(g) = \mathfrak h(x_1(g),...,x_m(g))$ is a
defining function of $\mathcal S$. Using the global exponential
coordinates, we obtain from Lemma \ref{L:zero}
\[
X_i\phi(g)\ =\ \frac{\partial \mathfrak h}{\partial x_i} \ ,
\]
hence $\nabh \phi = \nabla_x \mathfrak h$, which proves in
particular that $\Sigma_\mS = \varnothing$, and that $\nuX =
\frac{\nabla \mathfrak h}{|\nabla \mathfrak h|} = \n$. We next
observe that for every $g_0\in \bG$, the left-translated surface
$\tilde \mS = L_{g_0}(\mS)$ is again a vertical cylinder, with
defining function $\tilde{\mathfrak h}(x) = \mathfrak h(x(g_0) +
x(g)) =0$. As a consequence of this observation, if $g_0\in \mS$,
then by left-translation we can assume without restriction that $g_0
= e$. We thus immediately see that the horizontal plane $H_e =
\exp(V_1)$ is given by $t_1=...=t_k=...= x_{r,m_r}= 0$. Furthermore,
by an orthogonal transformation in the horizontal layer of the Lie
algebra, we can assume that the tangent space of $\mS$ in $e$ be
given by the hyperplane $x_m= 0$. Since thanks to Proposition
\ref{P:proj} the horizontal tangent space at $e$ is given by $H_e
\cap T_e\mS$, from the previous considerations we see that $HT_e\mS
= span\{\boldsymbol e_1,...,\boldsymbol e_{m-1}\}$, where
$\boldsymbol e_i = (\p/\p x_i)_e$. Since $[\p/\p x_i,\p/\p x_j] =
0$, $i,j=1,...,m-1$, we conclude that $[\boldsymbol e_i ,\boldsymbol
e_j]^H = 0$. In view of Proposition \ref{P:symm} we conclude that
$II^{H,\mS}$ is symmetric, thus completing the proof.

\end{proof}

\medskip

From the proof of Proposition \ref{P:rigate} one easily obtains the
following corollary.

\medskip

\begin{cor}\label{C:rigatoni}
Let $\mathcal S$ be a vertical cylinder as in \eqref{vertical}, then
the $H$-mean curvature at $g\in \mS$ is given by the formula
\begin{equation}\label{Hvertical}
\mathcal H(g)\ =\ (m-1)\ H(x(g))\ ,
\end{equation}
where $H(x(g))$ represents the Riemannian mean curvature of the
projection $\pi_{V_1}(\mathcal S)$ of $\mathcal S$ onto the
horizontal layer $V_1$. In particular, $\mathcal S$ is $H$-minimal
 if and only if $\pi_{V_1}(\mathcal S)$ is a classical minimal surface in $V_1 \simeq \R^m$.
\end{cor}

\medskip

\begin{dfn}\label{D:HMC}
We define the \emph{horizontal principal curvatures} as the real
eigenvalues $\kappa_1,...,\kappa_{m-1}$ of the symmetrized
operator
\[
\mathcal A^{H,\mS}_{sym}\ =\ \frac{1}{2}\left\{\mathcal A^{H,\mS}
+ (\mathcal A^{H,\mS})^t\right\}\ ,
\]
The $H$-mean curvature of $\mS$ at a non-characteristic point
$g_0\in \mS$ is defined as
\[
\mathcal H\ =\ -\ trace\ \mathcal A^{H,\mS}_{sym}\ =\
\sum_{i=1}^{m-1} \kappa_i\ =\ \sum_{i=1}^{m-1}
<\nabla^H_{\boldsymbol e_i} \boldsymbol e_i,\nuX>\ .
\]
If $g_0$ is characteristic, then we let
\[
\mathcal H(g_0)\ =\ \underset{g\to g_0, g\in \mathcal S\setminus
\Sigma_\mS}{\lim}\ \mathcal H(g)\ ,
\]
provided that such limit exists, finite or infinite. We do not
define the $H$-mean curvature at those points $g_0\in \Sigma_\mS$ at
which the limit does not exist. Finally, we call $\vec{\mathcal H} =
\mathcal H \nuX$ the $H$-\emph{mean curvature vector}.
\end{dfn}

\medskip

\begin{prop}\label{P:equalMC}
The $H$-mean curvature in Definition \ref{D:HMC} coincides with
the function defined in \eqref{mcmot}. In fact, the following
intrinsic identity holds
\begin{equation}\label{equal2}
\mathcal H\ =\ \sum_{i=1}^m\ \di\ <\nuX,X_i>\ =\  \sum_{i=1}^m\
\di\ \pb_i\ .
\end{equation}
\end{prop}

\begin{proof}[\textbf{Proof}]
In what follows to simplify the exposition we continue to indicate
with $\pb_1,...,\pb_m$ an $m$-tuple of $C^1$ extensions to the
whole of $\bG$ of the coefficients of the horizontal Gauss map
with respect to the basis $X_1,...,X_m$. We begin by observing
that using the horizontal Koszul identity \eqref{KoszulH}, one
easily recognizes that
\[
<\nabla_{\boldsymbol e_i}^H \boldsymbol e_j,\nuX>\ =\ -\
<\boldsymbol e_i, [\boldsymbol e_i, \nuX]^H>\ .
\]

From Definition \ref{D:HMC} we thus obtain
\begin{equation}\label{equalMC}
\mathcal H\ =\ -\ \sum_{i=1}^{m-1} <\nabla^H_{\boldsymbol e_i}
\boldsymbol e_i,\nuX>\ =\ \sum_{i=1}^{m-1} <\boldsymbol e_i,
[\boldsymbol e_i, \nuX]^H>\ .
\end{equation}

Recalling \eqref{up2}, we find
\[
[\boldsymbol e_i,\nuX]\ =\ \sum_{j=1}^m \boldsymbol e_i(\pb_j) X_j
+ \sum_{j=1}^m [\boldsymbol e_i,X_j]\ .
\]

To proceed in the calculations, we write
\[
\boldsymbol e_i\ =\ \sum_{\ell = 1}^m a^\ell_i X_\ell\
,\quad\quad\quad i =1,...,m-1\ ,
\]
with $\{a^\ell_i\}$ satisfying the orthogonality conditions
\begin{equation}\label{o1}
\sum_{\ell = 1}^m a^\ell_i\ \pb_\ell\ =\ 0\ ,\quad\quad\
\sum_{\ell = 1}^m a^\ell_i\ a^\ell_j\ =\ \delta_{ij}\ ,\quad \quad
i,j=1,...,m-1\ .
\end{equation}

We thus obtain
\begin{align*}
[\boldsymbol e_i,\nuX]\ & =\ \sum_{j=1}^m \boldsymbol e_i(\pb_j)
X_j - \sum_{\ell =1}^m \sum_{j=1}^m X_j(a^\ell_i) X_\ell
\\
& + \ \sum_{s=1}^k \left(\sum_{\ell =1}^m \sum_{j=1}^m a^\ell_i
b^s_{\ell j}\right) T_s\ .
\end{align*}

The latter identity gives
\[
[\boldsymbol e_i,\nuX]^H\  =\ \sum_{j=1}^m \boldsymbol e_i(\pb_j)
X_j - \sum_{\ell =1}^m \sum_{j=1}^m X_j(a^\ell_i) X_\ell\ ,
\]
and therefore,
\begin{align}\label{horsff}
\sum_{i=1}^{m-1} <[\boldsymbol e_i,\nuX]^H,\boldsymbol e_i>\ & =\
\sum_{j=1}^m \sum_{i=1}^{m-1} \boldsymbol e_i(\pb_j) a^j_i\ -\
\sum_{\ell,j=1}^m \sum_{i=1}^{m-1} a^\ell_i X_j(a^\ell_i)
\\
& =\ \sum_{j=1}^m \sum_{i=1}^{m-1} \boldsymbol e_i(\pb_j) a^j_i
\notag\\
& =\ \sum_{\ell,j=1}^m X_\ell(\pb_j) \sum_{i=1}^{m-1} a^\ell_i
a^j_i\ , \notag
\end{align}
where in the second to the last equality we have used \eqref{o1}.
We now observe that, since $\{\boldsymbol e_1,...,\boldsymbol
e_{m-1},\nuX\}$ is an orthonormal basis of $H_gT\mS$, we have
\begin{align*}
\sum_{i=1}^{m-1} a^\ell_i a^j_i\ & =\ \sum_{i=1}^{m-1}
<X_\ell,\boldsymbol e_i><X_j,\boldsymbol e_j>
\\
& =\ <X_\ell,X_j>\ -\ <X_\ell,\nuX> <X_j,\nuX>\ =\ \delta_{\ell
j}\ -\ \pb_\ell \pb_j\ .
\end{align*}

Substituting this identity in \eqref{horsff}, and recalling
\eqref{equalMC}, we finally have
\[
\mathcal H\ =\ \sum_{i=1}^{m-1} <[\boldsymbol
e_i,\nuX]^H,\boldsymbol e_i>\ =\ \sum_{j=1}^m X_j(\pb_j)\ -\
\sum_{\ell,j=1}^m \pb_j X_\ell(\pb_j)\ =\ \sum_{j=1}^m X_j(\pb_j)\
.
\]

This concludes the proof.

\end{proof}

\medskip

It is clear from Definition \ref{D:HMC} that $\mathcal H\in
C(\mS\setminus \Sigma)$.

\medskip

\begin{dfn}\label{D:MS}
An oriented $C^k$ hypersurface $\mathcal S \subset \bG$, $k\geq
2$, is said to have \emph{constant $H$-mean curvature} if
$\mathcal H \equiv const$ on $\mS$. We say that $\mS$ is
$H$-\emph{minimal} if its $H$-mean curvature $\mathcal H$ vanishes
everywhere as a continuous function on $\mS$.
\end{dfn}

\medskip

\begin{rmrk}\label{R:productgroup1}
Consider the product group $\hat{\bG} = \bG \times \R$ with the
canonical group law $(g,s) \circ (g',s') = (g g', s+ s')$, induced
by the one on $G$. The stratification of the Lie algebra for
$\hat{\bG}$ is then given by $\hat{V_1}\oplus\cdots\oplus \hat{V_r}$
where $\hat{V_1} = V_1 \times \R$, $\hat{V_j} = V_j \times \{0\}$
for $j = 2,3,...,r$. If $\{e_1,...,e_m\}$ is a basis for $V_1$, we
let $\hat{e_j} = (e_j,0)$ for $j=1,...,m$,
$\hat{e}_{m+1}=(0,...,0,1)$.  A basis for $\hat{V_1}$ is then given
by $\{\hat{e_1},...,\hat{e}_{m+1}\}$. This identifies a subbundle
$H\hat \bG$. We can naturally identify $\bG$ with the hypersurface
$\mathcal S = \bG \times \{0\} \subset \hat{\bG}$ with global
defining function $ \phi(g,s)=s$. Now observe that
$\nabla^{H,\mS}_{i}\phi = X_i\phi$, $i=1,...,m$ and
$\nabla^{H,\mS}_{m+1}\phi = 1$ and therefore, $\nuX = \hat e_{m+1}$.
As a consequence, we have $\nabla^{H,\mS}_{i}\,\nui = 0$ for $i =
1,...,m+1$. In view of Proposition \ref{P:equalMC} we conclude that
the $H$-mean curvature of $\bG$ in $\hat \bG$ is zero.
\end{rmrk}

\medskip

To state the next proposition we consider for a function $u : \bG
\to \mathbb R$ the \emph{symmetrized horizontal Hessian} of $u$
 at $g\in \bG$. This is the $m\times m$ matrix with entries
\begin{equation}\label{hessian}
u_{,ij}\ \overset{def}{=}\ \frac{1}{2}\ \bigg\{X_iX_j u\ +\ X_j X_i u\bigg\}\ ,\quad\quad\quad\quad i , j = 1 , ... , m\ .
\end{equation}

Setting $\nabla^2_H u = [u_{,ij}]$, the mapping $g\to \nabla^2_H
u(g)$ defines a $2$-covariant tensor on the subbundle $H\bG$. We
recall that the horizontal Laplacian associated with the basis
$\{e_1,...,e_m\}$ of $V_1$ is given by $\delh  u = tr\ \nabla^2_H\
u$. We also consider the following nonlinear operator
\begin{equation}\label{infty}
\delinf u\ \overset{def}{=}\ \sum_{i,j=1}^m u_{,ij}\ X_i u\ X_j u\
=\ \frac{1}{2} <\nabh(|\nabh u|^2), \nabh u>\ ,
\end{equation}
which, by analogy with its by now classical Euclidean ancestor, we
call the \emph{horizontal} $\infty$-\emph{Laplacian}. The reason for
introducing the operator $\delinf$ is in the following result which
is often useful in computing the $H$-mean curvature. We consider a
$C^2$ hypersurface $\mS\subset \bG$, and for a given $g_0\in
\mS\setminus \Sigma_\mS$, suppose that there exist a neighborhood
$\mathcal U$ of $g_0$, and $\phi\in C^2(\mathcal U)$, such that
$\mS\cap \mathcal U = \mS\cap \partial \{g\in \mathcal U\mid
\phi(g)<0\}$. We observe that the hypothesis that $g_0\not\in
\Sigma_\mS$ implies that $\nabh \phi(g_0)\not= 0$, and therefore, by
possibly restricting $\mathcal U$ we can assume that $\nabh
\phi(g)\not= 0$ for every $g\in \mS \cap \mathcal U$. We thus have
\begin{equation}\label{up3}
\Up(g)\ =\ \nabh \phi(g)\ , \quad\quad\quad\text{for every}\quad
g\in \mS \cap \mathcal U\ ,
\end{equation}
and therefore
\begin{equation}\label{nuX}
\nuX\ =\ \frac{\nabh\phi}{|\nabh\phi|}\ \quad\quad\quad\text{for
every}\quad g\in \mS \cap \mathcal U\ .
\end{equation}

\medskip

\begin{prop}\label{P:XMCG}
At every point of $\mS\cap \mathcal U$ one has
\begin{align*}
|\nabla_H \phi|^3\ \mathcal H\ & =\ \left\{|\nabla_H \phi|^2\
\Delta_H \phi\ -\ \Delta_{H,\infty} \phi\right\}\ .
\end{align*}
\end{prop}

\begin{proof}[\textbf{Proof}]
We use the summation convention over repeated indices. Invoking
Proposition \ref{P:equalMC} and \eqref{nuX}, we obtain at every
point in $\mS\setminus \Sigma$
\[
\mathcal H\ =\ \di \pb_i\ =\ X_i (\nu_{H,i})\ =\ X_i
\left(\frac{X_i\phi}{|\nabh\phi|}\right)\ =\
\frac{\delh\phi}{|\nabh\phi|}\ -\ \frac{\delinf
\phi}{|\nabh\phi|^3}\ .
\]

\end{proof}

\medskip

It is interesting to consider a nonlinear operator which
interpolates in an appropriate sense between the $H$-mean curvature
operator in Definition \ref{D:HMC}, and the operator $\delinf$.
Consider the one-parameter family of quasilinear operators defined
by
\begin{equation}\label{pilap}
\delp u\ =\ div_H(|\nabla^H u|^{p-2} \nabh u)\ =\ 0\
,\quad\quad\quad\quad 1 < p < \infty\ .
\end{equation}

Supposing that $|\nabla^H u| \not= 0$, we formally rewrite in the
more suggestive fashion
\begin{equation}\label{pilap2}
\delp u\ =\ (p-2)\ |\nabla^H u|^{p-4}\ \left\{\frac{1}{p-2}\
|\nabla^H u|^2\ \delh u\ +\ \delinf u \right\}\ .
\end{equation}

If $u_p$ is a solution to $\delp u_p = 0$, then equation
\eqref{pilap2} gives
\[
\frac{1}{p-2}\ |\nabla^H u_p|^2\ \delh u_p\ +\ \delinf u_p\ =\ 0\ .
\]

If we assume that $u_p \to u_\infty$ as $p \to \infty$, and that
$|\nabla^H u_p|^2 \delh u_p$ is bounded independently of $p$ large,
then by letting $p\to \infty$ we formally find that $u_\infty$ must
be a solution to $\delinf u_\infty = 0$. On the other hand, if we
know instead that $u_p \to u_1$ as $p\to 1$, then we discover from
\eqref{pilap2} and Proposition \ref{P:XMCG} that
\[
\delp u_p\ \ \underset{p\to 1}{\longrightarrow}\ \ -\ \mathcal
H(u_1)\ ,
\]
where $\mathcal H(u_1)$ is the $H$-mean curvature of the level sets
of $u_1$! This suggests that one should study the behavior as $p\to
1$ of the one-parameter family of quasilinear operators $\delp$.

\medskip

\noindent \textbf{Comparison with S. Pauls' notion of horizontal
mean curvature}. In the first Heisenberg group $\HH$ another notion
of horizontal mean curvature was introduced by S. Pauls in
\cite{Pa}. Such notion is based on the procedure of Riemannian
$\epsilon$-regularization defined in the proof of Theorem
\ref{T:rr}. Using \eqref{gijinverseep} one sees that, given a
function $\phi \in C^1(\Hn)$, its gradient with respect to the
metric $(g^\ep_{ij})$ is given by \begin{equation}\label{gradep}
\nabla_\ep \phi\ =\ X_1\phi\ X_1\ +\ X_2 \phi\ X_2\ +\ T_\ep \phi\
T_\ep\ .
\end{equation}

Let us note in passing that \eqref{gradep} gives
\begin{equation}\label{gradepnorm}
|\nabla_\ep \phi|^2\ =\ |\nabla_H\phi|^2\ +\ \ep (T\phi)^2\
,\quad\text{and}\quad \Delta_\ep \phi\ =\ \Delta_H \phi\ +\ \ep\
T^2\phi\ ,
\end{equation}
where we have denoted by $\Delta_\ep$ the Laplace-Beltrami
operator with respect to the metric $(g^\ep_{ij})$.

In \cite{Pa} the author defined the horizontal mean curvature of
$\mS$ at a point $g_0\in \mS\setminus \Sigma_\mS$ as follows
\begin{equation}\label{pauls}
\mathcal H_P(g_0)\ \overset{def}{=}\ \underset{\ep \to 0}{\lim}\
H^\ep_\mathcal R(g_0)\ ,
\end{equation}
where $H^\ep_\mathcal R$ indicates the mean curvature of $\mS$ in
the Riemannian metric \eqref{gijep}. We now recognize that such
notion coincides with the one introduced in Definition
\ref{D:HMC}.

\medskip

\begin{prop}\label{P:pauls}
The horizontal mean curvature defined by \eqref{pauls} coincides
with the one in Definition \ref{D:HMC}.
\end{prop}

\begin{proof}[\textbf{Proof}]
Let $\mS\subset \HH$ be a $C^2$ surface. The Riemannian Gauss map of
$\mS$ with respect to the metric \eqref{gijep} is given by $\n^\ep =
\frac{\bN^\ep}{|\bN^\ep|_\ep}$, where we have denoted by
$|\bN^\ep|_\ep$ the length of $\bN^\ep$ in such metric. Let us
notice that $|\bN^\ep|_\ep = \frac{1}{\sqrt \ep} \sqrt{W^2 + \ep\
\omega^2}$. Recalling \eqref{nunep} we see that
\begin{equation}\label{gmep}
\n_\ep\ =\ \alpha^\ep\ \big\{\pb X_1 + \qb X_2 +  \sqrt \ep \ob
T_\ep\big\}\ ,
\end{equation}
where we have let $\alpha^\eps = W/\sqrt{W^2 + \ep\ \omega^2}$. From
\eqref{gmep} and \eqref{coord} below, we recognize that the
expression of the Gauss map in the Cartesian coordinates $(x,y,t)$
is given by
\[
\n_\ep\ =\ \left(\alpha^\ep\ \pb\ ,\ \alpha^\ep\ \qb\ ,\ \ep\
\alpha^\ep\ \ob + \frac{\alpha^\ep}{2} (x \qb - y \pb)\right)\ .
\]

Using Proposition \ref{P:equalMC}, \eqref{divfree}, and the fact
that $det(g^\ep_{ij}) = \ep^{-1}$, we then see that at any $g_0\in
\mS$
\begin{align*}
H^\ep_\mathcal R(g_0)\ & =\  div_\ep\ \n_\ep\  =\ \p_x(\alpha^\ep\
\pb) + \p_y(\alpha^\ep\ \qb)\ +\ \p_t(\ep\
\alpha^\ep\ \ob + \frac{\alpha^\ep}{2} (x \qb - y \pb)) \\
& =\ X_1(\alpha^\ep\ \pb) + X_2(\alpha^\ep\ \qb)\ +\ \ep\
T(\alpha^\ep\ \ob)
\\
& =\ \alpha^\ep\ \mathcal H\ +\ \pb X_1(\alpha^\ep)\ +\ \qb
X_2(\alpha^\ep)\ +\ \ep (\alpha^\ep T \ob + \ob T(\alpha^\ep))\ ,
\end{align*}
where we have used the fact that $\mathcal H = X_1 \pb + X_2 \qb$,
see Definition \ref{D:HMC} and Proposition \ref{P:equalMC}. We now
claim that at any $g_0\in \mS\setminus \Sigma$ we have $\alpha^\ep
\to 1$ as $\ep \to 0$, and that furthermore the following
cancelation relations hold
\begin{equation}\label{eplimits}
X_1(\alpha^\ep)\ \to\ 0\ ,\quad\quad X_1(\alpha^\ep)\ \to\ 0\
,\quad\quad T(\alpha^\ep)\ \to\ 0\ ,\quad\quad\text{as}\quad \ep
\to 0\ .
\end{equation}

We only prove the first relation of \eqref{eplimits}, leaving the
analogous details of the remaining two to the reader. We have
\begin{align*}
 X_1 \alpha^\ep\ & =\ \frac{(W^2 + \ep \omega^2) X_1 W - W (W X_1W + \ep \omega X_1 \omega)}{(W^2 + \ep \omega^2)^{3/2}}
 \\
 & =\ \ep\ \omega\ \frac{\omega X_1 W - W X_1 \omega}{(W^2 + \ep
 \omega^2)^{3/2}}\ \longrightarrow \ 0\ .
 \end{align*}

From \eqref{eplimits} we conclude that $\mathcal H_P(g_0) =
\underset{\ep \to 0}{\lim}\ H^\ep_\mathcal R(g_0) = \mathcal
H(g_0)$, at every $g_0\in \mS\setminus \Sigma$.

\end{proof}

\medskip

\noindent \textbf{Comparison with the notion of pseudo-hermitian
mean curvature of Cheng, Hwang, Malchiodi and Yang}. In \cite{CHMY}
the authors have introduced the following notion of pseudo-hermitian
curvature for a surface $\mS \subset M$, where $(M,J,\Theta)$ is a
three-dimensional oriented CR manifold, with CR structure $J$, and
global contact form $\Theta$. At every point $g\in \mS \setminus
\Sigma_\mS$, they consider the one-dimensional space $HT_g\mS$. They
fix a unit vector field $\boldsymbol e_1\in HT\mS$ with respect to
the metric $G=\frac{1}{2} d\Theta(\cdot,J\cdot)$ associated with the
Levi form, and then define $\boldsymbol e_2 = J(\boldsymbol e_1)$.
They call $\boldsymbol e_2$ the Legendrian normal or Gauss map. They
denote by $\nabla^{p.h.}$ the pseudo-hermitian connection associated
with $J,\Theta)$. There exists a function $H^{p.h.}$ such that
\begin{equation}\label{chmy}
\nabla^{p.h.}_{\boldsymbol e_1} \boldsymbol e_1\ =\ H^{p.h.}\
\boldsymbol e_2\ .
\end{equation}

Such function $H_{p.h.}$ is called the pseudo-hermitian mean
curvature of $\mS$, see (2.1) in \cite{CHMY}.

\medskip

\begin{prop}\label{P:chmy}
Let $M$ be the Heisenberg group $\HH$, then the function $H^{p.h.}$
coincides (up to a choice of the orientation) with the horizontal
$H$-mean curvature in Definition \ref{D:HMC}.
\end{prop}

\begin{proof}[\textbf{Proof}]
One can check that in the Heisenberg group the pseudo-hermitian
connection $\nabla^{p.h.}$ is nothing but the horizontal Levi-Civita
connection $\nabla^H$ introduced in section \ref{S:lcc}. Since a
basis for $HT\mS$ is given by $\boldsymbol e_1 = \nup$, see
Corollary \ref{C:sffH}, and we clearly have $\boldsymbol e_2 =
\nuX$, from the horizontal Koszul identity \eqref{KoszulH} we obtain
for every vector field $X = a \boldsymbol e_1 + b \boldsymbol e_2 +
c T$
\begin{align}\label{covariant0}
2 <\nabla^{H}_{\boldsymbol e_1} \boldsymbol e_1,X>\ & =\ 2
\boldsymbol e_1<\boldsymbol e_1,X> - X<\boldsymbol e_1,\boldsymbol
e_1>
\\
& -\ 2 <\boldsymbol e_1,[\boldsymbol e_1,X]^H> + <X,[\boldsymbol
e_1,\boldsymbol e_1]^H>
\notag\\
& =\ 2 \boldsymbol e_1(a)\ -\ 2 <\boldsymbol e_1,[\boldsymbol
e_1,X]^H>\  . \notag
\end{align}

We now have
\[
[\boldsymbol e_1,X]\ =\ \boldsymbol e_1(a) \boldsymbol e_1 +
\boldsymbol e_1(b) \boldsymbol e_2 + b [\boldsymbol e_1,\boldsymbol
e_2] + \boldsymbol e_1(c) T + c [\boldsymbol e_1,T]\ .
\]

The commutators in the right-hand side of the latter equation have
been computed in section \ref{S:geomid} below, where the vector
fields $\boldsymbol e_1$, and $\boldsymbol e_2$ are respectively
denoted by $Z$ and $Y$. From Lemma \ref{L:comm} and Lemma
\ref{L:commT} we find
\[
[\boldsymbol e_1,X]\ =\ \boldsymbol e_1(a) \boldsymbol e_1 +
\boldsymbol e_1(b) \boldsymbol e_2 + b \bigg\{T + \mathcal H
\boldsymbol e_1 + (\qb \boldsymbol e_2(\pb) - \pb \boldsymbol
e_2(\qb)) \boldsymbol e_2\bigg\} + \boldsymbol e_1(c) T + c (\qb
T\pb - \pb T\qb)\boldsymbol e_2\ .
\]

From the latter expression we obtain
\[
[\boldsymbol e_1,X]^H\ =\ \boldsymbol e_1(a) \boldsymbol e_1 +
\boldsymbol e_1(b) \boldsymbol e_2 + b \bigg\{\mathcal H \boldsymbol
e_1 + (\qb \boldsymbol e_2(\pb) - \pb \boldsymbol e_2(\qb))
\boldsymbol e_2\bigg\} + c (\qb T\pb - \pb T\qb)\boldsymbol e_2\ ,
\]
and therefore \eqref{covariant0} gives
\[
<\nabla^{H}_{\boldsymbol e_1} \boldsymbol e_1,X>\ =\ \boldsymbol
e_1(a) - \boldsymbol e_1(a) - b \mathcal H\ =\ <- \mathcal H
\boldsymbol e_2,X>\ .
\]

From the arbitrariness of $X$ we conclude
\begin{equation}\label{covariant}
\nabla^{H}_{\boldsymbol e_1} \boldsymbol e_1\ =\ -\ \mathcal H\
\boldsymbol e_2\ .
\end{equation}

This proves the proposition.

\end{proof}

\vskip 0.6in


\section{\textbf{Sub-Riemannian calculus on hypersurfaces}}\label{S:IBP}
\vskip 0.2in

In this section we establish some basic integration by parts
formulas involving the tangential horizontal gradient on a
hypersurface, and the horizontal mean curvature of the latter. Such
formulas are reminiscent of the classical one, and in fact they
encompass the latter. However, an important difference is that the
ordinary volume form on the hypersurface $\mS$ is replaced by the
$H$-perimeter measure $d\sigma_H$. Furthermore, they contain
additional terms which are due to the non-trivial commutation
relations, which is reflected in the lack of torsion freeness of the
horizontal connection on $\mS$. Such term prevents the corresponding
horizontal Laplace-Beltrami operator from being formally
self-adjoint in $L^2(\mathcal S,d\sigma_H)$ in general. Since the
framework we work in does not lend itself to a preferred choice of
coordinates, for ease of computation we have developed an approach,
based on Federer's co-area formula, which is coordinate-free. In
fact, our proof slightly simplifies several of the classical
formulas for hypersurfaces in $\Rn$ which are derived by writing
$\mathcal S$ as a graph, see e.g. \cite{MM}, \cite{Gi}.

\medskip

\begin{thrm}[\textbf{First sub-Riemannian integration by parts formula}]\label{T:ibp}
Consider a $C^2$ oriented hypersurface in a Carnot group $\mathcal
S\subset \bG$. If $u\in C^1_0(\mS \setminus \Sigma_\mS)$, then we
have
\begin{equation}\label{ibp0}
\int_\mathcal S \di u\ d \sigma_H\ =\ \int_\mathcal S u\
\bigg\{\mathcal H\ \nui\ -\ \boldsymbol
c^{H,\mS}_{i}\bigg\}d\sigma_H\ ,\quad\quad\quad i = 1,...,m\ ,
\end{equation}
where the $C^1$ functions $\boldsymbol c^{H,\mS}_{i}$ on $\mathcal
S\setminus \Sigma$ are defined by
\begin{equation}\label{ibp3}
\boldsymbol c^{H,\mS}_{i}\ =\ \sum_{s=1}^k \big(\sum_{j=1}^m
b^s_{ij} \pb_j \big) \ob_s\ ,
\end{equation}
with $b^s_{ij}$ denoting the horizontal group constants defined in
\eqref{gc}. Moreover, the horizontal vector field $\cx =
\sum_{i=1}^m \boldsymbol c^{H,\mS}_{i} X_i$ is perpendicular to the
horizontal Gauss map $\nuX$, i.e., one has
\begin{equation}\label{ibp4}
<\cx,\nuX>\ =\ 0\ .
\end{equation}
As a consequence, we have $\cx \in C^1(\mS\setminus \Sigma_\mS,
HT\mS)$.
\end{thrm}

\begin{proof}[\textbf{Proof}]
Since the question is local, to prove the theorem we will assume,
without loss of generality, that $\mS$ is the level set of a $C^2$
defining function $\phi$, and that $\mS$ is oriented in such a way
that $\bN = \nabla \phi$. Furthermore, thanks to the assumption of
the support of $u$, we can also assume that $\mS$ be
non-characteristic. Using a partition of unity we can always reduce
ourselves to this situation. For every $\rho\in \R$, we define
$\mathcal U_\rho = \{g\in \bG \mid \phi(g) < \rho\}$. By the
non-characteristic assumption on $\mS$ we can assume that, if $\mS =
\p \mathcal U_{\rho_0}$, then for every $\rho$ sufficiently close to
$\rho_0$ the characteristic locus of $\p \mathcal U_\rho$ is empty.
Federer's coarea formula gives, see \cite{Fe},
\begin{equation}\label{ibp5}
\int_{\mathcal U_\rho} \di u\ W\ dg\ =\ \int_{-\infty}^\rho \int_{\p
\mathcal U_\tau} \di u\ \frac{W}{|\bN|}\ dH_{N-1}\ d\tau\ =\
\int_{-\infty}^\rho \int_{\p \mathcal U_\tau} \di u\ d\sigma_H\
d\rho\ .
\end{equation}

Recalling \eqref{permeasure} and \eqref{up3}, we obtain from \eqref{ibp5}
\begin{equation}\label{ibp6}
\int_{\p \mathcal U_\rho} \di u\ d\sigma_H \ =\ \frac{d}{d\rho}\
\int_{\mathcal U_\rho} \di u\ W\ dg\ .
\end{equation}

This crucial observation allows us to reduce the computation of
the surface integral to that of an integral over the solid region
$\mathcal U_\rho$. Recalling Definition \ref{D:delta}, we have
\begin{equation}\label{ibp7}
\int_{\mathcal U_\rho} \di u\ W\ dg\ =\ \int_{\mathcal U_\rho} X_i
u\ W\ dg\ -\ \int_{\mathcal U_\rho} X_ju\ \pb_j\ \pb_i\ W\ dg\ ,
\end{equation}
where we have adopted the summation convention over repeated indices. Integrating by parts in the first integral in the right-hand side of \eqref{ibp7} we find
\begin{align}\label{ibp8}
\int_{\mathcal U_\rho} X_i u\ W\ dg\ & =\ \int_{\p \mathcal
U_\rho} u <\bN,X_i> \frac{W}{|\bN|}\ dH_{N-1}\ -\ \int_{\mathcal
U_\rho} u\ X_i W\ dg\
\\
& =\ \int_{\p \mathcal U_\rho} u\ \pb_i\ W\ d\sigma_H\ -\
\int_{\mathcal U_\rho} u\ \frac{\pb_j\ X_i p_j}{W}\ W\ dg\ ,
\notag
\end{align}
where we have used $div\ X_i = 0$, see \eqref{divfree}. We next
integrate by parts in the second integral in the right-hand side
of
 \eqref{ibp7}, obtaining as in \eqref{ibp8}
\begin{align}\label{ibp9}
&\int_{\mathcal U_\rho} X_ju\ \pb_j\ \pb_i\ W\ dg\ =\ \int_{\p
\mathcal U_\rho} u\ \pb_j\ \pb_j\ \pb_i\ W\ d\sigma_H
\\
& -\ \int_{\mathcal U_\rho} u\ X_j(\pb_j\ \pb_i\ W)\ dg\
\notag\\
& =\ \int_{\p \mathcal U_\rho} u \ \pb_i\ W\ d\sigma_H\ -\
\int_{\mathcal U_\rho} u\ (X_j \pb_j)\ \pb_i\ W\ dg \notag\\
& - \int_{\mathcal U_\rho} u\ \frac{\pb_j X_j p_i}{W}\ W\ dg\ .
 \notag
\end{align}

Inserting \eqref{ibp8}, \eqref{ibp9} into \eqref{ibp7}, we see
that the boundary integrals disappear and we finally obtain
\begin{equation}\label{ibp10}
\int_{\mathcal U_\rho} \di u\ W\ dg\ =\ \int_{\mathcal U_\rho} u\
X_j(\pb_j)\ \pb_i\ W\ dg\ -\ \int_{\mathcal U_\rho} u\ \boldsymbol
c^{\mathcal S}_i\ W\ dg\ ,
\end{equation}
where we have let
\[
\boldsymbol c^\mS_i\ =\ \frac{\pb_j}{W}\big\{X_i p_j - X_j
p_i\big\}\ .
\]

Formula \eqref{ibp10} is the crucial point in the proof. Proceeding
now as in \eqref{ibp5}, and applying \eqref{ibp6}, we conclude for
every $\rho\in \R$ in a sufficiently small neighborhood of a given
$\rho_0\in \R$
\begin{equation}\label{ibp11}
\int_{\p \mathcal U_\rho} \di u\ d\sigma_H\ =\ \int_{\p \mathcal
U_\rho} u\ X_j(\pb_j)\ \pb_i\  d\sigma_H\ -\ \int_{\p \mathcal
U_\rho} u\ \boldsymbol c^\mS_i\ d\sigma_H\ .
\end{equation}

We now have \[ \boldsymbol c^\mS_i\ =\ \frac{\pb_j}{W}
\big\{X_iX_j \phi - X_j X_i \phi\big\}\ =\  \sum_{s=1}^k \left(
\sum_{j=1}^m b^s_{ij} \pb_j \right) \ob_s\ . \]

Recalling \eqref{equal2} in Proposition \ref{P:equalMC}, we conclude
that \eqref{ibp0} holds. Finally, \eqref{ibp4} follows from the
skew-symmetry of the matrix $\{b^s_{ij}\}_{i,j=1,...,m}$ defined by
\eqref{gc} which gives $\sum_{j=1}^m b^s_{ij} \pb_i \pb_j = 0$ for
every $s=1,...,k$. Hence,
\[
<\cx,\nuX>\  =\  \sum_{s=1}^k \big(\sum_{j=1}^m b^s_{ij} \pb_i
\pb_j\big) \ob_s\ =\ 0\ .
\]

This completes the proof.

\end{proof}

\medskip

\begin{rmrk}\label{R:commutator}
We emphasize that in the Abelian case $\bG = \R^m$, we have $X_i =
\p/\p x_i$, $i = 1,...,m$, and so $[X_i,X_j] = 0$ and thereby
$\boldsymbol c^\mS_{i} \equiv 0$. In this case formula \eqref{ibp1}
 recaptures the classical integration by parts formula on a hypersurface, see for instance \cite{MM}, \cite{Si}.
\end{rmrk}

\medskip

\begin{rmrk}\label{R:commutatorHn}
We note explicitly that when $\bG = \Hn$, the Heisenberg group, then
the horizontal vector field $\cx$ is given by
\begin{equation}\label{cx}
\cx\ =\  \ob\ J(\nuX)\ , \end{equation} where $J:H\Hn \to H\Hn$ is
the symplectic transformation which, in the orthonormal basis
$\{X_1,...,X_{2n}\}$ of $H\Hn$, is represented by the block matrix
\[ J\ =\
\begin{pmatrix}
0 &  I
\\
-I & 0
\end{pmatrix}\ .
\]
We thus obtain from \eqref{cx}
\begin{equation}\label{cxHn}
\cx\ =\  \ob\ (\nuX)^\perp\ =\ \ob\ \big(\pb_{n+1} X_1 + ... +
\pb_{2n} X_n - \pb_1 X_{n+1} - ... - \pb_n X_{2n}\big)\ .
\end{equation}
Therefore, for $\Hn$ formula \eqref{ibp0} reads
\begin{equation}\label{ibpHn}
\int_\mathcal S \nabla^{H,\mS} u\ d \sigma_H\ =\ \int_\mathcal S u\
\bigg\{\mathcal H\ \nuX\ -\ \ob\ \nup\bigg\}d\sigma_H\ .
\end{equation}

In particular, when $n=1$ then in the notation of section
\ref{S:geomid}, see also Remark \ref{R:pq}, we have $\cx = \ob Z$,
and we can write \eqref{ibpHn} as follows
\begin{equation}\label{ibpH1}
\int_\mathcal S \nabla^{H,\mS} u\ d \sigma_H\ =\ \int_\mathcal S u\
\bigg\{\mathcal H\ Y\ -\ \ob\ Z\bigg\}d\sigma_H\ .
\end{equation}
\end{rmrk}

\medskip

We have the following notable consequences of Theorem \ref{T:ibp}.

\medskip

\begin{thrm}\label{T:divH}
Let $\mathcal S\subset \bG$ be a $C^2$ oriented hypersurface, with
characteristic set $\Sigma$. If $\zeta\in C^1_0(\mS \setminus
\Sigma_\mS,HT\mS)$, then we have
\begin{equation}\label{ibp1}
\int_\mathcal S \left\{div_{H,\mS} \zeta\ +\ <\boldsymbol
c^{\mathcal S},\zeta>\right\}\ d \sigma_H\ =\  \int_\mathcal S
\mathcal H\ <\zeta,\nuX>\ d\sigma_H\ ,
\end{equation}
where we have let
\[
div_{H,\mS} \zeta\ =\ \sum_{i=1}^m \di \zeta_i\ .
\]
\end{thrm}

\begin{thrm}\label{T:nonchar}
In a Carnot group $\bG$ suppose that the hypersurface $\mathcal S$
is a vertical cylinder as in Proposition \ref{P:rigate}. If $u\in
C^1_0(\mS)$ we have
\begin{equation}\label{ibp12}
\int_\mathcal S \di u\ d \sigma_H\ =\  \int_\mathcal S u\ \mathcal
H\ \nui\ d\sigma_H\ .
\end{equation}
\end{thrm}

\begin{proof}[\textbf{Proof}]
First of all we notice that the assumption on $\mathcal S$
guarantees that the characteristic set $\Sigma_\mS$ is empty, see
Proposition \ref{P:rigate}. It is thereby legitimate to assume $u\in
C^1_0(\mS)$, instead of $u\in C^1_0(\mS\setminus \Sigma_\mS)$. Next,
we observe that since the defining function of $\mathcal S$ depends
only on the horizontal variables, then the normal has no component
along $V_2$, and therefore $\om_s = 0$ for $s=1,...,k$. This implies
$\cx \equiv 0$, see \eqref{ibp3}. The conclusion thus follows from
\eqref{ibp1}.

\end{proof}

We next establish another integration by parts formula which
involves differentiation along a special combination of the vector
fields $\nuX$ and $T_s$, $s=1,...,k$, where $T_s$ constitute the
orthonormal basis of the first vertical layer defined in \eqref{X}.
Such result plays a central role in the last two sections of this
paper.

\medskip

\begin{thrm}[\textbf{Second sub-Riemannian integration by parts formula}]\label{T:T&Y}
Let $\mS$ be a $C^2$ oriented hypersurface in a Carnot group $\bG$.
For every $f, \zeta \in C^1_0(\mathcal U\setminus \Sigma_\mS)$,
where $\mathcal U\subset \bG$ is an open neighborhood of $\mS$, then
one has for $s = 1,...,k$
\begin{align}\label{oneT&Y}
& \int_\mS f \ (T_s \zeta - \ob_s <\nabla \zeta, \nuX>)\ d\sigma_H\
=\ -\ \int_\mS \zeta \ (T_s f - \ob_s <\nabla f, \nuX>)\ d\sigma_H
\\
& +\ \int_\mS f\ \zeta \left\{\ob_s\ \mathcal H \ +\  \sum_{\ell =
1}^{m_3}\sum_{i=1}^m b^\ell_{is} \pb_i \ob_{3,\ell}\right\}\
d\sigma_H\ , \notag
\end{align} where $\ob_s$ are defined in
\eqref{omegas}, and for $i=1,...,m$, $s = 1,...,k$ and $\ell =
1,...,m_3$, we have let $b^\ell_{is} =
<[e_i,\epsilon_s],e_{3,\ell}>$. In particular, for a hypersurface
$\mS \subset \Hn$, we have $k=1$, and therefore letting $\ob_1 =
\ob$ and setting $Yf \overset{def}{=} <\nabla f,\nuX>$, see also
section \ref{S:geomid}, we obtain,
\begin{equation}\label{twoT&Y}
\int_\mS f \ (T  - \ob Y)\zeta\ d\sigma_H\ =\ -\ \int_\mS \zeta \ (T
 - \ob Y)f\ d\sigma_H\ +\ \int_\mS f \zeta \ob\ \mathcal H\
d\sigma_H\ .
\end{equation}
\end{thrm}

\begin{proof}[\textbf{Proof}]
We use the same idea of the proof of Theorem \ref{T:ibp}, except
that this time we consider
\begin{align*}
& \int_{\mathcal U_\rho} [T_s f - <\nabla(\ob_s f),\nuX>] W dg\ =\
\int_{\p \mathcal U_\rho} <T_s,\bN> f \frac{W}{|\bN|}\ d\sigma\ -\
\int_{\mathcal U_\rho} f\ div(W T_s)\ dg
\\
& -\ \int_{\p \mathcal U_\rho} <\nuX,\bN> \ob_s f \frac{W}{|\bN|}\
d\sigma\ +\ \int_{\mathcal U_\rho} \ob_s f\ div(W\nuX)\ dg
\\
& =\  \int_{\p \mathcal U_\rho} \omega_s f\ d\sigma_H\ -\
\int_{\mathcal U_\rho} f\ \frac{T_sW}{W}\ Wdg \\
& -\ \int_{\p \mathcal U_\rho} W \ob_s f\ d\sigma_H\ +\
\int_{\mathcal U_\rho} \ob_s f\ W\ div(\nuX)\ dg \ +\ \int_{\mathcal
U_\rho} \ob_s f\ \frac{<\nabla W,\nuX>}{W}\ Wdg\ ,
\end{align*}
where we have used the identity $<\nuX,\bN> = W$, see \eqref{ip}.
Since $\ob_s W = \omega_s$, the two boundary terms drop and we are
left with
\begin{align*}
& \int_{\mathcal U_\rho} [T_s f - <\nabla(\ob_s f),\nuX>] W dg\ =\
-\
\int_{\mathcal U_\rho} f\ \frac{T_sW}{W}\ Wdg \\
&  +\ \int_{\mathcal U_\rho} \ob_s f\ W\ div(\nuX)\ dg \ +\
\int_{\mathcal U_\rho} \ob_s f\ \frac{<\nabla W,\nuX>}{W}\ Wdg\ .
\end{align*}

Using the coarea formula as in the proof of Theorem \ref{T:ibp}, and
differentiating the resulting integrals, we obtain from the latter
identity
\begin{align}\label{ibpI}
& \int_\mS [T_s f - <\nabla(\ob_s f),\nuX>] d\sigma_H\ =\ -\
\int_{\mS} f\ \frac{T_sW}{W}\ d\sigma_H \\
&  +\ \int_\mS \ob_s f\ div(\nuX)\ d\sigma_H \ +\ \int_\mS \ob_s f\
\frac{<\nabla W,\nuX>}{W}\ d\sigma_H\ . \notag
\end{align}

Since Definition \ref{D:HMC} and Proposition \ref{P:equalMC} give
\[
\mathcal H\ =\ div_{H,\mS}(\nuX)\ =\ \sum_{i=1}^m X_i \pb_i\ =\
div(\nuX)\ ,
\]
we can re-write \eqref{ibpI} as follows
\begin{align}\label{yummy}
& \int_\mS [T_s f - \ob_s <\nabla f,\nuX>] d\sigma_H\ =\
\int_{\mS} f\ <\nabla \ob_s, \nuX>  d\sigma_H \\
&  +\ \int_\mS \ob_s f\ \mathcal H\ d\sigma_H \ -\ \int_\mS f\
\left(\frac{T_sW}{W} - \ob_s  \frac{<\nabla W,\nuX>}{W}\right)\
d\sigma_H\ . \notag
\end{align}

We now observe that \eqref{omegas} gives
\begin{equation}\label{Yobs}
<\nabla \ob_s,\nuX>\ =\ \frac{<\nabla \omega_s,\nuX>}{W}\ -\ \ob_s\
\frac{<\nabla W,\nuX>}{W}\ .
\end{equation}

On the other hand, we have
\begin{equation}\label{Yoms}
\frac{<\nabla \omega_s,\nuX>}{W}\ =\ \frac{T_s W}{W}\ +\ \sum_{\ell
= 1}^{m_3}\sum_{i=1}^m b^\ell_{is}\ \pb_i\ \ob_{3,\ell}\ .
\end{equation}

To prove \eqref{Yoms}, suppose, as we may, that $\mS$ is locally
described as the zero set of a $C^2$ function $\phi$, and that $\bN
= \nabla \phi = \sum_{i=1}^m p_i X_i + \sum_{s=1}^k \omega_s T_s +
\sum_{j=3}^r \sum_{\ell = 1}^{m_j} \omega_{j,\ell} X_{j,\ell}$. We
thus have
\begin{align*}
<\nabla \omega_s,\nuX>\ & =\ <\nabla(T_s\phi),\nuX>\  =\
\sum_{i=1}^m \pb_i X_i(T_s\phi)
\\
& =\ \sum_{i=1}^m \pb_i T_s(X_i\phi)\ +\ W\ \sum_{i=1}^m \sum_{\ell
= 1}^{m_3} b^\ell_{is} \pb_i \ob_{3,\ell}
\\
& =\ \sum_{i=1}^m \pb_i T_s(\pb_i W)\ \ +\ W\ \sum_{\ell =
1}^{m_3}\sum_{i=1}^m b^\ell_{is} \pb_i \ob_{3,\ell}
\\& =\ \left(\sum_{i=1}^m \pb_i^2\right) T_s W +  \left(\sum_{i=1}^m \pb_i
T_s\pb_i\right) W \ +\ W\ \sum_{\ell =
1}^{m_3}\sum_{i=1}^m b^\ell_{is} \pb_i \ob_{3,\ell}\\
& =\ T_sW\ +\ W\ \sum_{\ell = 1}^{m_3}\sum_{i=1}^m b^\ell_{is} \pb_i
\ob_{3,\ell}\ ,
\end{align*}
where we have used the commutation relations $[X_i,T_s] = \sum_{\ell
= 1}^{m_3} b^\ell_{is} X_{3,\ell}$. This proves \eqref{Yoms}.
Inserting \eqref{Yoms} in \eqref{Yobs}, and the resulting equation
in \eqref{yummy}, we reach the conclusion
\[
\int_\mS (T_s f - \ob_s <\nabla f,\nuX>)\ d\sigma_H\ =\  \int_\mS f
\ \ob_s\ \mathcal H\ d\sigma_H\ +\ \int_\mS f\  \sum_{\ell =
1}^{m_3}\sum_{i=1}^m b^\ell_{is} \pb_i \ob_{3,\ell}\ d\sigma_H\ .
\]

If we replace $f$ by $f \zeta$ in the latter integral identity we
obtain the sought for integration by parts formula.

\end{proof}

\vskip 0.6in


\section{\textbf{Tangential horizontal Laplacian}}\label{S:lb}
\vskip 0.2in

In this section we introduce a tangential partial differential
operator, $\delhs$ (and a modified version of the latter), which
constitutes the sub-Riemannian counterpart of the classical
Laplace-Beltrami operator on a hypersuface. In fact, as we will see,
it reduces to the latter when the group $\bG$ is Abelian. It has
however one aspect which distinguishes it from its classical
predecessor, and this is lack of self-adjointness in
$L^2(\mS,d\sigma_H)$. This phenomenon is caused by the presence of
the ``drift" term $\cx$ in the integration by parts formula in
Theorem \ref{T:ibp}. In the next section we will show that the
horizontal mean curvature flow recently proposed by Bonk and Capogna
\cite{BC} satisfies a nonlinear pde which involves the operator
$\delhs$, see Theorem \ref{T:mcf}.

\medskip

\begin{dfn}\label{D:lb}
Given a function $u\in C^2(\mS)$, the \emph{tangential horizontal
Laplacian} of $u$ on $\mathcal S$ is defined as follows at points of
$\mS\setminus \Sigma_\mS$
\begin{equation}\label{horlap}
\Delta_{H,\mS} u\ \overset{def}{=}\ \sum_{i=1}^m \di \di u\ .
\end{equation}
We also introduce the \emph{modified tangential horizontal
Laplacian} on $\mS$
\begin{equation}\label{lap2}
\lhs u\ \overset{def}{=}\ \delhs u\ +\ <\boldsymbol c^{\mathcal S},
\del u>\ ,
\end{equation}
where $\cx$ is given by \eqref{ibp3}.
\end{dfn}

\medskip

\begin{rmrk}\label{R:h}
One should keep in mind that when $\mathcal S$ is a vertical
cylinder given by \eqref{vertical}, then the operators $\delhs$ and
$\lhs$ coincide
\[
\lhs\ =\ \delhs \ .
\]
In such case it is easy to show from Theorem \ref{T:ibp} that
$\delhs$ is formally self-adjoint in $L^2(\mathcal S,d\sx)$.
\end{rmrk}

\medskip

One basic raison d'\^etre for the operator $\lhs$ is in the
following sub-Riemannian Stokes' theorem which follows from
Theorem \ref{T:ibp}.

\medskip

\begin{cor}\label{C:lapX}
Let $u\in C^2_0(\mS\setminus \Sigma_\mS)$, then we have
\begin{equation}\label{lap1}
\int_\mathcal S \lhs u\ d\sigma_H\ =\ 0\ .
\end{equation}
\end{cor}

\begin{proof}[\textbf{Proof}]
It suffices to take $\di u$ instead of $u$ in Theorem \ref{T:ibp},
and then add the resulting identities in $i=1,...,m$. Keeping in
mind the definition \eqref{lap2}, formula \eqref{ibp1} gives,
\[
\int_\mathcal S \lhs u\ d\sigma_H\ =\  \int_\mathcal S \mathcal H
<\del u, \nuX> d\sigma_H\ =\ 0\ ,
\]
since by \eqref{P:rigate} one has $<\del u, \nuX> = 0$ on $\mathcal
S\setminus \Sigma_\mS$.

\end{proof}

\medskip

\begin{cor}\label{C:qf}
Let $u\in C^1(\mS)$, then for every $\zeta \in C^2_0(\mS\setminus
\Sigma_\mS)$ we have
\begin{equation}\label{qf1}
\int_\mathcal S <\del u,\del \zeta> d\sigma_H\ =\ -\ \int_\mathcal S
u\ \lhs \zeta\ d\sigma_H\ .
\end{equation}

\end{cor}

\begin{proof}[\textbf{Proof}]
We take $u\ \di \zeta$, instead of $u$, in Theorem \ref{T:ibp}.

\end{proof}

\medskip

\begin{rmrk}\label{R:productgroup2}
In connection with Remark \ref{R:productgroup1} we see that for
the product group $\hat{\bG} = \bG \times \R$ one has $\lhs=
\delhs = \delh$ on $\mathcal S = \bG \times \{0\}$.
\end{rmrk}

\medskip

The following formulas are verified by direct computation from the
definition.

\medskip

\begin{lemma}\label{L:comp}
Let $u , v \in C^2(\mathcal O)$, $F\in C^2(\R)$, then we have on
$\mathcal O\setminus \Sigma$
\begin{equation}\label{product}
\lhs (u v)\ =\ u\ \lhs v\ +\ v\ \lhs u\ +\ 2\ <\del u, \del v>\ ,
\end{equation}
\begin{equation}\label{comp}
\lhs (F \circ u)\ =\ (F''\circ u)\ |\del u|^2\ +\ (F' \circ u)\ \lhs
u\ .
\end{equation}
\end{lemma}

\medskip

The next result provides a useful mean for computing the operators
$\delhs$ and $\lhs$ on $\mathcal S$, using the vector fields
$X_1,...,X_m$ in the ambient group $\bG$.

\medskip

\begin{prop}\label{P:deltau}
Let $u\in C^2(\mS)$, then we have on $\mS\setminus \Sigma$
\begin{equation}\label{tsl}
\delhs u\ =\ \delh \ou\ -\ <\nabla^2_H\ou\ \nuX , \nuX>\ -\ <\nabh
\ou,\nuX>\ \mathcal H\ ,
\end{equation}
\begin{equation}\label{mtsl}
\lhs u\ =\ \delh \ou\ +\ <\cx,\nabh \ou>\ -\ <\nabla^2_H\ou\ \nuX ,
\nuX>\ -\
 <\nabh \ou,\nuX>\ \mathcal H\ ,
\end{equation}
where $\ou$ denotes any extension of $u$. In the above formulas, the
notation $\nabla_H^2 \ou$ indicates the horizontal Hessian of $\ou$
introduced in \eqref{hessian}.
\end{prop}

\begin{proof}[\textbf{Proof}]
We begin with Definition \ref{D:delta} which gives
\[
\del u\ =\ \nabh \ou\ -\ <\nabh \ou,\nuX> \nuX\ ,
\]
where $\ou$ is any extension of $u$. Applying \eqref{horlap}, and
using the summation convention over repeated indices, we find
\begin{equation}\label{u1}
\di \di u\ =\ \di (X_i \ou)\ -\ \di(<\nabh \ou,\nuX>)\ \nui\ -\
<\nabh \ou,\nuX> \di\ \nui\ .
\end{equation}

We now compute the terms in the right-hand side of \eqref{u1}.
\begin{align}\label{u2}
\di(X_i \ou)\ &=\ X_i X_i \ou\ -\ <\nabh(X_i \ou),\nuX> \nui
\\
& =\ \delh \ou\ -\ X_j X_i \ou\ \nui\ \nuj
\notag\\
& =\ \delh \ou\ -\ <\nabla^2_H \ou\ \nuX, \nuX>\ . \notag
\end{align}

Next, equation \eqref{normone} gives
\begin{align}\label{u3}
\di (<\nabh \ou,\nuX>)\ \nui\ &=\ X_i (X_j \ou\ \nuj)\ \nui\ -\
<\nabh(<\nabh \ou,\nuX>),\nuX>\ \nui\ \nui
\\
& =\ 0\ .
\notag
\end{align}

Finally, we find from Proposition \ref{P:equalMC}
\begin{equation}\label{u4}
<\nabh \ou,\nuX> \di\ \nui\ =\  <\nabh \ou,\nuX>\ \mathcal H\ .
\end{equation}

We now substitute \eqref{u2}, \eqref{u3}, \eqref{u4} in
\eqref{u1}. To reach the desired conclusion we only need to
observe that thanks to \eqref{ibp3} one has
\[
<\cx, \del \ou>\ =\ <\cx, \nabh \ou>\ -\ <\nabh \ou,\nuX>
<\cx,\nuX>\ =\ <\cx, \nabh \ou>\ .
\]

\end{proof}

\medskip

The first elementary example of solutions of the tangential
operators $\delhs$ and $\lhs$ is provided by the following
consequence of Proposition \ref{P:deltau}.

\medskip

\begin{prop}\label{P:lconst}
If the function $u$ is constant on $\mathcal S$, then
\[
\delhs u\ =\ \lhs u\ =\ 0\ .
\]
\end{prop}

\begin{proof}[\textbf{Proof}]
First of all, let us notice that, since $\delhs u$ and $\lhs u$ only
depend on the values of $u$ on $\mS$, we can without restriction
assume that $\ou \equiv 1$ in $\bG$. Under such hypothesis the
conclusion now follows trivially from Proposition \ref{P:deltau}.

\end{proof}

\medskip

Another interesting consequence of Proposition \ref{P:deltau} and of the grading structure of a Carnot group is
the following.

\medskip
\begin{prop}\label{P:harmonicity}
Let $\mathcal S\subset \bG$ be a $H$-minimal hypersurface, then if $x(g) =
(x_1(g),...,x_m(g))$ denote the projection onto the horizontal
layer of the exponential coordinates of $g\in \bG$ (see
\eqref{coord2}), one has
\[
\delhs(x_i)\ =\ 0\ ,\quad\quad\quad\quad i = 1,...,m\ .
\]
\end{prop}

\begin{proof}[\textbf{Proof}]
From Proposition \ref{P:coord} we have $\delh(x_i) = 0$, and also
$\nabla^2_H(x_i) = 0$. The desired conclusion thus follows
immediately from \eqref{tsl}.

\end{proof}

\medskip

We next analyze a situation of special interest, namely when $\bG$
is a Carnot group of step $r=2$, and one has a hypersurface
$\mathcal S$ given as a graph over the first layer of the Lie
algebra. In such case, identifying via the exponential map $g =
\exp\ \xi(g)$ with $\xi(g) \cong (x(g),t(g))$, we can find an open
set $\Om \subset V_1$, and a $C^2$ function $h:\Om \to \mathbb R$,
such that for some $s\in \{1,...,k\}$, $\mathcal S$ can be written
as
\begin{equation}\label{graphS}
\mathcal S\ =\ \{(x(g),t(g))\in \bG \mid x(g)\in \Om\ ,\ t_s(g) =
h((x(g))\}\ .
\end{equation}

For instance, in the special case of the Heisenberg group $\Hn$ we
would be considering a graph over $\mathbb R^{2n}$, i.e., $\mathcal
S = \{ (x,y,t)\in \Hn \mid (z,y)\in \Om \subset \mathbb R^{2n}\ ,\ t
= h(x,y)\}$.

\medskip

\begin{thrm}\label{T:graph}
Let $\bG$ be a Carnot group of step $r=2$, and $\mathcal S\subset
\bG$ be a $H$-minimal hypersurface of the type \eqref{graphS}, then
outside the characteristic set $\Sigma_\mS$ the coordinate functions
$x_1,...,x_m, t_1,...,t_k$ are solutions of the tangential
sub-Laplacian on $\mathcal S$.
\end{thrm}

\begin{proof}[\textbf{Proof}]
For the horizontal coordinates $x_1,...,x_m$ the conclusion
follows from Proposition \ref{P:harmonicity}. We now recall
\eqref{coord3} in Proposition \ref{P:coord}
\begin{equation}\label{graph2}
X_i t_s\ =\ \frac{1}{2}\ <[\xi_1, e_i],\epsilon_s>\ ,\
\quad\quad\quad \delh t_s\ =\ 0\ ,\quad\quad\quad  i = 1, ... , m\
,\ s = 1, ... , k\ .
\end{equation}

The first equation in \eqref{graph2} can be written
\[
X_i t_s\ =\ \frac{1}{2}\ \sum_{j=1}^m x_j <[e_j, e_i],\epsilon_s>\
.
\]

Thanks to \eqref{coord1-1} this gives
\[
X_j X_i t_s\ =\ \frac{1}{2}\ <[e_j,e_i],\epsilon_s>\ =\ -\
\frac{1}{2}\ <[e_i,e_j],\epsilon_s>\ =\ -\ X_i X_j t_s\ ,
\]
and therefore for every $s=1,...,k$, one has
\begin{equation}\label{Hessys}
\nabla^2_H(t_s)\ =\ 0\ .
\end{equation}

Now Proposition \ref{P:deltau} gives for any $l\in \{1,...,k\}$,
with $l\not= s$
\[
\delhs t_l\  =\ \delh t_l\ -\ <\nabla^2_H(t_l)\ \nuX , \nuX>\ -\
<\nabh h,\nuX>\ \mathcal H\ =\ 0\ ,
\]
when $\mathcal S$ is $H$-minimal, thanks to \eqref{graph2} and
\eqref{Hessys}. We are left with proving that, if $\mathcal S$ is
$H$-minimal then $\delhs(t_s) = 0$. Since on $\mathcal S$ we have
$t_s = h(x)$, we need to show that $\delhs h = 0$ on $\mathcal S$.
With this objective in mind we begin by expressing the $H$-mean
curvature of $\mathcal S$ in terms of the function $h$. We consider
the function $\phi(g) = t_s - h(x)$ defining $\mathcal S$. According
to Proposition \ref{P:XMCG}, one has on $\mathcal S\setminus
\Sigma_\mS$,
\begin{equation}\label{graph1}
\mathcal H =\ \frac{1}{|\nabh\phi|^3} \left\{|\nabh\phi|^2 \delh
\phi\ -\ \delinf \phi\right\}\ .
\end{equation}

On the other hand, we have from Proposition \ref{P:deltau}
\begin{align}\label{tslh}
\delhs h\ & =\ \delh h\ -\ <\nabla^2_Hh\ \nuX , \nuX>\ -\  <\nabh h,\nuX>\ \mathcal H \\
& =\ \frac{1}{|\nabh\phi|^2} \left\{|\nabh\phi|^2 \delh h\ -\
<\nabla^2_H
h\ \nabh\phi, \nabh\phi>\right\} \ -\  <\nabh h,\nuX>\ \mathcal H \notag\\
& =\ \frac{1}{|\nabh\phi|^2} \big\{|\nabh\phi|^2 \delh(h-t_s) +
\delh t_s\ -\ <\nabla^2_H (h-t_s)\ \nabh\phi, \nabh\phi> \notag\\
&  -\  <\nabla^2_H (t_s)\ \nabh\phi, \nabh\phi> \big\} \ -\ <\nabh
h,\nuX>\ \mathcal H\ .\notag
\end{align}

If  in \eqref{tslh} we use \eqref{Hessys} and the second equation
in \eqref{graph2}, we obtain
\begin{equation}\label{graph3}
\delhs h\  =\ -\ \frac{1}{|\nabh\phi|^2} \left\{|\nabh\phi|^2
\delh \phi\ -\ \delinf \phi\right\}\ -\  <\nabh h,\nuX>\ \mathcal
H\ .
\end{equation}

We now compare \eqref{graph3} with \eqref{graph1} to reach the
following interesting conclusion
\begin{equation}\label{tslh2}
\delhs h\ =\ -\  \{|\nabh\phi|\ +\ <\nabh h,\nu_H>\}\ \mathcal H\
.
\end{equation}

It is now clear from \eqref{tslh2} that if $\mathcal H\equiv 0$,
then $\delhs h = 0$, and this completes the proof.

\end{proof}

\medskip

\begin{cor}\label{C:graphs}
In the Heisenberg group let \[ \mathcal S\ =\ \{(x,y,t)\in \Hn\mid
(x,y)\in \Om\ ,\ t = h(x,y)\}\ , \] where $\Om \subset \mathbb
R^{2n}$ is an open set, and $h\in C^2(\Om)$. Is $\mathcal S$ is
$H$-minimal, then the coordinate functions
$x_1,...,x_n,y_1,...,y_n,t$ are solutions of $\delhs$ on $\mathcal
S$.
\end{cor}

\medskip

\begin{cor}\label{C:laphorvar}
Let $\bG$ be a Carnot group, and consider the exponential
horizontal coordinates $x_1(g),...,x_m(g)$ in $\bG$, then
\[
\delhs(x_i)\ =\ -\ <\nuX,X_i>\ \mathcal H\ =\ -\ \pb_i\ \mathcal
H\ ,\quad\quad\quad i = 1,...,m\ .
\]
Consider the exponential coordinates $t_1(g),...,t_k(g)$ in the
first vertical layer $V_2$, then
\[
\delhs(t_s)\ =\ -\ \frac{1}{2}\ \sum_{i,j=1}^m b^s_{ij}\ x_i\
\pb_j\ \mathcal H\ ,\quad\quad\quad s = 1,...,k\ .
\]
In particular, when $\bG = \HH$, then
\[
\delhs(t)\ =\  -\ \frac{1}{2}\ (x \pb + y \qb)\ .
\]
\end{cor}

\medskip

We close this section with introducing the notions of $p$-Dirichlet
integral and of $p$-harmonic function on an hypersurface. Such
notions play a central role in the development of geometric
subelliptic pde's on hypersurfaces in Carnot groups.

\medskip

\begin{dfn}\label{D:DI}
Suppose that $\mS\subset \bG$ be a $C^2$ hypersurface, with
$\Sigma_\mS = \varnothing$. Given $1<p<\infty$ we define the
$p$-\emph{Dirichlet integral} of a function $u\in C^1_0(\mS)$ as
\[
\mathcal E_{H,\mathcal S}(u)\ =\ \frac{1}{p}\ \int_{\mathcal S}
|\del u|^p\ d\sigma_H\ .
\]
Suppose that $u\in L^p(\mS,d\sigma_H)$, and that moreover $\di u\in
L^p(\mS,d\sigma_H)$, for $i=1,...,m$. We say that $u$ is
$p$-\emph{subharmonic} (-\emph{superharmonic}) in $\mathcal S$ if
for every $\zeta\in C^1_0(\mS)$, $\zeta \geq 0$, one has
\[
\int_\mathcal S\ |\del  u|^{p-2}\ <\del u,\del \zeta> d\sigma_H\
\leq\ 0\ (\geq\ 0)\ .
\]
We say that $u$ is $p$-\emph{harmonic} in $\mathcal S$ if $u$ is
simultaneously $p$-subharmonic and $p$-superharmonic. When $p=2$
we simply say that $u$ is \emph{subharmonic, superharmonic} or
\emph{harmonic} in $\mathcal S$.
\end{dfn}

\medskip

According to Corollary \ref{C:qf} we can adopt the following
alternative notion of subharmonicity.

\medskip

\begin{dfn}\label{D:sh}
A function $u\in L^1_{loc}(\mathcal S,d\sigma_H)$ is called
subharmonic in $\mathcal S$ if
\begin{equation}\label{sh}
0\ \leq\ \int_\mathcal S u\ \lhs \zeta\ d\sigma_H\ ,\quad\quad\quad
\text{for every}\quad \zeta \in C^2_0(\mS)\ ,\ \zeta \geq 0\ .
\end{equation}
\end{dfn}

\vskip 0.6in


\section{\textbf{Flow by horizontal mean curvature}}{\label{S:MCF}

\vskip 0.2in

In connection with Proposition \ref{P:deltau}, we recall the
Riemannian counterpart of \eqref{tsl}
\begin{equation}\label{lb}
\Delta_M u\ =\ \Delta u\ -\ <\nabla^2 u\ \n,\n> \  -\ (n-1)
<\nabla_M u,\n>\  H\ ,
\end{equation}
where $\Delta_M$ and $\nabla_M$ respectively represent the
Laplace-Beltrami operator and the intrinsic gradient on an
$(n-1)$-dimensional Riemannian manifold $M$. Formula \eqref{lb}
plays a crucial role, for instance, in the derivation of the
equation for flow by mean curvature, see for instance \cite{E}. If
one considers a family of smooth embeddings $F(\cdot,t): M\to
\Rn$, then with $M_t = F(M,t)$, the equation of flow my mean
curvature is given by
\begin{equation}\label{mcf}
\frac{\p F}{\p t}(p,t)\ =\ -\ (n-1)\ H\ \n\ .
\end{equation}

If we write $\boldsymbol x = F(p,t)$, then using \eqref{lb} and
\eqref{mcf} we obtain the nonlinear partial differential equation
\begin{equation}\label{mcfpde}
\frac{\p \boldsymbol x}{\p t}\ =\ \Delta_{M_t} \boldsymbol x\ ,
\end{equation}
which is satisfied by the components $(x_1,...,x_n)$ of
$\boldsymbol x$. This can be readily recognized as follows.
Equation \eqref{lb} gives for each component $x_i$
\begin{align*}
\Delta_{M_t} (x_i)\ & =\ \Delta(x_i)\ -\ <\nabla^2(x_i)\ \n,\n>\
-\ (n-1)\ <\nabla(x_i),\n>\ H
\\
& =\ -\ (n-1)\ <e_i,\n>\ H\ =\ -\ (n-1)\ \n_i\ H\ .
\end{align*}

In other words, we have $\Delta_{M_t} \boldsymbol x = - (n-1) H
\n$. This equation, combined with \eqref{mcf}, proves
\eqref{mcf2}.

We next want to prove a sub-Riemannian analogue of \eqref{mcfpde}
for the mean curvature flow in the Heisenberg group recently
proposed by Bonk and Capogna in \cite{BC}. We consider a smooth
hypersurface in a Carnot group $\mS\subset \bG$, and a family of
smooth embeddings $F : \mS \times (0,T) \to \bG$. We will denote
by $S^\la = F(\mS,\la)$. The reader should note that we are using
the unconventional parameter $\la\in (0,T)$ to indicate time. The
reason is due to the fact that, to keep a homogeneous notation
with the Heisenberg group, we have already reserved the letter $t
= (t_1,...,t_k)$ to indicate the exponential coordinates in the
first vertical layer $V_2$ of the Lie algebra of $\bG$, see
\eqref{coord2}. In \cite{BC} the authors have introduced the
following definition of \emph{horizontal mean curvature flow} when
the group $\bG$ is $\Hn$. At any point $F(g,\la)\in \mS^\la
\setminus \Sigma^\la$ ($\Sigma^\la$ denotes the characteristic set
of $\mS^\la$), they require that
\begin{equation}\label{bonkcap}
<\frac{\p F}{\p \la},\bN>\ =\ -\ \mathcal H\ <\nuX,\bN>\  .
\end{equation}

We notice that it is important to project the flow along the normal
direction since the vector equation $\frac{\p F}{\p \la} = -
\mathcal H \nuX$ is meaningless: the right-hand side evolves in the
horizontal bundle $H\bG$, whereas the left-hand side has components
which move outside of it. Also, as noted in \cite{BC}, ``any
tangential component of the velocity field only gives rise to a
re-parametrization of the surface with no effect on the geometric
evolution". At characteristic points the equation \eqref{bonkcap} is
not defined and the way the authors circumvent this obstacle is by
restricting to $\mS$ the Riemannian $\ep$-regularization of the
sub-Riemannian metric of $\mS$ introduced in \eqref{gijep}. We refer
the reader to \cite{BC} for the relevant details. We want to next
prove the following result which underscores the interest of the
operator $\delhs$ introduced in the previous section. It should be
thought of as the sub-Riemannian analogue of \eqref{mcfpde}.

\medskip

\begin{thrm}\label{T:mcf}
Let $F : \mS \times (0,T) \to \bG$ be a $C^2$ solution of the
horizontal mean curvature flow \eqref{bonkcap}, then at any
non-characteristic point $F(g,\la)\in \mS^\la$ one has
\begin{equation}\label{geombc}
<\frac{\p F}{\p \la},\bN>\ =\ <\dla F,\bN>\ ,
\end{equation}
where the latter equation must be interpreted component-wise.
\end{thrm}

\begin{proof}[\textbf{Proof}]
To make our proof as transparent as possible we discuss in detail
the case of the first Heisenberg group $\HH$. The details of the
more general case, as well as some applications of \eqref{geombc},
will appear elsewhere. We consider $F(g,\la) =
(x(g,\la),y(g,\la),t(g,\la))$ and notice that we have from
\eqref{coord} below,
\begin{align}\label{mcf1}
\dla F\ & =\ (\dla(x),\dla(y),\dla(t))
\\
& =\ \dla(x) X_1 + \dla(y) X_2 + \left(\dla(t) + \frac{y \dla(x) -
x \dla(y)}{2}\right) T \notag
\end{align}

At this point we use \eqref{mcf1} and the fact that \[ \bN\ =\
\left(\pb\ X_1\ +\ \qb\ X_2\ +\ \frac{<\bN,T>}{W}\ T\right)\ W\ ,
\]
to discover that
\begin{align}\label{mcf2}
<\dla F,\bN>\ & =\ W\ \bigg\{\pb\ \dla(x)\ +\ \qb\ \dla(y)
\\
& +\ \frac{<\bN,T>}{W} \left(\dla(t) + \frac{y \dla(x) - x
\dla(y)}{2}\right)\bigg\}\ . \notag
\end{align}

We now use Corollary \ref{C:laphorvar}, which in the present
situation gives,
\begin{equation}\label{mcf3}
\dla(x)\ =\ -\ \pb\ \mathcal H\ ,\quad \dla(y)\ =\ -\ \qb\
\mathcal H\ ,\quad \dla(t)\ =\ -\ \frac{x \qb - y \pb}{2}\
\mathcal H\ .
\end{equation}

Substituting \eqref{mcf3} in \eqref{mcf2} we obtain the remarkable
conclusion
\begin{align}\label{mcf4}
<\dla F,\bN>\ & =\ W \bigg\{- \pb^2 \mathcal H - \qb^2 \mathcal
H \\
& + \ \frac{<\bN,T>}{W}\left(- \frac{x \qb - y \pb}{2} \mathcal H
+
\frac{x \qb - y \pb}{2} \mathcal H\right)\bigg\} \notag\\
& =\ -\ W\ \mathcal H\ . \notag
\end{align}

On the other hand, \eqref{ip} gives
\[
 <\nuX,\bN>\ \mathcal H\ =\ W\ \mathcal H\ .
 \]

 Combining the latter equation with \eqref{mcf4} we reach the
 conclusion
 \begin{equation}\label{mcf5}
 <\dla F,\bN>\ =\  -\  \mathcal H\ <\nuX,\bN>\ .
 \end{equation}

Finally, from \eqref{mcf5} and \eqref{bonkcap} we obtain
\eqref{geombc}

\end{proof}

\vskip 0.6in


\section{\textbf{Some geometric identities in the Heisenberg group}}\label{S:geomid}

\vskip 0.2in

In this section we collect several geometric identities in the
Heisenberg group $\HH$ which, besides their intrinsic interest,
play an important role in the development of the first and second
variation formulas in Section \ref{S:1&2var}. We note
preliminarily that
\begin{equation}\label{wedges} X_1 \wedge X_2\ =\ T\
,\quad\quad\quad X_2 \wedge T\ =\ X_1\ ,\quad\quad\quad X_1 \wedge
T\ =\ -\ X_2\ , \end{equation} where the wedge products are
computed with respect to the left-invariant Riemannian metric with
respect to which $\{X_1,X_2,T\}$ constitute an orthonormal basis.
We also observe that the passage from the orthonormal basis
$\{X_1,X_2,T\}$ to the standard rectangular coordinates of
$\mathbb R^3$ is given by the formula
\begin{equation}\label{coord}
a X_1\ +\ b X_2\ +\ c T\ =\ \left(a, b, c\ +\ \frac{bx
-ay}{2}\right)\ .
\end{equation}

Throughout this section $\mS\subset \HH$ denotes an oriented $C^2$
surface, with non-unit normal $\bN$, and Riemannian Gauss map
$\boldsymbol \nu$, we consider the functions $p_1, p_2$ and $W$ on
$\mathcal S$ defined in \eqref{ps}. As we have mentioned in Remark
\ref{R:pq}, for computational ease it will be convenient to adopt in
this and the next two sections the slightly different notation $p =
p_1, q = p_2$, i.e.,
\begin{equation}\label{pq}
p\ =\ <\bN, X_1>\ ,\ \quad q\ =\ <\bN , X_2>\ , \ \quad W\ =\
\sqrt{p^2 + q^2}\ .
\end{equation}

The horizontal Gauss map defined in \eqref{up2} is now given on
$\mathcal S\setminus \Sigma$ by
\begin{equation}\label{gm}
\nuX\ =\  \pb\ X_1\ +\ \qb\ X_2\ ,
\end{equation}
where we have let
\begin{equation}\label{pbar}
\pb\ =\ \frac{p}{W}\ ,\quad\quad\quad \qb\ =\ \frac{q}{W}\ ,
\quad\quad\text{so that}\quad\quad \pb^2\ +\ \qb^2\ \equiv\
1\quad\quad \text{on}\quad\quad \mathcal S \setminus \Sigma_\mS\ .
\end{equation}

We also introduce the notation
\begin{equation}\label{tc} \omega\
=\ <\bN,T>\,\quad\quad\quad \ob \ =\ \frac{\omega}{W}\ .
\end{equation}

We notice explicitly that if $\mS\in C^k$, $k\geq 2$, then $p, q,
\omega, \pb, \qb, \ob\in C^{k-1}(\mS\setminus \Sigma_\mS)$. We also
note that, thanks to Proposition \ref{P:equalMC}, the $H$-mean
curvature of $\mathcal S$ is presently given by the formula
\begin{equation}\label{HH}
\mathcal H\ =\ X_1\pb\ +\ X_2 \qb\ .
\end{equation}

Along with the horizontal Gauss map $\nuX$ we consider the vector
field
\begin{equation}\label{nup}
\nup\ =\ \begin{pmatrix} 0 & 1 \\ - 1 & 0 \end{pmatrix} \nuX\ =\ \qb
X_1 - \pb X_2\ , \end{equation} which, as already noticed in Section
\ref{S:hmc}, constitutes a basis of $HT\mS$. It will be convenient
to keep a different notation for the action of the vector fields
$\nuX$, $\nup$ on a function $\zeta\in C^1_0(\mathcal S \setminus
\Sigma_\mS)$. We thus set
\begin{equation}\label{Y}
Y\zeta\ \overset{def}{=}\ <\nabla \zeta, \nuX>\ =\ \pb\ X_1 \zeta\
+\ \qb\ X_2 \zeta\ ,
\end{equation}
\begin{equation}\label{Z}
Z\zeta\ \overset{def}{=}\ <\nabla\zeta, \nup>\ =\ \qb\ X_1\zeta\
-\ \pb\ X_2\zeta\ .
\end{equation}

We mention that in the right-hand sides of \eqref{Y}, \eqref{Z} the
vector fields $X_1$, $X_2$ act on an extension $\overline \zeta$ of
$\zeta$. However, for the sake of simplifying the notation we have
used, and will continue to do so below, the same notation for both
functions. It is worth observing that $\{Z, Y, T\}$ constitutes an
orthonormal frame on $\mS$. One has in fact
\begin{equation}\label{ZYT}
Z \wedge Y\ =\ T\ ,\quad\quad Y \wedge T\ =\ Z\ ,\quad\quad T
\wedge Z\ =\ Y\ .
\end{equation}

Moreover, the (Riemannian) divergence in $\HH$ of these vector
fields is given by
\begin{equation}\label{divY}
div\ Y\ =\ X_1 \pb\ +\ X_2 \qb\ =\ \mathcal H\ ,\quad\quad\quad div\ Z\ =\ X_1 \qb\ -\ X_2 \pb\ .
\end{equation}

Using Cramer's rule one easily obtains from
\eqref{Y} and \eqref{Z}
\begin{equation}\label{Xs}
X_1 \ =\ \pb\ Y\ +\ \qb\ Z\ ,\quad\quad\quad X_2\ =\ \qb\ Y\ -\ \pb\
Z\ .
\end{equation}

One also has
\begin{equation}\label{deltas}
\one \ =\ \qb\ Z\ ,\quad\quad\quad\quad \two \ =\ -\ \pb\ Z\ .
\end{equation}

To prove \eqref{deltas} we proceed as follows
\begin{align*}
& \one \zeta\ =\ X_1\zeta\ -\ <\nabh \zeta,\nuX>\boldsymbol \nu^H_1\
=\ X_1\zeta\ -\ \big(\pb\ X_1\zeta\ + \qb\ X_2\zeta\big)\ \pb
\\
& =\ X_1\zeta\ -\ \pb^2\ X_1\zeta\ - \pb\ \qb\ X_2\zeta\ =\ \qb^2\
X_1\zeta\ -\ \pb\ \qb\ X_2\zeta\ =\ \qb\ \big(\qb\ X_1\zeta\ -\
\pb\ X_2\zeta\big)
\\
& =\ \qb\ Z\zeta\ ,
\end{align*}
\begin{align*}
& \two \zeta\ =\ X_2\zeta\ -\ <\nabh \zeta,\nuX>\boldsymbol \nu^H_2\
=\ X_2\zeta\ -\ \big(\pb\ X_1\zeta\ + \qb\ X_2\zeta\big)\ \qb
\\
& =\ X_2\zeta\ -\ \pb\ \qb\ X_1\zeta\ -\  \qb^2\ X_2\zeta\ =\
\pb^2\ X_2\zeta\ -\ \pb\ \qb\ X_1\zeta\ =\ -\ \pb\ \big(\qb\
X_1\zeta\ -\ \pb\ X_2\zeta\big)
\\
& =\ -\ \pb\ Z\zeta\ .
\end{align*}

These formulas give \begin{equation}\label{delzeta}
 \del \zeta\ =\
\qb\ Z\zeta\ X_1\ -\ \pb\ Z\zeta\ X_2\ . \end{equation}

From \eqref{delzeta} and \eqref{pbar} we obtain
\begin{equation}\label{deltasquare}
|\del \zeta|^2\ =\ (Z\zeta)^2\ =\ (\qb X_1\zeta\ -\ \pb X_2
\zeta)^2\ .
\end{equation}

\medskip

We next establish some identities that will be used times and again
in Sections \ref{S:1&2var} and \ref{S:stab}.

\medskip

\begin{lemma}\label{L:zeroder}
One has on $\mathcal S \setminus \Sigma_\mS$
\begin{equation}\label{0der}
\pb\ Z\pb +\ \qb\ Z\qb\ =\ \pb\ Y\pb +\ \qb\ Y\qb\ =\ \pb\ T\pb +\
\qb\ T\qb\ =\ 0\ ,
\end{equation}
\begin{equation}\label{2derZ}
\pb\ Z^2 \pb\ +\ \qb\ Z^2\qb\ =\ -\ (Z\pb)^2\ -\ (Z\qb)^2\ .
\end{equation}
It is useful to note the following alternative expression of the
first two identities in \eqref{0der}
\begin{equation}\label{zerovargen}
\pb\ \qb\ X_1\pb - \pb^2 X_2\pb + \qb^2 X_1\qb - \pb\ \qb\ X_2\qb\
=\ 0\ , \end{equation}
\begin{equation}\label{zerovar2gen}
\pb^2\ X_1\pb\ +\ \pb\ \qb\ X_2\pb\ +\ \pb\ \qb\ X_1\qb\ +\ \qb^2\
X_2\qb\ =\ 0\ .
\end{equation}
\end{lemma}

\begin{proof}[\textbf{Proof}]
The proof of \eqref{0der} follows trivially by differentiating the
identity $\pb^2 + \qb^2 \equiv 1$, whereas \eqref{2derZ} follows
by differentiating $\pb Z\pb + \qb Z\qb = 0$ with respect to $Z$.
One has from \eqref{0der} and \eqref{Z}
\[
0\ =\ \pb\ Z\pb\ +\ \qb\ Z\qb\ =\ \pb (\qb\ X_1\pb - \pb\ X_2\pb)\
+\ \qb (\qb\ X_1\qb - \pb\ X_2\qb)\ ,
\]
which proves \eqref{zerovargen}. Similarly,
\[
0\ =\ \pb\ Y\pb\ +\ \qb\ Y\qb\ =\ \pb (\pb\ X_1\pb + \qb\ X_2\pb)\
+\ \qb (\pb\ X_1\qb + \qb\ X_2\qb)\ ,
\]
which implies \eqref{zerovar2gen}.

\end{proof}

\medskip

\begin{lemma}\label{L:basicid}
One has on $\mathcal S \setminus \Sigma_\mS$
\begin{equation}\label{Yphi}
<Z,\bN>\ =\ 0\ ,\quad\quad\quad <Y,\bN>\ =\ W\ ,
\end{equation}
\begin{equation}\label{YTphi}
Y\omega\ =\ TW\ ,
\end{equation}
\begin{equation}\label{zerovar4}
 \qb\ Y\pb\ -\ \pb\ Y\qb\ =\ X_2\pb\ -\ X_1\qb\ ,
\end{equation}
\begin{equation}\label{divZ}
\frac{ZW}{W}\ =\  \qb\ Y\pb\ -\ \pb\ Y\qb\ +\ \ob\ ,
\end{equation}
and
\begin{equation}\label{ZTphi/W}
\frac{Z\omega}{W}\ =\ \qb\,T\pb\ -\ \pb\,T\qb\ .
\end{equation}
\end{lemma}

\begin{proof}[\textbf{Proof}]
The first identity in \eqref{Yphi} is obvious, while the second one
is simply a reformulation of \eqref{ip}. The identity \eqref{YTphi}
is just a special case of \eqref{Yoms}. To prove \eqref{zerovar4},
it suffices to use \eqref{Xs} and \eqref{0der} to find
\[
X_2\pb - X_1\qb\ =\ \qb Y\pb - \pb Y\qb - (\pb Y\pb + \qb Y\qb)\
=\  \qb Y\pb - \pb Y\qb\ . \]

As for \eqref{divZ} we  have
\begin{equation*}
\omega\ =\ T\phi\ =\ X_1X_2\phi\ -\ X_2X_1\phi\ =\ X_1(\qb\ W)\ -\
X_2(\pb\ W)\ =\  -\ (X_2\pb\ -\ X_1\qb)\ W\ +\ ZW\
 , \end{equation*}
from which the desired conclusion follows immediately.

Finally, we turn to the proof of \eqref{ZTphi/W}.
Applying $T$ to both sides of \eqref{Yphi} we obtain

\begin{align*}
0 \ =\ T(Z\phi)
&\ =\ T(\qb X_1\phi\ -\ \pb X_2\phi)
\ =\ T\qb X_1\phi\ +\ \qb\, TX_1\phi\ -\ T\pb X_2\phi\ -\ \pb\, TX_2\phi\\
& \ =\ T\qb\,X_1\phi - T\pb\,X_2\phi\ +\ \qb X_1T\phi - \pb X_2T\phi
\ =\ pT\qb - q T\pb + Z(T\phi)\ .
\end{align*}
It follows that
\[
\frac{Z\omega}{W}\ =\ \frac{Z(T\phi)}{W}\ =\ \qb\, T\pb - \pb\,
T\qb\ .
\]
\end{proof}

\medskip

\begin{cor}\label{C:Zob}
One has on $\mS\setminus \Sigma_\mS$
\begin{align*}
\mathcal A\ \overset{def}{=}\ -\  Z\ob\ & =\ (\pb T\qb - \qb T\pb) +
\ob\ (\qb Y\pb - \pb Y\qb) + \ob^2
\\
& =\ \pb (T\qb - \ob Y\qb)\ -\ \qb (T\pb - \ob Y\pb) + \ob^2\ .
\end{align*}
\end{cor}

\begin{proof}[\textbf{Proof}]
We have
\[
Z\ob\ =\ \frac{Z\omega}{W}\ -\ \ob\ \frac{ZW}{W}\ ,
\]
so the desired result follows immediately from \eqref{divZ},
\eqref{ZTphi/W}.

\end{proof}

\medskip

The next lemma expresses a useful orthogonality property which
enters several times in the computations of section \ref{S:1&2var}.

\medskip

\begin{lemma}\label{L:perp}
Let $\mathcal X , \mathcal Y$ be smooth vector fields on $\mathcal
S$, then on the set $\mathcal S\setminus \Sigma_\mS$ one has
\[
\mathcal X \qb\ \mathcal Y \pb\ -\ \mathcal X \pb\ \mathcal Y \qb\
=\ 0\ .
\]
In particular, letting $\mathcal X = Y$ or $T$, and $\mathcal Y = Z$
or $Y$, we find
\[
Y \qb\ Z \pb\ -\ Y \pb\ Z \qb\ =\ 0\ ,
\]
\[
T \qb\ Z \pb\ -\ T \pb\ Z \qb\ =\ 0\ ,
\]
\[
T \qb\ Y \pb\ -\ T \pb\ Y \qb\ =\ 0\ .
\]
\end{lemma}

\begin{proof}[\textbf{Proof}]
To prove the lemma we note that
\[
\mathcal X W\ =\ \pb\ \mathcal X p\ +\ \qb\ \mathcal X q\
,\quad\quad\quad \mathcal Y W\ =\ \pb\ \mathcal Y p\ +\ \qb\
\mathcal Y q\ ,
\]
and proceed as follows
\begin{align*}
& \mathcal X\qb\ \mathcal Y\pb\ -\ \mathcal X\pb\ \mathcal Y \qb\ =\
\mathcal X(qW^{-1}) \mathcal Y(pW^{-1})\ -\ \mathcal X(pW^{-1})
\mathcal Y(qW^{-1})
\\
& =\ \frac{1}{W^2} \big\{(\mathcal X q - \qb \mathcal XW) (\mathcal
Y p - \pb \mathcal Y W)\ -\ (\mathcal X p - \pb \mathcal X
W)(\mathcal Y q - \qb \mathcal Y W)\big\}
\\
&  =\  \frac{ - p^2 \mathcal X q \mathcal Y p - p q \mathcal X q
\mathcal Y q
 - pq \mathcal X p \mathcal Y p - q^2 \mathcal X q \mathcal Y p + p q \mathcal X p \mathcal Y p + q^2 \mathcal X p \mathcal Y q +
  p^2 \mathcal X p \mathcal Y q
+ p q \mathcal X q \mathcal Y q}{W^4}
\\
&  +\ \frac{\mathcal X q \mathcal Y p - \mathcal X p \mathcal Y
q}{W^2}\ =\ \frac{q^2 \mathcal X q \mathcal Y p - p^2 \mathcal X p
\mathcal Y q - q^2 \mathcal X q \mathcal Y p
 + p^2 \mathcal X p \mathcal Y q}{W^4}\ =\ 0\ .
\end{align*}

\end{proof}

In the following lemma we collect some geometric identities
involving the $H$-mean curvature of $\mS$ which play an essential
role in the sequel.

\medskip

\begin{lemma}\label{L:zerovar}
One has on the set $\mathcal S \setminus \Sigma_\mS$
\begin{equation}\label{wedge7}
\qb^2 X_1\pb \  -\ \pb\ \qb \big(X_2 \pb + X_1\qb\big)
\ +\ \pb^2 X_2\qb\ =\ \mathcal H\ ,
\end{equation}
\begin{equation}\label{zerovar3}
\qb\ Z\pb\ -\ \pb\ Z\qb\ =\ \mathcal H\ ,
\end{equation}
\begin{equation}\label{Zp&Zq}
Z\pb\ =\ \qb\ \mathcal H\ ,\quad\quad\quad Z\qb\ =\ -\ \pb\
\mathcal H\ . \end{equation}
 The following formula is dual
to \eqref{wedge7}, \eqref{zerovar3},
\begin{equation}\label{zerovar3ter}
 \pb\ \qb\ X_1\pb\ +\ \qb^2 X_2\pb\ -\ \pb^2
X_1\qb\ -\ \pb\ \qb\ X_2 \qb\ =\ X_2\pb\ -\ X_1 \qb\ ,
\end{equation}
We also have the following expressions for the derivatives of the
$H$-mean curvature along $Y$ and $T$
\begin{equation}\label{YH}
\qb\ Y(Z\pb)\ -\ \pb\
Y(Z\qb)\ =\ Y \mathcal H\ , \end{equation}
\begin{equation}\label{TH}
\qb\ T(Z\pb)\ -\ \pb\ T(Z\qb)\ =\ T \mathcal H\ .
\end{equation}
\end{lemma}

\begin{proof}[\textbf{Proof}]
In view of \eqref{HH} one has that \eqref{wedge7} is equivalent to
\[
\pb^2 X_2\qb + \qb^2 X_1\pb - \pb\ \qb \big(X_2 \pb + X_1
\qb\big)\ =\ X_1\pb\ +\ X_2\qb\ ,
\]
which is in turn equivalent to
\[
\qb^2 X_2\qb + \pb^2 X_1\pb + \pb\ \qb \big(X_2 \pb + X_1
\qb\big)\ =\ 0\ ,
\]
and this is nothing but \eqref{zerovar2gen}. We now use
\eqref{wedge7} to prove \eqref{zerovar3} as follows
\begin{align*}
 & \qb\ Z\pb\ -\ \pb\ Z\qb\ =\ \qb^2\ X_1\pb\ -\ \pb\ \qb\ X_2\pb\ -\
 \pb\ \qb\ X_1\qb\ +\ \pb^2\ X_2\qb\  =\ \mathcal H\ .
\end{align*}

The proof of \eqref{Zp&Zq} immediately follows from the equation
$\pb Z\pb + \qb Z\qb = 0$, from \eqref{zerovar3},
 and from Cramer's rule. Next, it
is easy to recognize that \eqref{zerovar3ter} is equivalent to
\eqref{zerovargen}. The proof of \eqref{YH} follows from
differentiating \eqref{zerovar3} with respect to $Y$, upon using
Leibniz rule and Lemma \ref{L:perp}. Similarly, we establish
\eqref{TH} by differentiating \eqref{zerovar3} with respect to $T$
and then using Lemma \ref{L:perp}.

\end{proof}

\medskip

We now establish a result which says that one of the two
horizontal principal
 curvatures is zero. We stress that this phenomenon, whose Riemannian
counterpart is obviously not true, reflects the fact that
$H$-minimal surfaces are ruled surfaces, see \cite{GP}, \cite{CHMY}.

\medskip

\begin{prop}\label{P:surprise}
One has on $\mathcal S \setminus \Sigma_\mS$
\[
|\one \boldsymbol{\nu}^H_{1}|^2\ +\ |\two \boldsymbol{\nu}^H_{2}|^2\
=\ (Z\pb)^2\ +\ (Z\qb)^2\ \equiv\ \mathcal H^2 \ .
\]
In particular, if $\mS$ is $H$-minimal, we have
\[
|\one \boldsymbol{\nu}^H_{1}|^2\ =\ |\two \boldsymbol{\nu}^H_{2}|^2\
=\ 0\ .
\]
\end{prop}

\begin{proof}[\textbf{Proof}]
According to \eqref{deltasquare}, \eqref{zerovar3} and
\eqref{Zp&Zq}, we have
\begin{align*}
& \mathcal H^2\ -\ |\one \boldsymbol{\nu}^H_{1}|^2\ -\ |\two
\boldsymbol{\nu}^H_{2}|^2\   =\ (\qb Z\pb - \pb Z\qb)^2\ -\
(Z\pb)^2\ -\ (Z\qb)^2
\\
& =\ \qb^2 (Z\pb)^2\ +\ \pb^2 (Z\qb)^2\ - 2\ \pb\ \qb\ Z\pb Z\qb \
-\ (Z\pb)^2\ -\ (Z\qb)^2
\\
& =\  -\ (\pb^2 (Z\pb)^2\ +\ \qb^2 (Z\qb)^2\ +\ 2\ \pb\ \qb\ Z\pb
Z\qb)\ =\ -\ (\pb Z\pb + \qb Z\qb)^2\ =\ 0\ ,
\end{align*}
where the last equation follows from Lemma \ref{L:zeroder}.

\end{proof}

\medskip

\begin{lemma}\label{L:HZY} One has on
$\mathcal S \setminus \Sigma_\mS$
\[
Z\pb\ X_1 \ +\ Z\qb\ X_2 \ =\  \mathcal H\ Z\ ,\quad\quad\quad Z\qb\
X_1\ -\ Z\pb\ X_2\ =\ -\ \mathcal H\  Y\ . \] \end{lemma}

\begin{proof}[\textbf{Proof}]
One easily obtains from the equations \eqref{Xs}
\begin{align*}
Z\pb\ X_1 \ +\ Z\qb\ X_2 \ & =\ (\pb Z\pb + \qb Z\qb)\ Y\ +\
(\qb Z\pb - \pb Z\qb)\ Z \\
&  =\ (\qb Z\pb - \pb Z\qb)\ Z\ =\ \mathcal H\ Z\ ,
\end{align*}
where in the second to the last equality we have used
\eqref{0der}, and in the last one we have used \eqref{zerovar3}.
The proof of the second identity is similar.

\end{proof}

\medskip

The next commutator formulas will be useful in the sequel.

\medskip

\begin{lemma}\label{L:comm}
One has on $\mS\setminus \Sigma_\mS$
\begin{equation*}
[Z,Y] \ =\ T\ +\ \mathcal H\ Z \  +\ (\qb Y \pb - \pb Y \qb)\ Y\ .
\end{equation*}
\end{lemma}

\begin{proof}[\textbf{Proof}]
To compute the commutator between $Z$ and $Y$ we use the equations
\eqref{Y} and \eqref{Z} to find
\begin{align*}
& [Z,Y] \ =\ Z(Y)\ -\ Y(Z)\\
& =\ \qb X_1(\pb X_1 + \qb X_2 ) - \pb X_2(\pb X_1  + \qb X_2 ) -
\pb X_1(\qb X_1  - \pb X_2 ) - \qb X_2(\qb X_1  - \pb X_2 )
\\
& =\ X_1 X_2  - X_2 X_1  + \big(\qb X_1 \pb - \pb X_1 \qb\big)\
X_1 \ +\ \big(\qb X_2 \pb - \pb X_2 \qb\big)\ X_2  \notag\\
& =\ T\ +\ \big(Z\pb + \pb (X_2 \pb - X_1 \qb)\big) X_1 \ +\
\big(Z\qb + \qb (X_2 \pb - X_1 \qb)\big) X_2  \notag\\
& =\ T\ +\ Z\pb\ X_1 \ +\ Z\qb\ X_2 \ +\ (X_2 \pb - X_1 \qb)\ Y\
,\notag
\end{align*}
where we have repeatedly used \eqref{0der} along with the identity $\pb^2 +
\qb^2 = 1$. We now appeal to Lemma \ref{L:HZY} and to
\eqref{zerovar3ter} to reach the desired conclusion.

\end{proof}

\medskip

\begin{lemma}\label{L:commT}
On $\mathcal S \setminus \Sigma_\mS$, one has
\begin{equation*}
[Z,T] \ =\ (\qb T \pb - \pb T \qb)\ Y\ .
\end{equation*}
\end{lemma}

\begin{proof}[\textbf{Proof}]
Using the trivial commutation relations $[X_i,T] = 0$, $i=1,2$ we
obtain
\[
T(Z)\ =\ T\qb\ X_1\ -\ T\pb\ X_2 \ +\ Z(T)\ .
\]

From this identity, and from \eqref{Xs}, we obtain
\begin{align*}
& [Z,T]\ =\  T\pb\ X_2\ -\ T\qb\ X_1 \ =\ T\pb (\qb Y - \pb Z)\ -\
T\qb (\qb Z + \pb Y)
\\
& =\ (\qb T\pb - \pb T\qb) Y\ -\ (\pb T\pb + \qb\ T\qb)\ Z\ =\ (\qb
T\pb - \pb T\qb) Y\ ,
\end{align*}
where in the last equality we have used Lemma \ref{L:zeroder}.

\end{proof}

\medskip

\begin{cor}\label{C:mixedcomm}
One has on $\mS\setminus \Sigma_\mS$
\[
[T - \ob Y,Z]\ =\ \ob\ \bigg\{(T - \ob\ Y) + \mathcal H\ Z\bigg\}\ .
\]
In particular, if $\mS$ is $H$-minimal, then
\[
[T - \ob Y,Z]\ =\ \ob\ (T - \ob\ Y)\ .
\]
\end{cor}

\begin{proof}[\textbf{Proof}]
One has from Lemmas \ref{L:comm} and \ref{L:commT}
\begin{align*}
[T - \ob Y,Z]\ & =\ [T,Z] - \ob\ [Y,Z] + Z\ob\ Y
\\
& -\ (\qb T\pb - \pb T\qb) Y + \ob\ (T + (\qb Y\pb - \pb Y \qb)) Y -
\mathcal A\ Y
\\
& =\ \ob\ T\ +\ (\ob\ (\qb Y\pb - \pb Y \qb) - (\qb T\pb - \pb
T\qb)\ -\ \mathcal A)\ Y\ ,
\end{align*}
where we have used the hypothesis that $\mS$ be $H$-minimal. Using
Corollary \ref{C:Zob} we reach the desired conclusion.

\end{proof}

\medskip

\begin{cor}\label{C:za}
If $\mS$ is $H$-minimal, one has
\[
Z\mathcal A\ =\ \ob (\ob^2 - 3 \mathcal A)\ .
\]
For an arbitrary $C^2$ surface we have instead on $\mS\setminus
\Sigma_\mS$
\[
Z\mathcal A\ =\ \ob \bigg\{(\ob^2 - 3 \mathcal A) + \mathcal
H^2\bigg\}\ .
\]
\end{cor}

\begin{proof}[\textbf{Proof}]
From the assumption of $H$-minimality of $\mS$, \eqref{Zp&Zq} and
from Corollary \ref{C:Zob} we obtain
\begin{align*}
Z \mathcal A\ & =\ \pb Z(T\qb - \ob Y\qb) - \qb Z(T\pb - \ob Y\pb) +
2 \ob Z\ob
\\
& =\ \pb [Z,T - \ob Y]\qb - \qb [Z,T - \ob Y]\pb - 2 \ob \mathcal A
\\
& =\ -\ \ob \bigg[\pb (T\qb - \ob Y\qb)  - \qb (T\pb - \ob Y\pb) + 2
\mathcal A\bigg]\ =\ \ob (\ob^2 - 3 \mathcal A)\ .
\end{align*}

The proof of the second part of the corollary is based on a longer
computation which, in addition, exploits also Corollary
\ref{C:mixedcomm}, \eqref{0der}, \eqref{zerovar3} and \eqref{Zp&Zq}.
We leave the details to the interested reader.

\end{proof}

\medskip

\begin{cor}\label{C:comm} One has on $\mathcal S \setminus \Sigma_\mS$
\begin{equation*}
[Z,Y] \pb\ =\ T\pb\ +\ Z \pb\ \mathcal H\ +\ (\qb Y \pb - \pb Y
\qb)\ Y\pb\  ,
\end{equation*}
\[
[Z,Y]\qb\ =\ T\qb\ +\  Z\qb\ \mathcal H\ +\ (\qb Y\pb - \pb Y
\qb)\ Y\qb\ ,
\]
\[
[Z,T] \pb\ =\ (\qb T \pb - \pb T \qb)\ Y\pb\ ,\quad\quad\quad
[Z,T] \qb\ =\ (\qb T \pb - \pb T \qb)\ Y\qb\ .
\]
\end{cor}

\medskip

\begin{lemma}\label{L:Zcomm}
On $\mathcal S \setminus \Sigma_\mS$ one has
\[
Z(\qb Y\pb\ -\ \pb Y\qb)\ -\ (\qb Y\pb\ -\ \pb Y\qb)^2\ -\ (\qb
T\pb - \pb T\qb)\ =\ Y\mathcal H\ +\ \mathcal H^2\  .
\]
\end{lemma}

\begin{proof}[\textbf{Proof}]
Thanks to Lemma \ref{L:perp} we have
\begin{align*}
& Z(\qb Y\pb\ -\ \pb Y\qb)\ =\ \qb\ Z( Y\pb)\ -\ \pb\ Z(Y\qb)
\\
& =\ \qb\ Y( Z\pb)\ -\ \pb\ Y(Z\qb)\ +\ \qb\ [Z,Y]\pb\ -\ \pb\
[Z,Y]\qb \notag\\
& =\ Y\mathcal H\ +\ \qb \big\{T\pb + Z\pb\ \mathcal H\ + (\qb
Y\pb - \pb Y\qb) Y\pb\big\} \ -\ \pb \big\{T\qb + Z\qb\ \mathcal
H\ + (\qb Y\pb - \pb Y\qb) Y\qb\big\} \notag\\
& =\ Y\mathcal H\ +\ (\qb T\pb - \pb T\qb)\ +\ (\qb Z\pb - \pb
Z\qb)\ \mathcal H\ +\ \big(\qb Y\pb - \pb
Y\qb\big)^2 \notag\\
& =\ Y\mathcal H\ +\ \mathcal H^2\ +\ (\qb T\pb - \pb T\qb)\ +\
\big(\qb Y\pb - \pb Y\qb\big)^2\ , \notag
\end{align*}
where we have used \eqref{YH} and Corollary \ref{C:comm} in the
third equality, and \eqref{zerovar3} in the second to the last
equality.

\end{proof}

\medskip

\begin{lemma}\label{L:Zcomm2}
One has on the set $\mathcal S \setminus \Sigma_\mS$
\[
Z(\qb\ T\pb\ -\ \pb\ T\qb)\ =\ T\mathcal H\ +\ (\qb\ T\pb\ -\ \pb\
T\qb)\ (\qb\ Y\pb\ -\ \pb\ Y\qb)\ .
\]
\end{lemma}

\begin{proof}[\textbf{Proof}]
Using Lemma \ref{L:perp} we obtain
\begin{align*}
& Z(\qb\ T\pb\ -\ \pb\ T\qb)\ =\ \qb\ Z(T\pb)\ -\ \pb\ Z(T\qb)
\\
& =\ \qb\ T(Z\pb)\ -\ \pb\ T(Z\qb)\ +\ \qb\ [Z,T]\pb\ -\ \pb\
[Z,T]\qb\ .
\end{align*}

Now using Lemma \ref{perp} again we obtain from \eqref{zerovar3}
\[
T\mathcal H\ =\ \qb\ T(Z\pb)\ -\ \pb\ T(Z\qb)\ ,
\]
which, substituted in the above equation gives, along with
Corollary \ref{C:comm}, the desired result.

\end{proof}

\medskip

\vskip 0.6in


\section{\textbf{First and second variation of the $H$-perimeter
in the Heisenberg group}}\label{S:1&2var}

\vskip 0.2in

A fundamental tool in Riemannian geometry are the first and second
variation formulas for the area functional. Consider a $C^2$
oriented hypersurface $\mathcal S \subset \mathbb R^{n}$, with Gauss
map $\nu :\mathcal S \to \mathbb S^{n-1}$, and denote by $\mathcal
S^\lambda = G_\lambda(\mathcal S)$ the hypersurface obtained by
deforming $\mathcal S$ in the normal direction with the
one-parameter family of local diffeomorphisms $G_\lambda(x) = x +
\lambda \zeta(x) \nu(x)$, where $\zeta\in C^\infty_0(\mS)$, and
$\lambda \in \R$ is small. One has the following theorem, see for
instance (10.12), (10.13) in \cite{Gi}, or also \cite{MM},
\cite{Si}, \cite{CM} and \cite{BGG}.

\begin{thrm}\label{T:MM}
The first variation of the area of $\mathcal S$ is given by the
formula \begin{equation}\label{fvclass} \frac{d}{d\lambda}
H_{n-1}(G_\lambda(\mathcal S))\Bigl|_{\lambda = 0}\ =\
\int_\mathcal S  H\ \zeta\ dH_{n-1}\ ,
\end{equation}
where $H = \kappa_1 + ... +\kappa_{n-1}$ indicates the sum of the
principal curvatures of $\mathcal S$. The second variation is
given by
\begin{equation}\label{svclass}
\frac{d^2}{d\lambda^2} H_{n-1}(G_\lambda(\mathcal S))\Bigl|_{\lambda
= 0}\ =\ \int_\mathcal S \bigg\{ |\nabla \zeta|^2\ +\  \zeta^2\
\bigg(H^2\ -\ \sum_{i=1}^{n} |\nabla \nu_i|^2\bigg)\ \bigg\}\
dH_{n-1}\ ,
\end{equation}
where $\nabla$ denotes the Levi-Civita connection on $\mathcal S$,
and it can be shown that $\sum_{i=1}^{n} |\nabla \nu_i|^2$ is the
sum of the squares of the principal curvatures of $\mS$.
\end{thrm}

\medskip

In this section we consider an oriented surface in the Heisenberg
group $\HH$, with non-unit Riemannian normal $\bN$, and horizontal
Gauss map $\nuX$, and compute the first and second variation for
general deformations of $\mathcal S$. We observe that, in view of
applications to the fundamental question of stability of $H$-minimal
surfaces, it is important to be able to treat general deformations,
versus deformations along a specific direction. More precisely, we
consider deformations of $\mathcal S$ given by $\mS \to \mS^\la =
J_\la(S)$, where
\begin{equation}\label{def}
J_\lambda(\mS)\ =\ \mS + \la \mathcal X\ =\ \mS + \lambda \bigg(a
X_1 + b X_2 + k T\bigg)\ ,
\end{equation}
where $a, b , k\in C^1_0(\mathcal S\setminus \Sigma_\mS)$, and $\la
\in \R$ is a small parameter. Throughout this section, and the
following one, we will continue to use the notations of section
\ref{S:geomid}. We recall that $\HH$ is endowed with a
left-invariant Riemannian metric with respect to which $\{X_1, X_2,
T\}$ constitute an orthonormal basis with inner product
$<\cdot,\cdot>$. Since no other inner product will be used, there
will not be any confusion, for instance, with the standard Euclidean
inner product of $\mathbb R^3$.

\medskip

\begin{dfn}\label{D:variations}
Let $\mathcal S \subset \HH$ be an oriented $C^2$ surface, consider
the family of vector fields $\mathcal X = a X_1 + b X_2 + k T$, with
$a, b , k\in C^2_0(\mathcal S\setminus \Sigma_\mS)$, and the family
of surfaces $\mS^\la$. We define the \emph{first variation} of the
$H$-perimeter with respect to the deformation \eqref{def} as
\[
\fv\ =\ \frac{d}{d\lambda}~ P_H(\mathcal S^\lambda)\Bigl|_{\lambda =
0}\ .
\]
If $\Sigma_\mS = \varnothing$, then we say that $\mS$ is
\emph{stationary} if $\fv = 0$, for every $\mathcal X$.
\end{dfn}

\medskip

Classical minimal surfaces are stationary points of the perimeter
(the area functional for graphs). It is natural to ask what is the
connection between the notion of $H$-minimal surface and that of
$H$-perimeter. The answer to this question is contained in the
following result. To simplify the formulas we introduce the
following notation
\begin{equation}\label{inner} F\ \overset{def}{=}\ \pb a + \qb b +
\ob k\ =\ \frac{<\mathcal X,\bN>}{<\nuX,\bN>}\ .
\end{equation}

\medskip

 \begin{thrm}\label{T:variations}
Let $\mS\subset \HH$ be an oriented $C^2$ surface, then
\begin{equation}\label{fvH}
\mathcal V^H_I(\mS;\mathcal X)\ =\
 \int_{\mathcal S}
\mathcal H\ F\ d\sigma_H\ .
\end{equation}
In particular, $\mathcal S$ is stationary if and only if it is
$H$-minimal.
\end{thrm}

\medskip

Versions of Theorem \ref{T:variations} have also been obtained
independently by other people. An approach based on motion by
$H$-mean curvature can be found in \cite{BC}. When $\mathcal X = a
\nuX + k T$, then a proof based on CR-geometry can be found in
\cite{CHMY}, and \cite{RR1}, \cite{RR2}. We mention that Hladky and
Pauls have recently proved in \cite{HP} (with a different approach
which does not directly use the first variation) that, for a wide
class of sub-Riemannian spaces, a non-characteristic $C^2$
hypersurface is a critical point of the $H$-perimeter if and only if
it is $H$-minimal.

\medskip

\begin{dfn}\label{D:minimal}
Given an oriented $C^2$ surface $\mS \subset \HH$, we define the
\emph{second variation} of the $H$-perimeter with respect to the
deformation \eqref{def} as
\[
\sv\ =\ \frac{d^2}{d\lambda^2}~ P_H(\mathcal
S^\lambda)\Bigl|_{\lambda = 0}\ .
\]
\end{dfn}

\medskip

Our main result in this section is the following theorem.

\medskip

\begin{thrm}\label{T:2varfinal}
The second variation of the $H$-perimeter with respect to the
deformation of $\mS$ given by \eqref{def} is expressed by the
formula
\begin{align}\label{uffa9}
\sv\ & =\
  \int_\mS \bigg\{ 2\ (\qb Za - \pb Zb)  \left(Tk -
\ob Yk\right) \\
& +\ \left(Ta - \ob Ya\right) \bigg[- 2 \qb Zk - \qb
(a \pb + b \qb) - \pb (a \qb - b \pb)\bigg] \notag\\
& +\ \left(Tb - \ob Yb\right)\bigg[2 \pb Zk + \pb (a
\pb + b \qb) - \qb (a \qb - b \pb)\bigg] \notag\\
& +\ 2\ (a \qb - b \pb) (\qb Za - \pb Zb)\ob
\notag\\
& +\ \left(Za + \pb\ \ob\ Zk\right)^2\ +\ \left(Zb + \qb\ \ob\
Zk\right)^2
\notag\\
& +\  (a^2 + b^2)\ \ob^2
\notag\\
& +\ 2\ \ob (a Za + b Zb)\ +\ 2\
\ob^2 (a \pb + b \qb) Zk \notag\\
& -\ \big(\qb Za - \pb Zb + (a \qb  - b \pb) \ob\big)^2\bigg\}\
d\sigma_H\ . \notag
\end{align}
\end{thrm}

\medskip

\medskip

In order to prove Theorems \ref{T:variations} and \ref{T:2varfinal}
we develop some preliminary material which constitutes the necessary
geometric backbone. We begin by deriving from Theorem \ref{T:ibp}
two integration by parts formulas which play a fundamental role in
this section and in the following one.

\medskip

\begin{lemma}\label{L:Zf}
Let $\zeta\in C_0^1(\mS \setminus \Sigma_\mS)$, then
\[
\int_{\mathcal S} Z\zeta\ d\sigma_H\ =\ -\ \int_{\mathcal S} \zeta\
\ob\  d\sigma_H\ .
\]
\end{lemma}

\begin{proof}[\textbf{Proof}]
We begin by noting that, thanks to \eqref{deltas}, \eqref{Y} and
\eqref{Z}, we can rewrite the two identities in \eqref{ibpH1} in the
form
\begin{equation}\label{intbyparts1}
\int_\mathcal S \left\{\qb\ Zu\ +\ \qb\ u\ \ob\right\}\ d\sigma_H\
=\ \int_\mathcal S \pb\ u\ \mathcal H\ d\sigma_H\ ,
\end{equation}
\begin{equation}\label{intbyparts2}
\int_\mathcal S \left\{\pb\ Zu\ +\ \pb\ u\ \ob \right\}\ d\sigma_H\
=\ -\ \int_\mathcal S \qb\ u\ \mathcal H\ d\sigma_H\ .
\end{equation}

Choosing $u = \qb \zeta$ in \eqref{intbyparts1}, and $u = \pb \zeta$
in \eqref{intbyparts2}, and adding the resulting equations, we
obtain
\[
\int_\mathcal S \left\{(\pb Z\pb + \qb Z\qb)\ \zeta\ +\ (\pb^2 +
\qb^2)\ Z\zeta\right\}\ d\sigma_H\ +\ \int_\mathcal S (\pb^2 +
\qb^2)\ \zeta\ \ob\ d\sigma_H\ =\ 0\ .
\]

From the latter equation, and from \eqref{pbar}, \eqref{0der}, we
immediately reach the conclusion.

\end{proof}

\medskip

\begin{rmrk}\label{R:equivalence}
We have above derived Lemma \ref{L:Zf} from Theorem \ref{T:ibp}
specialized to $\HH$. The two results are in fact equivalent. To see
this suppose that the identity in Lemma \ref{L:Zf} hold. Applying it
twice, once with the choice $\zeta = \qb u$, and the other with
$\zeta= \pb u$, with $u\in C_0^1(\mS \setminus \Sigma)$, we obtain
\eqref{intbyparts1} and \eqref{intbyparts2} if we use the identities
$Z\pb = \qb \mathcal H$, $Z\qb = -  \pb \mathcal H$, in
\eqref{Zp&Zq}.
\end{rmrk}

\medskip

Another crucial integration by parts formula which we will need is
\eqref{twoT&Y} in Theorem \ref{T:T&Y}. For the reader's convenience
we combine this formula and Lemma \ref{L:Zf} into a single
statement.

\medskip

\begin{lemma}\label{L:ibpcombined}
Let $f\in C^1(\mS)$, $\zeta\in C^1_0(\mS\setminus \Sigma_\mS)$, then
\[
\int_{\mathcal S} f\ Z\zeta\ d\sigma_H\ =\ -\ \int_\mS \zeta\ Zf\
d\sigma_H\ -\ \int_{\mathcal S} f\ \zeta\ \ob\ d\sigma_H\ .
\]
\[
 \int_\mathcal S f \left(T\zeta - \ob Y\zeta\right)
  d\sigma_H\ =\  -\ \int_\mS \zeta \left(Tf - \ob Yf\right) d\sigma_H\ +\ \int_\mathcal S f \zeta
 \ob \mathcal H\ d\sigma_H\ .
 \]
In particular, if $\mS$ is $H$-minimal, we find
\[
\int_\mathcal S f \left(T\zeta - \ob Y\zeta\right)
  d\sigma_H\ =\  -\ \int_\mS \zeta \left(Tf - \ob Yf\right) d\sigma_H\
  .
  \]
 \end{lemma}

 \medskip

After these preliminaries we turn to the proofs of the main results
in this section. We first recall the representation formula for the
$H$-perimeter of $\mathcal S$ given in \eqref{permeasure} of
Definition \ref{D:permeas}
\begin{equation}\label{per1}
\sigma_H(\mathcal S)\ =\ \int_{\mathcal S} \frac{W}{|\bN|}\ d\sigma\
=\ \int_{\mathcal S} \frac{\sqrt{p^2 + q^2}}{|\bN|}\ d\sigma\ ,
\end{equation}
where $d\sigma$ represents the standard surface measure on $\mathcal
S$. Next, we establish a simple lemma which provides the general
expression for the first and second variation of the $H$-perimeter.
We consider, for small values of $\la \in \R$, a deformation of
$\mathcal S$ of the type $\mS^\la = J_\lambda(\mS) = \mS + \lambda
\mathcal X$, where $\mathcal X\in C^1_0(\mS\setminus
\Sigma_\mS;\HH)$. We denote by $\bN^\lambda$ the non-unit Riemannian
normal on $\mS^\la$. Letting
\[
X_1^\lambda(g)\ =\ X_1(J_\la(g))\ ,\quad\quad\quad X_2^\lambda(g)\
=\ X_2(J_\la(g))\ ,\quad\quad\quad g\in \mS\ ,
\]
we consider the functions
\begin{equation}\label{plambda0}
p^\lambda\ =\ <\bN^\lambda,X_1^\lambda>\ ,\quad\quad  q^\lambda\
=\ <\bN^\lambda, X_2^\lambda>\ ,\quad\quad W^\lambda\ =\
\sqrt{(p^\lambda)^2 + (q^\lambda)^2}\ .
\end{equation}

We stress that we are assuming that $\mathcal X$ is compactly supported away from the characteristic set
of $\mS$, so that the angle function $W$ for $\mS$ never vanishes on the support of $\mathcal X$,
see \eqref{csaf}.

\medskip

\begin{lemma}\label{L:1&2var}
The first variation of the $H$-perimeter along the deformation
$J_\lambda(\mS) = \mS + \lambda \mathcal X$ is given by the
formula
\begin{equation}\label{1vargen}
\fv\ =\ \int_\mathcal S \frac{\frac{d
W^\lambda}{d\lambda}}{W^\lambda}\Bl\ d\sigma_H\ =\ \int_\mathcal S
\frac{\left(p^\lambda\ \frac{d p^\lambda}{d\lambda} +  q^\lambda\
\frac{d q^\lambda}{d\lambda}\right)\Bl}{ W^2}\ d\sigma_H\ .
\end{equation}
The second variation is given by
\begin{align}\label{2vargen}
\sv\ & =\ \int_\mathcal S \frac{\left( p^\lambda\ \frac{d^2
p^\lambda}{d\lambda^2} + q^\lambda\ \frac{d^2
q^\lambda}{d\lambda^2}\right)\Bl}{W^2}\ d\sigma_H \\
& +\ \int_\mS \frac{\left(\left(\frac{d
p^\lambda}{d\lambda}\right)^2 + \left(\frac{d
q^\lambda}{d\lambda}\right)^2\right)\Bl }{W^2}\ d\sigma_H
\notag\\
& -\ \int_\mathcal S \frac{\left(p^\lambda\ \frac{d
p^\lambda}{d\lambda} +  q^\lambda\ \frac{d
q^\lambda}{d\lambda}\right)^2\Bl}{W^4} \ d\sigma_H\ . \notag
\end{align}

\end{lemma}

\begin{proof}[\textbf{Proof}]
Because of the assumptions on $\mathcal X$, the integrals in the
right-hand sides of \eqref{1vargen}, \eqref{2vargen} are performed
on a compact set which does not intersect $\Sigma_\mS$. By a
partition of unity we can thus reduce the analysis to a set $\mS\cap
\mathcal O$, where $\mathcal O\subset \HH$ is an open neighborhood
of a point $g_0\in \mS\setminus \Sigma_\mS$. We can thus assume that
there exist an open set $\Om \subset \mathbb R^2_{u,v}$ and a
parametrization $\theta:\Om \to \HH$ such that $\mS\cap \mathcal O
=\theta(\Om)$. We suppose that the orientation of $\mS$ is given by
$\bN = \theta_u \wedge \theta_v$. Recalling that $d\sigma =
|\theta_u \wedge \theta_v| du\wedge dv$, we can, after projecting
$\mathcal S$ onto $\Om$, rewrite \eqref{per1} as follows
\begin{equation}\label{pergen}
 \sigma_H(\mathcal
S\cap \mathcal O)\ =\ \int_\Om W\ d\sigma\ =\ \int_\Om \sqrt{p^2 +
q^2} \ du \wedge dv\ .
\end{equation}

According to \eqref{pergen} we have
\begin{equation*}
 \sigma_H(\mathcal
S^\lambda\cap \mathcal O)\ =\ \int_\Om W^\lambda\ du\wedge dv\ .
\end{equation*}

Observing that
\[
\frac{d W^\lambda}{d\lambda}\ =\ \frac{ p^\lambda\ \frac{d
p^\lambda}{d\lambda} +  q^\lambda\ \frac{d q^\lambda}{d\lambda}}{
W^\lambda}\ ,
\]
and that $W^\lambda\Bl = W$, we find
\[
\fv\ =\ \int_\Om \frac{\left( p^\lambda\ \frac{d
p^\lambda}{d\lambda} +  q^\lambda\ \frac{d
q^\lambda}{d\lambda}\right)\Bl}{ W}\ du\wedge dv\ =\ \int_\mS
\frac{\left( p^\lambda\ \frac{d p^\lambda}{d\lambda} +  q^\lambda\
\frac{d q^\lambda}{d\lambda}\right)\Bl}{ W^2}\ d\sigma_H\ ,
\]
which gives \eqref{1vargen}. To obtain \eqref{2vargen} we proceed analogously, observing that
\begin{align*}
\frac{d^2 W^\lambda}{d\lambda^2}\ & =\ \frac{ p^\lambda\ \frac{d^2
p^\lambda}{d\lambda^2} +  q^\lambda\ \frac{d^2 q^\lambda}{d\lambda^2}}{
W^\lambda}\ +\ \frac{\left(\frac{d
p^\lambda}{d\lambda}\right)^2 +  \left(\frac{d q^\lambda}{d\lambda}\right)^2}{
W^\lambda}
\\
& -\ \frac{\left( p^\lambda\ \frac{d
p^\lambda}{d\lambda} +  q^\lambda\ \frac{d q^\lambda}{d\lambda}\right)^2}{
(W^\lambda)^3}\ .
\end{align*}

\end{proof}

\medskip

We now turn to the proof of the main results. Given an open set $\Om
\subset \mathbb R^2$, we denote by $\theta:\Om \to \HH$ a $C^2$
parametrization of an oriented surface
\begin{equation}\label{Spar}
\mathcal S\ =\ \{\theta(u,v) = (x(u,v),y(u,v),t(u,v)) \in \HH \mid
(u,v) \in \Om \}\ . \end{equation}

We assume throughout that the orientation of $\mS$ is given by the
non-unit normal $\bN = \theta_u \wedge \theta_v$. Using
\eqref{coord} we see that \begin{equation}\label{thetaH}
\theta(u,v)\ =\ x(u,v) X_1(\theta(u,v))\ +\ y(u,v)
X_2(\theta(u,v))\ +\ t(u,v) T\ .
\end{equation}

From this equation we find
\[
\theta_u\ =\ x_u X_1 + y_u X_2 + t_u T + x X_{1,u} + y X_{2,u}\ ,
\]
with a similar expression for $\theta_v$. Keeping in mind that
\begin{equation}\label{xuxv}
\begin{cases} X_{1,u}\ =\ - \frac{y_u}{2}\
T\ ,\quad\quad\quad X_{2,u}\ =\  \frac{x_u}{2}\ T\ , \\
X_{1,v}\ =\ - \frac{y_v}{2}\ T\ ,\quad\quad\quad X_{2,v}\ =\
\frac{x_v}{2}\ T\ , \end{cases}
 \end{equation}
 we obtain
\begin{equation}\label{dertheta}
\begin{cases}
\theta_u\ =\  x_u X_1 + y_u X_2 + \left(t_u + \frac{y x_u - x
y_u}{2}\right) T\ , \\
\theta_v\ =\  x_v X_1 + y_v X_2 + \left(t_v + \frac{y x_v - x
y_v}{2}\right) T\ . \end{cases}
\end{equation}

From \eqref{wedges} and \eqref{xuxv} we find
\begin{equation}\label{wedgesX}
\begin{cases}
\theta_u \wedge X_1\ =\ \left(t_u + \frac{y x_u - x y_u}{2}\right)
X_2 - y_u T\ ,
\\
\theta_u \wedge X_2\ =\ - \left(t_u + \frac{y x_u - x
y_u}{2}\right) X_1 + x_u T\ ,
\\
\theta_u \wedge T\ =\ y_u X_1 - x_u X_2\ ,
\\
\theta_v \wedge X_1\ =\ \left(t_v + \frac{y x_v - x y_v}{2}\right)
X_2 - y_v T\ ,
\\
\theta_v \wedge X_2\ =\ - \left(t_v + \frac{y x_v - x
y_v}{2}\right) X_1 + x_v T\ ,
\\
\theta_v \wedge T\ =\ y_v X_1 - x_v X_2\ .
\end{cases}
\end{equation}

The non-unit outer Riemannian normal to $\mathcal S$ is thus given
by
\begin{align}\label{nun}
\bN\ =\  \theta_u\ \wedge\ \theta_{v}\ & =\ \left\{y_u \left(t_v +
\frac{y x_v - x y_v}{2}\right) - y_v \left(t_u + \frac{y x_u - x
y_u}{2}\right)\right\}\ X_1
\\
& +\ \left\{x_v \left(t_u + \frac{y x_u - x y_u}{2}\right) - x_u
\left(t_v + \frac{y x_v - x y_v}{2}\right)\right\}\ X_2 \notag\\
& +\ \big\{x_u y_v - x_v y_u\big\}\ T \notag\\
& =\ \left(y_u t_v - y_v t_u\ -\ \frac{y}{2} (x_u y_v - x_v
y_u)\right)\ X_1 \notag\\
& +\ \left(x_v t_u - x_u t_v\ +\ \frac{x}{2} (x_u y_v - x_v
y_u)\right)\ X_2 \notag\\
& +\ (x_u y_v - x_v y_u)\ T . \notag
\end{align}

 We denote by $\boldsymbol \nu = \bN/|\bN|$ the Riemannian Gauss map of
$\mathcal S$. Keeping in mind \eqref{pq}, we see from \eqref{nun}
that
\begin{equation}\label{ppar}
\begin{cases}
p\  =\ y_u t_v - y_v t_u\ -\ \frac{y}{2} (x_u y_v - x_v y_u)\ ,
\\
q\  =\ x_v t_u - x_u t_v\ +\ \frac{x}{2} (x_u y_v - x_v y_u)\ ,
\\
\om\ =\ x_u y_v - x_v y_u\ . \end{cases}
\end{equation}

We note at this moment that, given the assumption $\theta\in
C^2(\Om)$, the functions $p, q , \omega, W$ are of class $C^1(\Om)$,
and that moreover $\pb, \qb$, $\ob$ are of class $C^1(\mS\setminus
\Sigma_\mS)$. In what follows, given a function $\zeta$ defined in a
neighborhood of $\mS$, we will by abuse of notation denote with
$\zeta(u,v) = \zeta\circ\theta(u,v) = \zeta(x(u,v),y(u,v),t(u,v))$.
The chain rule gives
\begin{equation}\label{crh}
\zeta_u\ =\ x_u\ \zeta_x\ +\ y_u\ \zeta_y\ +\ t_u\ \zeta_t\
,\quad\quad \zeta_v\ =\ x_v\ \zeta_x\ +\ y_v\ \zeta_y\ +\ t_v\
\zeta_t\ .
\end{equation}

Using \eqref{Hn2}, we obtain from \eqref{crh}
\begin{equation}\label{zuzv}
\begin{cases}
\zeta_u\ =\ x_u\ X_1\zeta\ +\ y_u\ X_2\zeta\ +\ \left(t_u +
\frac{y x_u - x y_u}{2}\right) \ T\zeta\ ,
\\
\zeta_v\ =\ x_v\ X_1\zeta\ +\ y_v\ X_2\zeta\ +\ \left(t_v +
\frac{y x_v - x y_v}{2}\right) \ T\zeta\ .
\end{cases}
\end{equation}

\medskip

In the sequel, it will be convenient to also have the expression
of $\zeta_u$, $\zeta_v$ with respect to the orthonormal frame
$\{Z,Y,T\}$, where $Y$ and $Z$ are like in \eqref{Y}, \eqref{Z}.
From \eqref{zuzv} and \eqref{Xs}, we have
\begin{equation}\label{zuzvYZ}
\begin{cases}
\zeta_u\ =\ (x_u \pb + y_u \qb)\ Y\zeta\ +\ (x_u \qb - y_u \pb)\
Z\zeta\ +\ \left(t_u + \frac{y x_u - x y_u}{2}\right) \ T\zeta\ ,
\\
\zeta_v\ =\ (x_v \pb + y_v \qb)\ Y\zeta\ +\ (x_v \qb - y_v \pb)\
Z\zeta\ +\ \left(t_v + \frac{y x_v - x y_v}{2}\right) \ T\zeta\ .
\end{cases}
\end{equation}

We now fix functions $a, b , k \in C^\infty_0(\mS\setminus
\Sigma_\mS)$, and consider the vector field
\begin{equation}\label{vf}
\mathcal X\ =\ a\ X_1\ +\ b\ X_2\ + k\ T\ . \end{equation}

For small values of $\lambda \in \R$, we let $\mathcal S^\lambda$
be the surface obtained by deforming $\mathcal S$ through the map
$J_\lambda = Id + \lambda \mathcal X$, so that
\begin{equation}\label{Jlambda}
J_\la(g)\ =\ g\ +\ \lambda\ (a\ X_1\ +\ b\ X_2\ +\ k\ T)\
,\quad\quad\quad\quad g\in \mS\ .
\end{equation}

The parametric representation of $\mathcal S^\lambda$ is given by
\begin{equation}\label{tl0}
\theta^\la\ =\ \theta\ +\ \la\ \mathcal X\ ,
\end{equation}
so that
\begin{equation}\label{dtl}
\theta_u^\la\ =\ \theta_u\ +\ \la\ \mathcal X_u\ ,\quad\quad\quad
\theta_v^\la\ =\ \theta_v\ +\ \la\ \mathcal X_v\ ,
\end{equation}
and therefore the non-unit Riemannian normal to the surface
$\mS^\la = J_\la(\mS)$ is given by
\begin{equation}\label{Nl}
\bN^\la\ =\ \theta_u^\la \wedge \theta_v^\la\ =\ \bN\ +\ \la\
(\theta_u \wedge\mathcal X_v - \theta_v \wedge \mathcal X_u)\ +\
\la^2\ \mathcal X_u \wedge \mathcal X_v\ .
\end{equation}

From \eqref{Nl} we obtain
\begin{equation}\label{dNl}
\frac{d\bN^\la}{d\la}\Bl\ =\ \theta_u \wedge\mathcal X_v -
\theta_v \wedge \mathcal X_u\ ,\quad\quad\quad
\frac{d^2\bN^\la}{d\la^2}\Bl\ =\ 2\ \mathcal X_u \wedge \mathcal
X_v\ .
\end{equation}

Using \eqref{xuxv} we find
\begin{equation}\label{Xu}
\mathcal X_u\ =\ a_u\ X_1\ +\ b_u\ X_2\ +\ \left(k_u + \frac{b x_u
- a y_u}{2}\right)\ T\ ,
\end{equation}
\begin{equation}\label{Xv}
\mathcal X_v\ =\ a_v\ X_1\ +\ b_v\ X_2\ +\ \left(k_v + \frac{b x_v
- a y_v}{2}\right)\ T\ .
\end{equation}

From \eqref{dNl}, \eqref{Xu}, \eqref{Xv} and \eqref{wedgesX} one
has
\begin{align}\label{dNl2}
& \frac{d\bN^\la}{d\la}\Bl \\
& = \bigg\{\left(t_v + \frac{y x_v - x y_v}{2}\right) b_u -
\left(t_u + \frac{y x_u - x y_u}{2}\right) b_v + (y_u k_v - y_v
k_u) - \frac{b}{2} (x_u y_v - x_v y_u)\bigg\} X_1
\notag\\
& +\ \bigg\{\left(t_u + \frac{y x_u - x y_u}{2}\right) a_v -
\left(t_v + \frac{y x_v - x y_v}{2}\right) a_u + (x_v k_u - x_u
k_v) + \frac{a}{2} (x_u y_v - x_v y_u)\bigg\} X_2
\notag\\
& +\ \bigg\{(y_v a_u - y_u a_v) + (x_u b_v - x_v b_u)\bigg\} T\ .
\notag
\end{align}

Equations  \eqref{dNl}, \eqref{Xu}, \eqref{Xv} also give
\begin{align}\label{d2Nl}
\frac{d^2\bN^\la}{d\la^2}\Bl & = 2\ \bigg\{\bigg[(b_u \left(k_v +
\frac{b x_v - a y_v}{2}\right) - b_v \left(k_u + \frac{b x_u - a
y_u}{2}\right)\bigg] X_1
\\
& + \bigg[ a_v \left(k_u + \frac{b x_u - a y_u}{2}\right) - a_u
\left(k_v + \frac{b x_v - a y_v}{2}\right)\bigg] X_2 \notag\\
& + (a_u b_v - a_v b_u) T\bigg\}\ . \notag
\end{align}

 We now let
\begin{equation}\label{xl}
\begin{cases}
X^\la_1\ =\ X_1(\theta^\la)\  =\ X_1\ -\ \la\ \frac{b}{2}\ T\ ,
\\
X^\la_2\ =\ X_2(\theta^\la)\ =\ X_2\ +\ \la\ \frac{a}{2}\ T\ ,
\end{cases}
\end{equation}
for which we clearly have
\begin{equation}\label{derxl}
\frac{dX^\la_1}{d\la}\ =\ -\ \frac{b}{2}\ T\ ,\quad\quad
\frac{dX^\la_2}{d\la}\ =\ \frac{a}{2}\ T\ ,\quad\quad
\frac{d^2X^\la_1}{d\la^2}\ =\ 0\ ,\quad\quad
\frac{d^2X^\la_2}{d\la^2}\ =\ 0\ . \end{equation}

Consider the quantities in \eqref{plambda0}. From \eqref{derxl} we
find
\begin{equation}\label{1derpl}
\frac{d p^\la}{d\la}\Bl\ =\ <\frac{d\bN^\la}{d\la}\Bl, X_1>\ -\
\frac{b}{2}\ <\bN,T>\ ,
\end{equation}
\begin{equation}\label{2derpl}
\frac{d^2 p^\la}{d\la^2}\Bl\ =\ <\frac{d^2\bN^\la}{d\la^2}\Bl,
X_1>\ -\  b\  <\frac{d\bN^\la}{d\la}\Bl, T>\ .
\end{equation}

Similarly, we find
\begin{equation}\label{1derql}
\frac{d q^\la}{d\la}\Bl\ =\ <\frac{d\bN^\la}{d\la}\Bl, X_2>\ +\
\frac{a}{2}\ <\bN,T>\ ,
\end{equation}
\begin{equation}\label{2derql}
\frac{d^2 q^\la}{d\la^2}\Bl\ =\ <\frac{d^2\bN^\la}{d\la^2}\Bl,
X_2>\ +\  a\ <\frac{d\bN^\la}{d\la}\Bl, T>\ .
\end{equation}

\medskip

\begin{lemma}\label{L:dls}
Let $p^\la$ and $q^\la$ relative to the surface $\mS^\la =
J_\la(\mS)$, where $J_\la$ is defined by \eqref{Jlambda}, then
\[
\frac{d p^\la}{d\la}\Bl =\ W \left\{- \left(Zb + b\ \ob\right) -
\qb\ \ob\ Zk\ +\ \pb \left(Tk - \ob Yk\right)\right\}\ ,
\]
\[
\frac{d q^\la}{d\la}\Bl\ =\ W \left\{\left(Za + a\ \ob\right) + \pb\
\ob\ Zk\ +\ \qb \left(Tk - \ob Yk\right)\right\}\ .
\]
\end{lemma}

\begin{proof}[\textbf{Proof}]
Using \eqref{dNl2}, \eqref{1derpl} and the third equation in
\eqref{ppar}, we obtain
\begin{equation}\label{dpl}
\frac{d p^\la}{d\la}\Bl \ =\ \left(t_v + \frac{y x_v - x
y_v}{2}\right) b_u - \left(t_u + \frac{y x_u - x y_u}{2}\right) b_v
+ (y_u k_v - y_v k_u)  - b\ \omega\ .
\end{equation}

We now use the equations \eqref{zuzv} to express the derivatives
$b_u, b_v, k_u, k_v$ in terms of derivatives $Zb, Yb,  Tb$ with
respect to the orthonormal frame $\{Z,Y,T\}$. Ordering terms one
finds
\begin{align*}\label{dpl2}
\frac{d p^\la}{d\la}\Bl \ & =\ \left[\left(t_v + \frac{y x_v - x
y_v}{2}\right) (x_u \pb + y_u \qb) - \left(t_u + \frac{y x_u - x
y_u}{2}\right)(x_v \pb + y_v \qb)\right]\ Yb
\\
& +\ \left[\left(t_v + \frac{y x_v - x y_v}{2}\right) (x_u \qb -
y_u \pb) - \left(t_u + \frac{y x_u - x y_u}{2}\right)(x_v \qb -
y_v \pb)\right]\ Zb \notag\\
& +\ \bigg[y_u (x_v \pb + y_v \qb) - y_v (x_u \pb + y_u
\qb)\bigg]\
Yk \notag\\
& +\ \bigg[y_u (x_v \qb - y_v \pb) - y_v (x_u \qb - y_u
\pb)\bigg]\
Zk \notag\\
& +\ \bigg[y_u t_v - y_v t_u - \frac{y}{2} (x_u y_v - x_v
y_u)\bigg]\ Tk\ -\ b\ \omega\ . \notag
\end{align*}

Simplifying in the latter equation, gives
\begin{align*}
\frac{d p^\la}{d\la}\Bl \ & =\ \bigg[- \left(x_v t_u - x_u t_v +
\frac{x}{2} (x_u y_v - x_v y_u)\right) \pb + \left(y_u t_v - y_v
t_u - \frac{y}{2} (x_u y_v - x_v y_u)\right) \qb\bigg]\ Yb
\\
& +\ \bigg[ - \left(x_v t_u - x_u t_v + \frac{x}{2} (x_u y_v - x_v
y_u)\right) \qb - \left(y_u t_v - y_v t_u - \frac{y}{2} (x_u y_v -
x_v y_u)\right) \pb\bigg]\ Zb
\\
& -\ (x_u y_v - x_v y_u) \pb Yk - (x_u y_v - x_v y_u) \qb Zk + \pb W
Tk - b\ \omega\ ,
\end{align*}
where we have used \eqref{ppar}. At this point we notice that, in
view of \eqref{ppar} again, the coefficient of $Yb$ vanishes, and we
obtain from the remaining terms the following espression
\begin{equation*}
\frac{d p^\la}{d\la}\Bl \  =\ W\ \left\{- Zb - b \ \ob\ +\ \pb\ Tk\
-\ \ob (\pb Yk + \qb Zk)\right\}\ ,
\end{equation*}
which gives the first equation in the thesis of the lemma. In a
similar fashion, we obtain the desired expression of $\frac{d
q^\la}{d\la}\Bl$.

\end{proof}

\medskip

From Lemma \ref{L:dls} we immediately obtain the following crucial
result.

\medskip

\begin{lemma}\label{L:pdpqdq}
In the situation of Lemma \ref{L:dls} we have
\[
\left(p^\la \frac{d p^\la}{d\la} + q^\la \frac{d
q^\la}{d\la}\right)\Bl\ =\ W^2\ \left\{\left(Tk - \ob Yk\right) +
(\qb Za - \pb Zb) + (\qb a - \pb b) \ob\right\}\ .
\]
\end{lemma}

\medskip

With Lemma \ref{L:pdpqdq} in hands we are ready to give the proof
of Theorem \ref{T:variations}.

\medskip

\begin{proof}[\textbf{Proof of Theorem \ref{T:variations}}]
Substituting the equation in Lemma \ref{L:pdpqdq} in
\eqref{1vargen} of Lemma \ref{L:1&2var}, we obtain
\begin{align}\label{fvpar2}
\fv\ & =\ \int_\mS \left\{Tk\ -\ \ob\ Yk\right\} d\sigma_H
\\
& +\ \int_\mS \left\{(\qb Za - \pb Zb) + (\qb a - \pb b) \ob\right\}
d\sigma_H\  .\notag
\end{align}

In order to extract the geometry from \eqref{fvpar2} we need to
convert the two integrals in the right-hand side into ones which
involve only the functions $a, b$ and $k$, and not their covariant
derivatives along the orthonormal frame $\{Z, Y, T\}$. This is where
we use Lemma \ref{L:ibpcombined} for the first time. Applying
\eqref{intbyparts1}, \eqref{intbyparts2} we find
\begin{align}\label{ibp12}
& \int_\mathcal S \left\{\qb\ Za\ +\ \qb\ \ob\ a\right\}\ d\sigma_H\
=\ \int_\mathcal S \pb\ a\ \mathcal H\ d\sigma_H\ ,
\\
& \int_\mathcal S \left\{\pb\ Zb\ +\ \pb\ \ob\ b\right\}\ d\sigma_H\
=\ -\ \int_\mathcal S \qb\ b\ \mathcal H\ d\sigma_H\ . \notag
\end{align}

Furthermore, Lemma \ref{L:ibpcombined} gives
\begin{equation}\label{fvpar3}
 \int_\mS \left\{Tk\ -\ \ob\ Yk\right\} d\sigma_H
 \ =\
 \int_\mS \ob\ k\ \mathcal H\ d\sigma_H\
,
\end{equation}

Combining \eqref{ibp12}, \eqref{fvpar3} with \eqref{fvpar2} we
obtain
\begin{equation}\label{fvnonangle}
\fv\ =\ \int_\mS \mathcal H\ \left\{ \pb\ a\ +\ \qb\ b\ + \ob\ k
\right\}\ d\sigma_H\ .
\end{equation}

Recalling the definition \eqref{inner} of $F$, we reach the desired
conclusion.

\end{proof}

\medskip

With Lemma \ref{L:dls} and some elementary computations, we obtain
the following result which is useful in the proof of Theorem
\ref{T:2varfinal} since it provides the integrand of the second
addend in the right-hand side of \eqref{2vargen}.

\medskip

\begin{lemma}\label{L:squares}
In the situation of Lemma \ref{L:dls} we have
\begin{align}\label{squares1}
\frac{\left(\frac{d p^\la}{d\la}\right)^2\Bl\ +\ \left(\frac{d
q^\la}{d\la}\right)^2\Bl}{W^2}\ & =\ \left(Za + \ \pb\ \ob\
Zk\right)^2\ +\ \left(Zb + \ \qb\ \ob\ Zk\right)^2
\\
& +\ \left(Tk - \ob Yk\right)^2\ +\ (a^2 + b^2)\ \ob^2
\notag\\
& +\ 2\ \left(Tk - \ob Yk\right)\left(\qb Za - \pb
Zb + \ob (a \qb - b \pb)\right) \notag\\
& +\ 2\ \ob (a Za + b Zb)\ +\ 2\ \ob^2 (a \pb + b \qb) Zk\ . \notag
\end{align}
\end{lemma}

\medskip

We finally turn to the computation of the first addend in the
right-hand side of \eqref{2vargen}, i.e., $\left( p^\lambda\
\frac{d^2 p^\lambda}{d\lambda^2} + q^\lambda\ \frac{d^2
q^\lambda}{d\lambda^2}\right)\Bl$. From \eqref{2derpl} and from
\eqref{dNl2}, \eqref{d2Nl} we obtain
\begin{align}\label{d2nl2}
& \left( p^\lambda\ \frac{d^2 p^\lambda}{d\lambda^2} + q^\lambda\
\frac{d^2 q^\lambda}{d\lambda^2}\right)\Bl\ =\ W\
\left\{<\frac{d^2N^\la}{d\la^2}\Bl,\nuX> + (a \qb - b \pb)
<\frac{dN^\la}{d\la}\Bl,T>\right\}
\\
& =\ W\ \bigg\{\bigg[2 k_v (\pb b_u - \qb a_u)\ +\ 2 k_u (\qb a_v
- \pb b_v)\bigg]
\notag\\
& +\ \bigg[(b x_v - a y_v) (\pb b_u - \qb a_u) + (b x_u - a y_u)
(\qb
a_v - \pb b_v)\bigg] \notag\\
& +\ (a \qb - b \pb) \bigg[(y_v a_u - y_u a_v) + (x_u b_v - x_v
b_u)\bigg]\bigg\}\ . \notag
\end{align}

Now we compute the three expressions in square brackets in the
right-hand side of \eqref{d2nl2}. Using \eqref{zuzvYZ} we obtain
\begin{align}\label{uffau}
(\pb b_u - \qb a_u)\ & =\  (x_u \pb + y_u \qb) (\pb Yb - \qb Ya)\
+\ (x_u \qb - y_u \pb) (\pb Zb - \qb Za)
\\
& +\  \left(t_u + \frac{y x_u - x y_u}{2}\right) (\pb Tb - \qb
Ta)\ ,\notag
\end{align}

\begin{align}\label{uffav}
(\qb a_v - \pb b_v)\ & =\  -\ (x_v \pb + y_v \qb) (\pb Yb - \qb
Ya)\ -\ (x_v \qb - y_v \pb) (\pb Zb - \qb Za)
\\
& -\  \left(t_v + \frac{y x_v - x y_v}{2}\right) (\pb Tb - \qb
Ta)\ .\notag
\end{align}

These formulas, combined with \eqref{ppar}, give
\begin{align}\label{uffa2}
& (b x_v - a y_v) (\pb b_u - \qb a_u) + (b x_u - a y_u) (\qb a_v -
\pb b_v)
\\
& =\ W\ \bigg\{\pb (a \pb + b \qb) \left(Tb - \ob Yb\right) - \qb (a
\pb + b \qb) \left(Ta - \ob
Ya\right) \notag\\
& +\ (b \pb - a \qb) (\pb Zb - \qb Za) \ob\bigg\}\ . \notag
\end{align}

Again from \eqref{zuzvYZ} and \eqref{ppar}, we obtain
\begin{align}\label{uffa3}
& y_v a_u - y_u a_v\ =\ \bigg[- \pb \left(Ta - \ob Ya\right) + \qb
\ob Za\bigg]\ W\ ,
\end{align}
and
\begin{align}\label{uffa4}
& x_u b_v - x_v b_u\ =\ \bigg[- \qb \left(Tb - \ob Yb\right) - \pb
\ob Zb\bigg]\ W\ .
\end{align}

Formulas \eqref{uffa3}, \eqref{uffa4} give
\begin{align}\label{uffa5}
& (a \qb - b \pb) \bigg[(y_v a_u - y_u a_v) + (x_u b_v - x_v
b_u)\bigg]
\\
& =\ (a \qb - b \pb) \bigg[- \pb \left(Ta - \ob
Ya\right) - \qb \left(Tb - \ob Yb\right) \notag\\
& + (\qb Za - \pb Zb)  \ob \bigg]\ W\ . \notag
\end{align}

Finally,  a (long) computation, based on \eqref{zuzvYZ},
\eqref{ppar}, \eqref{uffau}, \eqref{uffav}, and the identity $\pb^2
+ \qb^2 = 1$, give
\begin{align}\label{uffa6}
& 2 k_v (\pb b_u - \qb a_u) + 2 k_u (\qb a_v - \pb b_v)
\\
& =\ 2\ W\ \bigg\{ - (\pb Zb - \qb Za)  \left(Tk -
\ob Yk\right) \notag\\
& +\ Zk \bigg[\pb \left(Tb - \ob Yb\right) - \qb \left(Ta - \ob
Ya\right)\bigg]\bigg\}\ . \notag
\end{align}

From \eqref{d2nl2}, \eqref{uffa2}, \eqref{uffa5} and \eqref{uffa6}
we finally conclude
\begin{align}\label{uffa7} & \frac{\left( p^\lambda\ \frac{d^2
p^\lambda}{d\lambda^2} + q^\lambda\ \frac{d^2
q^\lambda}{d\lambda^2}\right)\Bl}{W^2}
\\
& =\ -\ 2\ (\pb Zb - \qb Za)  \left(Tk -
\ob Yk\right) \notag\\
& +\ 2\ Zk \bigg[\pb \left(Tb - \ob Yb\right) - \qb
\left(Ta - \ob Ya\right)\bigg] \notag\\
& +\ \pb (a \pb + b \qb) \left(Tb - \ob Yb\right) - \qb (a \pb + b
\qb) \left(Ta - \ob
Ya\right) \notag\\
& +\ (b \pb - a \qb) (\pb Zb - \qb Za) \ob
 \notag\\
 & +\ (a \qb - b \pb) \bigg[- \pb \left(Ta - \ob
Ya\right) - \qb \left(Tb - \ob Yb\right) \notag\\
& + (\qb Za - \pb Zb)  \ob \bigg]\  . \notag
\end{align}

We have thus proved the following lemma.

\medskip

\begin{lemma}\label{L:uffa8}
In the situation of Lemma \ref{L:dls} we have
\begin{align}\label{uffa8} & \frac{\left( p^\lambda\ \frac{d^2
p^\lambda}{d\lambda^2} + q^\lambda\ \frac{d^2
q^\lambda}{d\lambda^2}\right)\Bl}{W^2}
\\
& =\ -\ 2\ (\pb Zb - \qb Za)  \left(Tk -
\ob Yk\right) \notag\\
& +\ \left(Ta - \ob Ya\right) \bigg[- 2 \qb Zk - \qb
(a \pb + b \qb) - \pb (a \qb - b \pb)\bigg] \notag\\
& +\ \left(Tb - \ob Yb\right)\bigg[2 \pb Zk + \pb (a
\pb + b \qb) - \qb (a \qb - b \pb)\bigg] \notag\\
& +\ 2\ (a \qb - b \pb) (\qb Za - \pb Zb)\ob \ . \notag
\end{align}
\end{lemma}

\medskip

We can finally give the proof of Theorem \ref{T:2varfinal}.

\medskip

\begin{proof}[\textbf{Proof of Theorem \ref{T:2varfinal}}]
Combining \eqref{2vargen} in Lemma \ref{L:1&2var} with Lemmas
\ref{L:pdpqdq}, \ref{L:squares} and \ref{L:uffa8}, we obtain the
desired conclusion.

\end{proof}

\vskip 0.6in

\section{\textbf{The stability of $H$-minimal surfaces}}\label{S:stab}

\vskip 0.2in

Unfortunately, in its present form Theorem \ref{T:2varfinal} is not
as useful as one would wish. A completely analogous situation occurs
in the Riemannian case, where one still needs to carefully use
intrinsic integration by parts to extract the geometry, see
\cite{BGG}. The main objective of this section is to give a
geometric meaning to the second variation formula of Theorem
\ref{T:2varfinal}. We stress that for the sake of simplicity, and
because of its relevance in the applications to stability, we state
it for stationary points of the $H$-perimeter functional
($H$-minimal surfaces), but a more general formula containing the
$H$-mean curvature $H$, along with its covariant derivatives, can be
obtained with some additional work if we use the full form of the
geometric identities in section \ref{S:geomid}. We begin with the
relevant definition.

\medskip

\begin{dfn}\label{D:stable}
Given an oriented $C^2$ surface $\mS \subset \HH$, with $\Sigma_\mS
= \varnothing$, we say that $\mS$ is \emph{stable} if it is
stationary (i.e., $H$-minimal), and if
\[
\sv\ \geq \ 0\ ,\quad\quad\quad\text{for every}\quad \mathcal X \in
C^2_0(\mS,\HH) .
\]
If there exists $\mathcal X \neq 0$ such that $\sv < 0$, then we say
that $\mS$ is \emph{unstable}.
\end{dfn}

\medskip

Our main result concerning the stability is contained in the
following theorem.

\medskip

\begin{thrm}\label{T:svgeometric}
Let $\mS\subset \HH$ be $H$-minimal, then
\[
\sv\ =\ \int_\mS \bigg\{|\del F|^2\ +\ (2\mathcal A - \ob^2)
F^2\bigg\} d\sigma_H\ ,
\]
where $F$ is as in \eqref{inner}. As a consequence, $\mS$ is stable
if and only if the following stability inequality of Hardy type
holds on $\mS$
\[
\int_\mS (\ob^2 - 2\mathcal A) F^2 d\sigma_H\ \leq\ \int_\mS |\del
F|^2\  d\sigma_H\ .
\]
\end{thrm}

\medskip

\begin{cor}\label{C:stabilityvc}
Every vertical plane $\mS = \{(x,y,t)\in \HH\mid \alpha x + \beta y
= \gamma\}$, with $\alpha^2 + \beta^2 \not= 0$, is stable.
\end{cor}

\begin{proof}[\textbf{Proof}]
Consider the defining function $\phi(x,y,t) = \alpha x + \beta y -
\gamma$. One has $\omega = T\phi \equiv 0$, and therefore $\omega =
\mathcal A \equiv 0$. Since every plane in $\HH$ is $H$-minimal, we
can apply Theorem \ref{T:svgeometric}, to find for every vector
field $\mathcal X = a X_1 + b X_2 + k T \in C^2_0(\mS,\HH)$
\[
\sv\ =\ \int_\mS |\del F|^2\ d\sigma_H\ \geq \ 0\ .
\]

This proves the stability of $\mS$. We note explicitly that in the
present situation
\[
F\ =\ \frac{\alpha a + \beta b}{\sqrt{\alpha^2 + \beta^2}}\ .
\]

\end{proof}

\medskip

Another interesting consequence of Theorem \ref{T:svgeometric} is
the following stability inequality for intrinsic graphs. We recall
that a $C^2$ surface $\mS\subset \HH$ is called an intrinsic
$X_1$-graph according to \cite{FSS4} provided that there exist an
open set $\Om \subset \R^2_{(u,v)}$ and a function $\phi\in
C^2(\Om)$ such that $\mS$ can be described  by $(x,y,t) = (0,u,v)
\circ \phi(u,v) e_1 = (0,u,v) \circ (\phi(u,v),0,0)$. This means
that $\mS$ admits the parametrization
\[ \theta(u,v)\ =\
\left(\phi(u,v),u,v - \frac{u}{2} \phi(u,v)\right)\ , \quad\quad
(u,v)\in \Om\ .\]

If instead $\mS$ can be parametrized by
\[ \theta(u,v)\ =\
\left(u, \phi(u,v),v + \frac{u}{2} \phi(u,v)\right)\ , \quad\quad
(u,v)\in \Om\ ,\] then we say that $\mS$ is an intrinsic
$X_2$-graph. We only discuss the case of an intrinsic $X_1$-graph,
leaving to the reader to provide the trivial changes necessary to
treat the case of $X_2$-graphs.

Given a function $F$ denote by $\mathcal B_\phi(F) = F_u + \phi F_v$
the linear transport equation, so that $\mathcal B_\phi(\phi) =
\phi_u + \phi \phi_v$ indicates the nonlinear inviscid Burger
operator acting on $\phi$. Since from \eqref{ppar} we obtain
\begin{equation}\label{pqig}
p\ =\ 1\ ,\quad\quad q\ =\ - \mathcal B_\phi(\phi)\ ,\quad\quad
\omega\ =\ - \phi_v\ ,\quad\quad W\ =\ \sqrt{1 + \mathcal
B_\phi(\phi)^2}\ ,
\end{equation}
we see that the Riemannian normal to an intrinsic $X_1$-graph is
given by
\begin{equation}\label{Nig2}
\bN\ =\ X_1\ - \mathcal B_\phi(\phi) X_2 - \phi_v T\ .
\end{equation}

As a consequence of the first equality in \eqref{pqig} we deduce
that an intrinsic $X_1$-graph always has empty characteristic locus.
Furthermore, again from \eqref{pqig}, and from \eqref{Nig2}, we see
that if $\Om$ is bounded then the $H$-perimeter of $\mS$ is
expressed by the functional
\[ \sigma_H(\mS)\ =\ \mathcal P(\phi)\ =\ \int_\Om \sqrt{1 + \mathcal B_\phi(\phi)^2}\ du dv\
,
\]
so that
\[
d\sigma_H\ =\ \sqrt{1 + \mathcal B_\phi(\phi)^2}\ du dv\ ,
\]
see also \cite{ASV}.

\medskip

\begin{cor}\label{C:svig}
Let $\mS$ be a $C^2$ $H$-minimal, intrinsic $X_1$-graph, then $\mS$
is stable if and only if
\[
\int_\Om \frac{\phi_v^2 + 2 \mathcal B_\phi(\phi_v)}{\sqrt{1 +
\mathcal B_\phi(\phi)^2}}\ F^2\ du dv\ \leq\ \int_\Om \frac{\mathcal
B_\phi(F)^2}{\sqrt{1 + \mathcal B_\phi(\phi)^2}}\ du dv\ ,
\]
where $F$ is as in \eqref{inner}.
\end{cor}

\begin{proof}[\textbf{Proof}]
We begin by observing that, thanks to \eqref{pqig} we have
\begin{equation}\label{YZig}
\begin{cases}
Y\ =\ \nuX\ =\ \frac{1}{\sqrt{1 + \mathcal B_\phi(\phi)^2}} X_1\ -\
\frac{\mathcal B_\phi(\phi)}{\sqrt{1 + \mathcal B_\phi(\phi)^2}}
X_2\ ,
\\
Z\ =\ \nup\ =\ -\ \frac{\mathcal B_\phi(\phi)}{\sqrt{1 + \mathcal
B_\phi(\phi)^2}} X_1\ -\ \frac{1}{\sqrt{1 + \mathcal
B_\phi(\phi)^2}} X_2\ .
\end{cases}
\end{equation}

Given a function $f$ on $\mS$, by abuse of notation we continue to
indicate with the same letter the function $f(u,v) =
f(\theta(u,v))$. A simple use of the chain rule as in \eqref{crh}
gives
\[
f_u\ =\ \phi_u\ X_1 f\ +\ X_2 f\ -\ \phi\ Tf\ ,\quad\quad f_v\ =\
\phi_v\ X_1 f\ +\ Tf\ ,
\]
where in the right-hand sides of the latter equations we have
written $X_1f$ for $X_1f\circ \theta$, and similarly for $X_2f, Tf$.
Using the latter two equations and the second equation in
\eqref{YZig}, we obtain
\begin{equation}\label{Zig}
Zf\ =\ -\ \frac{\mathcal B_\phi(f)}{\sqrt{1 + \mathcal
B_\phi(\phi)^2}}\ .
\end{equation}

We now use \eqref{Zig} to compute $\mathcal A = - Z\ob$. From the
latter two equations in \eqref{pqig}, one has
\begin{align*}
\mathcal A\ & =\ Z\left(\frac{\phi_v}{\sqrt{1 + \mathcal
B_\phi(\phi)^2}}\right)
\\
& = \frac{Z(\phi_v)}{\sqrt{1 + \mathcal B_\phi(\phi)^2}}\ -\ \phi_v\
\frac{Z(\sqrt{1 + \mathcal B_\phi(\phi)^2})}{1 + \mathcal
B_\phi(\phi)^2}
\\
  & =\ -\ \frac{\mathcal B_\phi(\phi_v)}{1 + \mathcal
  B_\phi(\phi)^2}\ +\ \phi_v\ \frac{\mathcal B_\phi(\phi) \mathcal
  B_\phi(\mathcal B_\phi(\phi))}{(1 + \mathcal B_\phi(\phi)^2)^2}\ .
\end{align*}

One can now recognize that the $H$-mean curvature of $\mS$ is given
by
\begin{equation}\label{MCig}
 \mathcal B_\phi\left(\frac{\mathcal
B_\phi(\phi)}{\sqrt{1 + \mathcal B_\phi(\phi)^2}}\right)\ =\ -\
\mathcal H\ , \end{equation} see \cite{GS}, and also \cite{BSV}.
Using \eqref{MCig}, after some simple computations, we obtain that
the condition that $\mS$ be $H$-minimal is expressed by
\[
\mathcal B_\phi(\mathcal B_\phi(\phi))\ =\ 0\ .
\]

Substituting this equation in the above formula for $\mathcal A$ we
conclude that
\[
\mathcal A\ =\ -\ \frac{\mathcal B_\phi(\phi_v)}{1 + \mathcal
B_\phi(\phi)^2}\ .
\]

Again from \eqref{pqig} we finally obtain
\[
\ob^2\ -\ 2\ \mathcal A\ =\ \frac{\phi_v^2 + 2 \mathcal
B_\phi(\phi_v)}{1 + \mathcal B_\phi(\phi)^2}\ .
\]

To reach the desired conclusion we are left with using the latter
equation in the stability inequality in Theorem \ref{T:svgeometric},
in combination with the expression of $d\sigma_H$ and with
\eqref{Zig}.

\end{proof}

\medskip

In \cite{DGN3} it was conjectured that the only $C^2$ stable
intrinsic graphs in $\HH$ are the vertical planes. Using also the
results in \cite{DGN3}, in \cite{BSV} the authors have provided a
positive answer to this conjecture. We next turn to the proof of
Theorem \ref{T:svgeometric}.

\medskip

\begin{proof}[\textbf{Proof of Theorem \ref{T:svgeometric}}]
Since we want to extract a more geometrically meaningful formula
from the general expression in Theorem \ref{T:2varfinal}, we will
now make several reductions.  First, expanding the three squares,
and regrouping terms using repeatedly $\pb^2 + \qb^2 = 1$, we find
for the integrand in the right-hand side of \eqref{uffa9}
\begin{align*}
\mathcal Integrand\ & =\ 2\ (\qb Za - \pb Zb)  \left(Tk - \ob
Yk\right)
\\
& +\ \left(Ta - \ob Ya\right) \bigg[- 2 \qb Zk - \qb
(a \pb + b \qb) - \pb (a \qb - b \pb)\bigg] \notag\\
& +\ \left(Tb - \ob Yb\right)\bigg[2 \pb Zk + \pb (a
\pb + b \qb) - \qb (a \qb - b \pb)\bigg] \notag\\
& +\ (\pb Za + \qb Zb)^2\ +\ \ob^2
(Zk)^2\ +\ 2(\pb Za + \qb Zb)\ob Zk \notag\\
&  +\ 2 \ob^2 (a \pb + b \qb) Zk\ +\ 2
\ob(a Za + b Zb) \notag\\
& +\ (a^2 \pb^2 + b^2 \qb^2 + 2 \pb\ \qb\ a b) \ob^2 \notag
\end{align*}

We easily obtain from the latter equation
\begin{align}\label{uffa10}
\sv\  & =\ 2\ \int_\mS (\qb Za - \pb Zb) \left(Tk - \ob Yk\right)\
d\sigma_H
\\
& +\ 2\ \int_\mS  \bigg[ - \qb  \left(Ta - \ob Ya\right) + \pb
\left(Tb - \ob Yb\right)\bigg] Zk\ d\sigma_H \notag\\ & +\ \int_\mS
\bigg\{\left(Ta - \ob Ya\right) \bigg[ - \qb
(a \pb + b \qb) - \pb (a \qb - b \pb)\bigg] \notag\\
& +\ \left(Tb - \ob Yb\right)\bigg[\pb (a
\pb + b \qb) - \qb (a \qb - b \pb)\bigg] \notag\\
& +\ \left(\pb Za + \qb Zb +  \ob Zk\right)^2
\notag\\
&  +\ 2 \ob^2 (a \pb + b \qb) Zk\ +\ 2
\ob (a Za + b Zb) \notag\\
& +\ (a \pb + b \qb)^2 \ \ob^2\bigg\}\ d\sigma_H\ . \notag
\end{align}

Our final objective is to remove all derivatives from the functions
$a, b$ and $k$ from the right-hand side of \eqref{uffa10}. This is
somewhat delicate and involves some effort. The final product will
be achieved by a repeated use of the basic integration by parts
Lemma \ref{L:ibpcombined} and of the geometric identities in section
\ref{S:geomid}. We begin by observing that, thanks to \eqref{Zp&Zq},
the $H$-minimality of $\mS$ implies that
\[
\qb Za - \pb Zb\ =\ Z(\qb a - \pb b)\ .
\]

Using this observation, Lemma \ref{L:ibpcombined}, and Corollary
\ref{C:mixedcomm}, we obtain
\begin{align}\label{monsterone}
& \int_\mS \bigg\{ 2\ (\qb Za - \pb Zb)  \left(Tk - \ob Yk\right)\
d\sigma_H
\\
& =\ -\ 2\ \int_\mS  (\qb a - \pb b) Z(T - \ob Y)k \ d\sigma_H \ -\
2\ \int_\mS \ob (\qb a - \pb b) (T - \ob Y)k\
d\sigma_H\notag\\
& =\ -\ 2\ \int_\mS (\qb a - \pb b) (T - \ob Y)Zk\ d\sigma_H \ +\ 2\
\int_\mS (\qb a - \pb b) [T - \ob Y, Z]k\ d\sigma_H
\notag\\
& -\ 2\ \int_\mS \ob (\qb a - \pb b)(Tk - \ob Yk)\ d\sigma_H
\notag\\
& =\ 2\
\int_\mS Z k\ (T - \ob Y)(\qb a - \pb b)\ d\sigma_H \notag\\
& -\ 2\ \int_\mS \ob k\ (T - \ob Y)(\qb a - \pb b)\ d\sigma_H\ +\ -\ 2\ \int_\mS \ob k\ (T - \ob Y)(\qb a - \pb b)\ d\sigma_H\notag\\
& =\ 2\ \int_\mS \bigg[\qb (Ta - \ob Ya) - \pb
(Tb - \ob Yb)\bigg] Zk \ d\sigma \notag\\
& +\ 2\ \int_\mS a Zk\ (T\qb - \ob Y\qb)\ d\sigma\ -\ 2\ \int_\mS b
Zk\ (T\pb - \ob Y\pb)\ d\sigma\ . \notag
\end{align}

Substituting \eqref{monsterone} in \eqref{uffa10}, we find
\begin{align}\label{uffafa10}
\sv\  & =\ 2\ \int_\mS a Zk (T\qb - \ob Y\qb)\ d\sigma\ -\ 2\
\int_\mS b Zk (T\pb - \ob Y\pb)\ d\sigma
\\
& +\ \left(Ta - \ob Ya\right) \bigg[- \qb
(a \pb + b \qb) - \pb (a \qb - b \pb)\bigg] \notag\\
& +\ \left(Tb - \ob Yb\right)\bigg[\pb (a
\pb + b \qb) - \qb (a \qb - b \pb)\bigg] \notag\\
& +\ \left(\pb Za + \qb Zb +  \ob Zk\right)^2
\notag\\
&  +\ 2 \ob^2 (a \pb + b \qb) Zk\ +\ 2
\ob (a Za + b Zb) \notag\\
& +\ (a \pb + b \qb)^2 \ \ob^2\bigg\}\ d\sigma_H\ . \notag
\end{align}

It should be clear to the reader that we have made some interesting
progress, since we have eliminated terms containing products of
derivatives of $a, b$ and $k$. However, we are still far from our
final goal. We next use the $H$-minimality of $\mS$, and
\eqref{Zp&Zq} again, to see that
\begin{align*}
(ZF)^2\ & =\  \bigg(\left(\pb Za + \qb Zb +  \ob Zk\right) -
\mathcal A k\bigg)^2
\\
& =\ \bigg(\pb Za + \qb Zb + \ob Zk\bigg)^2\ +\ \mathcal A^2 k^2\ -\
2 \mathcal A k (\pb Za + \qb Zb)\ -\ \ob \mathcal A Z(k^2)\ ,
\end{align*}
where $F$ is as in \eqref{inner} and $\mathcal A$ is the function
defined in Corollary \ref{C:Zob}. This identity, the fact that $\pb
Za + \qb Zb = Z(a\pb + b \qb)$ (the $H$-minimality of $\mS$),
Corollary \ref{C:za} and Lemma \ref{L:ibpcombined}, give
\begin{align}\label{Fsquare}
& \int_\mS \bigg(\pb Za + \qb Zb + \ob Zk\bigg)^2\ d\sigma_H \ =\
\int_\mS (ZF)^2\ d\sigma_H\
\\
& +\ 2\ \int_\mS \mathcal A k Z(a \pb + b \qb)\ d\sigma_H\ -\
\int_\mS (\ob^2 - 2 \mathcal A) \ob^2 k^2\ d\sigma_H
\notag\\
& =\ \int_\mS (ZF)^2\ d\sigma_H\ +\ \int_\mS (2 \mathcal A - \ob^2)
(\ob^2 k^2 + 2 \ob k (a\pb + b \qb))\ d\sigma_H \notag\\
& -\ 2\ \int_\mS \mathcal A (a \pb + b \qb) Zk\ d\sigma_H\ .\notag
\end{align}

Substitution of \eqref{Fsquare} into \eqref{uffafa10} gives
\begin{align}\label{uffafa1}
\sv\  & =\  \int_\mS (ZF)^2\ d\sigma_H\ +\ \int_\mS (2 \mathcal A -
\ob^2) (\ob^2 k^2 + 2 \ob k (a\pb + b \qb))\ d\sigma_H
\\
& +\ 2\ \int_\mS a Zk\bigg[(T\qb - \ob Y\qb) + \pb\ \ob^2 - \pb
\mathcal A\bigg] \ d\sigma_H \notag\\
& \ +\ 2\ \int_\mS b Zk\bigg[- (T\pb - \ob Y\pb) + \qb\ \ob^2 - \qb
\mathcal A\bigg] \ d\sigma_H
\notag\\
& +\ \int_\mS \bigg\{\left(Ta - \ob Ya\right) \bigg[- \qb
(a \pb + b \qb) - \pb (a \qb - b \pb)\bigg] \notag\\
& +\ \left(Tb - \ob Yb\right)\bigg[\pb (a
\pb + b \qb) - \qb (a \qb - b \pb)\bigg]\bigg\} d\sigma_H \notag\\
& +\  \int_\mS \ob Z(a^2 + b^2) \ d\sigma_H\ +\ \int_\mS (a \pb + b
\qb)^2 \ \ob^2\ d\sigma_H\ . \notag
\end{align}

We now claim \[ (T\qb - \ob Y\qb) + \pb\ \ob^2 - \pb \mathcal A \ =
\ - \ (T\pb - \ob Y\pb) + \qb\ \ob^2 - \qb \mathcal A\ =\ 0\ .
\]

We only check the first of the two equations, leaving it to the
reader to provide the details of the second one. Corollary
\ref{C:Zob} and the identity $\pb^2 + \qb^2 = 1$ give
\begin{align*}
\pb \mathcal A\ & =\ \pb^2(T\qb - \ob Y\qb) - \pb\ \qb (T\pb - \ob
Y\pb) + \pb\ \ob^2\\
& =\ (T\qb - \ob Y\qb) - \qb^2(T\qb - \ob Y\qb) - \pb\ \qb (T\pb -
\ob
Y\pb) + \pb\ \ob^2\\
& =\ (T\qb - \ob Y\qb) + \pb\ \ob^2 - \qb\bigg[\qb (T\qb - \ob Y\qb)
+ \pb (T\pb - \ob Y\pb)\bigg]
\\
& =\ (T\qb - \ob Y\qb) + \pb\ \ob^2 \ ,
\end{align*}
which proves the first identity. Using the claim in \eqref{uffafa1},
with a further application of Lemma \ref{L:ibpcombined}, we obtain
\begin{align}\label{uffafa2}
\sv\  & =\  \int_\mS (ZF)^2\ d\sigma_H\ +\ \int_\mS (2 \mathcal A -
\ob^2) (\ob^2 k^2 + 2 \ob k (a\pb + b \qb))\ d\sigma_H
\\
& +\ \int_\mS \bigg\{\left(Ta - \ob Ya\right) \bigg[- \qb
(a \pb + b \qb) - \pb (a \qb - b \pb)\bigg] \notag\\
& +\ \left(Tb - \ob Yb\right)\bigg[\pb (a
\pb + b \qb) - \qb (a \qb - b \pb)\bigg]\bigg\} d\sigma_H \notag\\
& +\  \int_\mS (\mathcal A - \ob^2) (a^2 + b^2) \ d\sigma_H\ +\
\int_\mS (a \pb + b \qb)^2 \ \ob^2\ d\sigma_H\ . \notag
\end{align}

Completing the square in the second integral in the right-hand side
of \eqref{uffafa2}, we find
\begin{align}\label{uffafa3}
\sv\  & =\  \int_\mS (ZF)^2\ d\sigma_H\ +\ \int_\mS (2 \mathcal A -
\ob^2) F^2\ d\sigma_H
\\
& +\ \int_\mS \bigg\{\left(Ta - \ob Ya\right) \bigg[- 2 a \pb\ \qb
 + b (\pb^2 - \qb^2)\bigg] \notag\\
& +\ \left(Tb - \ob Yb\right)\bigg[2 b \pb\ \qb + a (\pb^2 - \qb^2)\bigg]\bigg\} d\sigma_H \notag\\
& +\  \int_\mS (\mathcal A - \ob^2) (a^2 + b^2) \ d\sigma_H\ -\ 2\
\int_\mS (\mathcal A - \ob^2)(a \pb + b \qb)^2 \ d\sigma_H\ . \notag
\end{align}

The proof of the theorem will be completed if we can establish the
following crucial claim
\begin{align}\label{claim}
& \int_\mS \bigg\{\left(Ta - \ob Ya\right) \bigg[- 2 a \pb\ \qb
 + b (\pb^2 - \qb^2)\bigg] \\
& +\ \left(Tb - \ob Yb\right)\bigg[2 b \pb\ \qb + a (\pb^2 - \qb^2)\bigg]\bigg\} d\sigma_H \notag\\
& +\  \int_\mS (\mathcal A - \ob^2) (a^2 + b^2) \ d\sigma_H\ -\ 2\
\int_\mS (\mathcal A - \ob^2)(a \pb + b \qb)^2 \ d\sigma_H\ =\ 0\ .
\notag
\end{align}

To prove \eqref{claim} we proceed as follows. First, Lemma
\ref{L:ibpcombined} gives \begin{align*}
 & -\ 2\
\int_\mS a \pb\ \qb (Ta - \ob Ya)\ d\sigma_H\ =\ \int_\mS a^2 (T -
\ob Y)(\pb\ \qb)\ d\sigma_H \ .\notag
\end{align*}

Therefore, the coefficient of $a^2$ in the left-hand side of
\eqref{claim} is given by
\[
(T - \ob Y)(\pb\ \qb) + (\mathcal A - \ob^2) - 2 (\mathcal A -
\ob^2) \pb^2\ =\ (T - \ob Y)(\pb\ \qb) + (\mathcal A - \ob^2)(\qb^2
- \pb^2)\ =\ 0\ ,
\]
where we have used the identity $\pb^2 + \qb^2 = 1$. Similarly, we
have
\begin{align*}
 &  2\
\int_\mS b \pb\ \qb (Tb - \ob Yb)\ d\sigma_H\ =\ -\ \int_\mS b^2 (T
- \ob Y)(\pb\ \qb)\ d\sigma_H\  ,\notag
\end{align*}
hence the coefficient of $b^2$ in the left-hand side of
\eqref{claim} is given by
\[
- (T - \ob Y)(\pb\ \qb) + (\mathcal A - \ob^2) - 2 (\mathcal A -
\ob^2) \qb^2\ =\ -\ \bigg[(T - \ob Y)(\pb\ \qb) + (\mathcal A -
\ob^2)(\qb^2 - \pb^2)\bigg]\ =\ 0\ .
\]

Finally, we have
\begin{align*}
 &
\int_\mS \bigg\{b (\pb^2 - \qb^2) (Ta - \ob Ya) + a (\pb^2 - \qb^2)
(Tb - \ob Yb)\bigg\} d\sigma_H\ =\ -\ \int_\mS a b (T - \ob Y)(\pb^2
- \qb^2)\ d\sigma_H\ .\notag
\end{align*}

We thus see that the coefficient of $ab$ in the left-hand side of
\eqref{claim} is given by
\[
4 (\ob^2 - \mathcal A) \pb\ \qb + (T - \ob Y)(\qb^2 - \pb^2)\ =\ 0\
.
\]

From these considerations, the claim \eqref{claim} follows. We have
thus completed the proof.

\end{proof}

\end{document}